\documentclass[11pt]{article}

\usepackage{graphicx}
\usepackage{amsmath}
\usepackage{amssymb}
\usepackage{latexsym}
\usepackage{subfigure}
\usepackage{crop}
\usepackage[colorlinks = true, pdfstartview = FitV, linkcolor = blue, citecolor = blue, urlcolor = blue]{hyperref}
\usepackage{float}

\definecolor{turquoise}{cmyk}{0.65,0,0.1,0.1}
\definecolor{purple}{rgb}{0.65,0,0.65}
\definecolor{green}{rgb}{0, 0.5, 0}
\definecolor{blue}{rgb}{0, 0, 1}
\definecolor{orange}{rgb}{0.8, 0.6, 0.2}
\definecolor{red}{rgb}{0.8, 0.2, 0.2}
\definecolor{brown}{rgb}{0.5, 0.16, 0.16}

\begin{document}

\title{Data driven recovery of local volatility surfaces}

\author{Vinicius Albani\thanks{Computational Science Center, University of Vienna, 
1090 Vienna, Austria, \href{mailto:vvla@impa.br}{\tt vvla@impa.br}} \,
	Uri M. Ascher\thanks{Dept. of Computer Science, University of British Columbia, Canada, \href{mailto:ascher@cs.ubc.ca}{\tt ascher@cs.ubc.ca}} \,
Xu Yang\thanks{IMPA, 
Rio de Janeiro, Brazil, \href{mailto:xuyang@impa.br}{\tt xuyang@impa.br}} \, 
and Jorge P. Zubelli\thanks{IMPA, 
Rio de Janeiro, Brazil, \href{mailto:zubelli@impa.br}{\tt zubelli@impa.br}}}

\maketitle

\begin{abstract}
This paper examines issues of data completion and location uncertainty, popular in many practical 
PDE-based
inverse problems, in the context of option calibration via recovery of local volatility surfaces.
While real data is usually more accessible for this application than for many others, the data 
is often given only at
a restricted set of locations. We show that attempts to ``complete missing data'' 
by approximation or interpolation,
proposed in the literature, may produce results that are inferior to treating the data as scarce.
Furthermore, model uncertainties may arise which translate to uncertainty in data locations, 
and we show how a model-based
adjustment of the asset price may prove advantageous in such situations. 
We further compare a carefully calibrated Tikhonov-type regularization approach 
against a similarly adapted EnKF method,
in an attempt to fine-tune the data assimilation process.
The EnKF method offers reassurance as a different method
for assessing the solution in a problem where information about the true solution is
difficult to come by.
However, additional advantage in the latter 
approach turns out to be limited in our context.


\end{abstract}


\section{Data manipulation and local volatility surfaces}
\label{sec:volatility}

The basic setup of data assimilation and inverse problems for model calibration
consists of an assimilation of dynamics, defined for instance as a discretized PDE, and observed data~\cite{reco}. 
The celebrated Kalman filter, for example, does a forward pass on a weighted least squares problem
fitting both model dynamics and data, and it guarantees variance minimization for a linear problem
with Gaussian noise.

However, in practice, the ``statistical sanctity'' of the data is often violated before the assimilation
process commences. This can happen for various reasons and in different circumstances:
\begin{enumerate}
\item
In cases where the data is {\em scarce}, in the sense that it is observed only at a small set
of locations compared to the size of a reasonable discretization mesh of a physical domain,
there would be many models (solutions to the inverse problem) that 
explain the data (e.g.,~\cite{haasol}).
It is then tempting to ``complete'' the data by some interpolation or other approximation (e.g.,~\cite{kahale}),
whereupon the role of an ensuing regularization as a prior is less crucial.
\item 
There may be a hidden uncertainty in the locations of data, 
not only in data values (e.g.,~\cite{huas,grha}).
For instance, engineers often prefer to see data given at regular mesh nodes, 
so a quiet constant interpolation, moving data items to the nearest cell vertex, is common practice.
\item
Data completion may be necessary to obtain a more efficient algorithm~\cite{rodoas1,ksaamrh}. 
\item
A quiet data completion is often assumed by mathematicians in order to enable building theory
for inverse problems. This includes assumptions of available data on continuous boundary segments~\cite{ehn1}, 
or of observed (measured) relationships between unknown functions that are presumed to hold everywhere in a physical domain.
\item
There are situations where some form of data completion and other manipulation is necessary 
because no one knows how to solve the problem otherwise~\cite{ksaamrh}.
\end{enumerate}

These observations raise the following questions: 
(i) when (and in what sense) is it practically acceptable to perform such data manipulations?
(ii) in which circumstances can one gain advantage by treating the observed data more carefully? and 
(iii) how should one assess correctness
of a solution that has been obtained with such manipulated data?
Our general observation is that researchers occasionally, but not always, seem to get away with such ``crimes'', 
in the sense of producing agreeable results.
For instance, in~\cite{rodoas1} the authors obtained agreeable reconstructions so long as the percentage of
completed data did not exceed about 50\%, but not more.
Such empirical evidence is relatively rare in the literature, however, and it depends on the problem at hand.
More insight is therefore required, and such may be gained by considering applied case studies.

In this article we focus on
a model calibration problem 
that has had tremendous impact 
in mathematical finance. 
It concerns the determination of the so-called local volatility surface,
making use  of derivative prices. A good model for the volatility is crucial 
for many applications ranging from risk management to hedging and pricing of financial instruments.
This setting features both scarce data and uncertainty regarding data location,
and it allows us to work with real data, often available through the internet.


The classical Black-Scholes-Merton model had subsumed a constant volatility model $\sigma$ \cite{blsc} 
in a simplifying stochastic dynamics for the underlying process~\cite{kornbook}.
However, the constant volatility assumption was quickly contradicted by
the actual derivative prices observed in the market. The disagreement between the Black-Scholes model-implied prices for different expiration dates and negotiated strike prices became known as the {\em smile} effect. 
A number of practical parametric as well as nonparametric
models have been proposed in this context; see~\cite{volguide} and references therein. The parametric ones try to fulfill different phenomenological features of the observed prices. 
Yet, in a ground breaking paper Dupire~\cite{dupire}
proposed the use of a 
function $\sigma$ that depends on time and the {price at that time}.
For the case of the European call contracts, he replaced the Black-Scholes equation by a PDE
of the form
\begin{eqnarray}
\frac{\partial C}{\partial \tau} &=&  \frac 12 \sigma^2(\tau,K)K^2 \frac{\partial^2 C}{\partial K^2} - 
b K \frac{\partial C}{\partial K} , \quad \tau > 0, K \geq 0 , \label{x1} 
\end{eqnarray}
with initial and boundary conditions (for calls) given by 
\begin{eqnarray}
C(\tau = 0, K) &=& (S_0 - K)^+ , \label{x1bc}\\
\lim_{K \rightarrow \infty} C(\tau,K) &=& 0, \nonumber \\
\lim_{K \rightarrow 0} C(\tau,K) &=& S_0, \nonumber
\end{eqnarray}
where $\tau$ is time to maturity, $K$ is the strike price, and
$C = C(\tau,K)$ is the value of the European call option with expiration date $T=\tau$.
The parameter $S_0$ is the asset's price at a given date.
%
One ensuing complication in using \eqref{x1}, however, is in the {\em calibration} of this model,
by which we mean finding a plausible volatility surface $\sigma(\tau,K)$ that matches, or explains, given market data
on call option values.
The task here is significantly more challenging than in the case where $\sigma$ is a constant. 
This paper deals with the computational challenges that this {\em inverse}, or inference problem
gives rise to.

The {\em forward} problem involves finding
the values of $C$ satisfying the 
differential problem~\eqref{x1}-\eqref{x1bc}
for given $\sigma(\tau,K)$ and $S_0$, evaluated at the points $(\tau,K)$ where data values are available.
A major difficulty here is that the data are {\em scarce}. To explain what we mean by
this, suppose we have discretized the PDE
using, say, the Crank-Nicolson method (cf.~\eqref{x30} below)
on a rectangular mesh that is reasonable in the sense that the essence of the differential
solution is retained in the discrete solution. 
Then the data is scarce in that
the number $M$ of degrees of freedom in the discrete $C$ typically far
exceeds the number of given data values $l$:  $l \ll M$. 
Moreover, the available data in some typical situations are given
at locations that are far from the boundaries of the (truncated) domain on which the
approximated PDE problem is solved.
See Figures~\ref{fig_pbrdata} and~\ref{fig_spxdata} below for 
examples of such (real) data sets.

Now, if the local volatility surface is discretized, or injected, on the same mesh as that
of the forward problem, then
there are roughly $M$ degrees of freedom in $\sigma$, 
which is again potentially far larger than the number of data constraints.
We can of course discretize $\sigma$ on a coarser sub-mesh 
(which in the extreme case would have only one
point, thus leading back to a constant volatility), or parameterize the surface 
in a more involved manner; 
see~\cite{HKPS-2007, Hof-Kra-2005,achpir,eggeng,boylethangaraj2000,achdoupironneau2005,DCJPZ2013,AlbaniZubelli2014} and references therein for further detail.
Here, however, we stick to a straightforward nodal representation of this surface on the full
$C$-mesh in the hope of retaining flexibility and detail, while avoiding artifacts that may arise from
restrictive simplifying assumptions. 
This approach has worked well in geophysical exploration problems~\cite{haasol}, among others.

Thus, the problem of finding a volatility surface that explains the data is often significantly
under-constrained in practice. This does not make it easy to solve, 
however, as the ultimate goal is to
obtain plausible volatility surfaces that can be worked with, and not just to match data.
Our task is therefore to assimilate the data information with the information contained in the
PDE model~\eqref{x1}, using any plausible a priori information as a {\em prior} in the assimilation process.
Such a priori information can vary significantly, addressing concerns of adherence to the financial
model, relative smoothness of the volatility surface, and numerical stability issues, among others.

One approach that has been relatively popular in financial circles is to apply to the data
special interpolation/extrapolation methods that take into account the a priori information
of the financial model (e.g.,~\cite{kahale}). This is used to obtain data values at all points of
the rectangular mesh on which $C$ is defined, and subsequently the ``new data'' is assimilated
with the information that the solutions of the Dupire equation for different $\sigma$'s
yield to calibrate~\eqref{x1}.
An advantage with this {\em data completion} approach is that the data is no longer scarce when we get to the 
redefined inverse problem. This allows for developing some existence and uniqueness theory as well;
see~\cite{acpaper}
and references therein.
However, a disadvantage is that such data completion constitutes a ``statistical crime'', as the errors in
the new data may no longer be considered as independent random variables, see~\cite{rodoas1}. In fact, we get two solutions $C$ that are 
in a sense competing rather than completing one another,
since the one satisfying~\eqref{x1}, even for the ``best'' $\sigma$, does not necessarily satisfy the data
interpolation conditions and vice versa.

In Sections~\ref{sec:scarce} and~\ref{sec:test}
we therefore examine the performance of this data completion approach 
against that of a scarce data approach that is based on a carefully tuned Tikhonov-type
regularization. Using both synthetic and real data sets, we show that the scarce data
approach can give better and more reliable results; in our reported experiments this has happened especially
for the real data applications.

We then continue with the scarce data approach.

The maximum a posteriori (MAP) functional considered in Section~\ref{sec:scarce} is based
on the statistical assumption that the data error covariance matrix is a scalar multiple
of the identity. In Section~\ref{sec:enkf} we subsequently consider an algorithm,
based on an approach recently proposed in~\cite{iglast} and~\cite{johns2008two}, 
where we attempt to learn more about
the error covariance matrix as the iterative process progresses, using ensemble Kalman filter (EnKF)
techniques. Although our problem is time-dependent, the time variable here does not really differ
from the other independent variable in the usual sense. 
In particular, the unknown surface $\sigma$
depends on both $K$ and $\tau$, 
unlike for instance a material function in reservoir simulation~\cite{iglast} or electromagnetic data
inversion~\cite{haasol}, which are independent of time.
Thus, the EnKF-like methods considered 
use an {\em artificial time}~\cite{ashudo}.
%
In Section~\ref{sec:enkf} we find that the EnKF algorithm can be
improved in our context by adding smoothing prior penalties, just like in Section~\ref{sec:scarce}.
The probabilistic setup, although general, is not fully effective as a substitute for
prior knowledge that is available in no uncertain terms.

The problem setting used in Sections~\ref{sec:scarce} and~\ref{sec:enkf} regards the asset price 
$S_0$ as a known parameter. However, in practice there is uncertainty in this parameter.
In fact, we have an observed value which is in the best case an average over a day of trading,
so $S_0$ should be treated as an unknown with an observed mean value and a variance that is relatively
easy to estimate.
This in turn affects the calibration problem and its solution process. Section~\ref{sec:S0}
deals with this additional complication, which translates into uncertainty in the data locations
of a transformed formulation for~\eqref{x1}.

In Section~\ref{sec:test} we collect our numerical tests,
addressing and assessing the various aspects of the methods desribed earlier.
We use synthetic data in Section~\ref{sec:adjusts0scar} to show the advantage in applying
the method of Section~\ref{sec:S0} for problems with uncertainty in the price $S_0$.
In Section~\ref{sec:equity} we use market equity data to fine-tune our regularization functional,
as well as to compare Tikhonov-type regularization vs the modified artificial time EnKF.
In Section~\ref{sec:datacompcom}
we use oil and gas commodity market data to further investigate data completion approaches,
showing that the scarce data approach is superior.
Conclusions are offered in Section~\ref{sec:conclusions}.





\section{Two approaches for handling scarce data}
\label{sec:scarce}

Below we assume that the parameter $S_0$ is given.\footnote{In the related {\em online calibration}
setting, data are given for several values of $S_0$. 
For each such value of $S_0$ we then make a variable transformation and find a volatility surface.}
This assumption will be modified in Section~\ref{sec:S0}.
We then apply a standard transformation changing the independent variable $K$
to the so-called {\em log moneyness} variable $y = \log (K/S_0)$. 
This is followed by changing the dependent variables of the
forward and inverse problems to
$u(\tau,y) = C (\tau,S_0 \exp(y))$ and $a(\tau,y) = \frac 12 \sigma(\tau,K(y))^2$, respectively.
We obtain the dimension-less parabolic PDE with no unbounded coefficients
\begin{eqnarray}
-\frac{\partial u}{\partial \tau} + a\left( \frac{\partial^2 u}{\partial y^2} - 
\frac{\partial u}{\partial y}\right)  + b \frac{\partial u}{\partial y} &=& 0, \quad \tau > 0, y \in \Re , \label{x2}
\end{eqnarray}
subject to the side conditions
\begin{subequations}
\begin{eqnarray}
u(\tau = 0, y) &=& S_0(1 - \exp (y))^+ , \label{x21a} \\
\lim_{y \rightarrow \infty} u(\tau,y) &=& 0, \label{x21b} \\
\lim_{y \rightarrow -\infty} u(\tau,y) &=& S_0. \label{x21c}
\end{eqnarray}
We can write \eqref{x2}--\eqref{x21} as $L(a)u = q$, with the linear differential operator $L$ depending on $a$ and operating on $u$
and with the right hand side $q = q(S_0)$ given.
Thus, the forward problem involves finding $u$ satisfying this parabolic linear differential problem
for a given {\em local variance surface} $a$ and price $S_0$.
\label{x21}
\end{subequations}

To find a numerical solution for $L(a)u = q$ we first approximate the domain in $y$
by a finite interval, restricting $l_y \leq y \leq r_y$ for two real values satisfying $l_y < 0 < r_y$.
The boundary conditions~\eqref{x21b} and~\eqref{x21c} are then required to hold at $r_y$ and $l_y$, respectively.
Next, we
discretize the PDE problem on a mesh with a fixed step $\Delta \tau$ in $\tau$ and a fixed
step $\Delta y$ in $y$. Denote by $u_{i,j}$ the approximation of $u(i\Delta \tau, l_y+j\Delta y)$
and by $a_{i,j}$ the injection of the surface $a$ at $(i\Delta \tau, l_y+j\Delta y)$, $i = 1, \ldots , M_\tau+1$,
$j = 0, 1, \ldots , M_y+1$, where $(M_y+1)\Delta_y = r_y - l_y, \ (M_\tau + 1)\Delta \tau = T$. 
Then, using the Crank-Nicolson method~\cite{achpir},
we have the difference relations
\begin{eqnarray}
-\frac{u_{i+1,j}-u_{i,j}}{\Delta \tau} &+& \frac{a_{i,j}+a_{i+1,j}}{4\Delta y^2}
\left( u_{i+1,j+1}-2u_{i+1,j}+u_{i+1,j-1} + u_{i,j+1}-2u_{i,j}+u_{i,j-1}\right) \nonumber \\ 
&-&  \frac{a_{i,j}+a_{i+1,j} - 2b}{4\Delta y} 
\left( u_{i+1,j+1}-u_{i+1,j-1} + u_{i,j+1}-u_{i,j-1}\right) \nonumber \\
&=& 0, \quad i = 1, \ldots, M_\tau, \ j = 1, \ldots , M_y .
\label{x30}
\end{eqnarray}
An obvious treatment of the initial and (Dirichlet) boundary conditions closes this system of $M = M_\tau M_y$
equations that are linear for the variables $u_{i,j}$.
The mesh function $u_h$ for the approximation $u$ can be conveniently reshaped (say, ordered by column) into a vector 
$u_h \in \Re^M$, retaining the same notation without confusion.
Similarly, we obtain the mesh function $a_h$ as an injection of $a(\tau,y)$, reshaped into a vector if need be.
Then we can write~\eqref{x30} as
\begin{eqnarray}
L_h (a_h ) u_h = q_h , \quad {\rm or~} u_h = L_h (a_h )^{-1}q_h , \label{x3}
\end{eqnarray}
where $L_h$ is a sparse, nonsingular $M \times M$ matrix and $q_h$ is the mesh injection of $q$.
 
\begin{figure}[htb]
\begin{center}
\includegraphics[scale=.45]{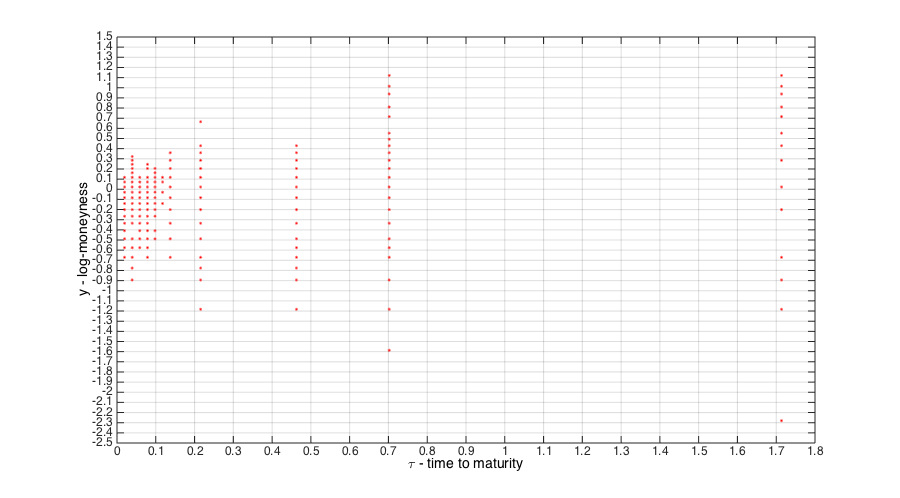}
\caption{Data locations for a PBR (Petrobraes, an oil company) 
set in the $(\tau,y)$ domain with our coarsest mesh
in the background. \label{fig_pbrdata}}
\end{center}
\end{figure}

The inverse problem is to find a volatility surface $\sigma$, approximated through $a_h$, 
that explains given observed data $d \in \Re^l$. 
These data values approximate $u$ at $l$ locations in the rectangular domain on which the problem~\eqref{x2}-\eqref{x21}
is defined. Typically, these locations are far from the boundaries and $l \ll M$; see Figures~\ref{fig_pbrdata}
and~\ref{fig_spxdata}.
Thus, the data set is {\em sparse}, or {\em scarce}, and the $l \times M$ matrix $P$ which maps grid locations for $u_h$
to those of $d$, using bilinear interpolation as necessary, has many more columns than rows.
The forward operator, which predicts the data for a given $a_h$, is the matrix-vector product (or projection) $Pu_h$.
The inverse problem is to find a plausible $a_h$ for which the predicted and observed data are sufficiently close,
as described below.

Note that the approximate solution of the inverse problem must be positive at all mesh points.
This positivity constraint turns out to hold automatically in all our reported calculations
(i.e., it is not an active constraint in the encountered optimization problems).

 
\subsection{Regularizing the inverse problem}
\label{sec:tikh} 

If the data $d$ has Gaussian noise $\sim \mathcal{N}(0,\Gamma)$, where $\Gamma$ is a
symmetric positive definite (SPD)
error covariance matrix,
then the maximum likelihood (ML) {\em data misfit} function is
\begin{equation}
		\phi(a_h) = \| P u_h (a_h) - d \|_{\Gamma^{-1}}^{2} , \label{x4}
		\end{equation}
where for an SPD matrix $C$ we define the vector energy norm $\| x \|_C = \sqrt{x^TCx} $. 
For instance, if $\Gamma = \alpha_0^{-1} I, \ \alpha_0 > 0$, then the {\em discrepancy principle}
(see, e.g., \cite{vogelbook,ehn1}) yields the stopping criterion
(i.e., we should find $a_h$ to reduce $\phi$ until)  
\begin{eqnarray} \phi(a_h) \leq \rho , \quad {\rm where~~} \rho = \alpha_0^{-1} l . \label{x5} \end{eqnarray}

However, there are in general many surface meshes $a_h$ that would satisfy the conditions~\eqref{x5}
(i.e., explain the data). 
We therefore introduce a regularization operator $R(a_h)$ which is a {\em prior}, in the sense that it represents prior
knowledge or belief about our sought surface.
We then take steps to minimize the MAP {\em merit function}
		\begin{equation} 
		\phi_{R}(a_h) =  \phi(a_h) +  R(a_h) . \label{x6}
		\end{equation}
%
%
Here we introduce the Tikhonov-type penalty function
\begin{eqnarray}
R(a_h) &=& \alpha_1 \sum_i\sum_j (a_{i,j} - a_0)^2 \nonumber \\
&+& \frac{\alpha_2}{\Delta \tau^2} \sum_j \sum_{i=1}^{M_\tau} (a_{i,j}-a_{i-1,j})^2
+ \frac{\alpha_3}{\Delta y^2} \sum_i \sum_{j=1}^{M_y} (a_{i,j}-a_{i,j-1})^2 , \label{x7}
\end{eqnarray}
where $\alpha_1, \alpha_2, \alpha_3 \geq 0$ are parameters to be determined.

The terms involving $\alpha_2$ and $\alpha_3$ correspond to smoothness of $a(\tau,y)$ in $\tau$ and $y$, respectively.
Thus, the values $a_{i,j}$ that are elements of the 2D mesh function $a_h$ 
should not fluctuate randomly, as they ought to form a reasonably smooth surface.
So, we insist on selecting $\alpha_2, \alpha_3 > 0$, and will modify a given EnKF algorithm accordingly as well. 
This penalizes lack of smoothness (or, {\em roughness}) in the reconstructed local variance surface.
Actual values for these parameters, determined by experimentation, are reported in Section~\ref{sec:test}. 

Next, the term governed by the weight $\alpha_1$ penalizes the distance from our sought surface
to an a priori function $a_0$ that we take to be a constant (not knowing better).
A reasonable value for this constant can often be estimated based on the risk category of the asset under consideration.
The importance of having this term is that without it the prior may tend to favour surface flatness,
whether or not this is realistic. But for the case of scarce data there could be too much freedom
in seeking merely a relatively flat surface which still explains the data.
Setting $\alpha_1 = 0$ might therefore cause the
solution process to behave unstably, as demonstrated further in Section~\ref{sec:test}.

Having determined values $\alpha_i, \ i = 1,2,3$ for the ensuing experiments
the next question is choosing the relative weight of the
data misfit function and the prior, which amounts to determining a value for $\alpha_0$
in the expression $\Gamma = \alpha_0^{-1} I$.
For this we have used the well-known L-curve method~\cite{vogelbook,ehn1}.



\subsection{The data completion approach}
\label{sec:datacomp}

Another approach to deal with the scarce data is to first interpolate/extrapolate the observed data to all 
$M$ locations of the $u$-mesh.
This can be done using the Kahale algorithm~\cite{kahale}, which represents a prior that reflects
financial considerations (namely, maintaining the {\em smile}, and more).
The inverse problem is subsequently defined on the thus-enriched data set,
the matrix $P$ in~\eqref{x4} becoming the identity,
and the Tikhonov-type regularization previously described is applied for its solution.

This data completion subsequently allows for the development of a more solid theory for the solution of the
corresponding regularized inverse problem.
 
A potential difficulty with this approach, however, is that the ensuing matching of the interpolated 
``data values'' to a field $u_h$ that approximately
satisfies the PDE problem~\eqref{x2}--\eqref{x21} is applied to data that can no longer be 
considered to have independent, random noise.
The interpolation/extrapolation operator and the PDE discrete Green's function operator could be in conflict. 
In Sections~\ref{sec:equity} and~\ref{sec:datacompcom} we describe examples
using real data sets from different markets which compare data completion using two different
techniques to the scarce data approach.
The obtained plots in Figures~\ref{spx6}, \ref{spx6a}, \ref{fig:com_hh} and \ref{fig:com_wti}
clearly demonstrate that the data completion approach can give inferior results.
In particular, the reconstructions using data completion
often agree neither with other curves nor with intuition.






\section{Artificial time EnKF-type methods}
\label{sec:enkf}

The homotopy, or continuation, approach of embedding a given problem in a larger one with artificial
time, subsequently defining an iterative method by ``advancing'' in the artificial time,
is old; see references in~\cite{ashudo}.
In our context there have been recent efforts that use such an iterative process to also adaptively learn
and improve knowledge of the error covariance~\cite{iglast,caerso}.
A Kalman filter setting is obtained as in~\cite{iglast} by defining for the ML data misfit function 
$\phi(a_h)$ of~\eqref{x4} the augmented state vector
\begin{eqnarray}
\hat{a}_h=\Psi(a_{h})= \begin{pmatrix}
a_h \\
Pu_h(a_h)
\end{pmatrix}, \label{enkf1}
\end{eqnarray}
together with the artificial dynamics (or prediction) $\hat{a}_{h}^{(n+1)} =  \Psi(a_{h}^{(n)})$.
The observations $d$ in~\eqref{x4} are then approximately matched by 
{ $H \hat {a}_{h}^{(n+1)}$}, with $H=\begin{pmatrix} 0 & I \end{pmatrix}$.
Further, the ensemble Kalman filter (EnKF) approach applies Monte Carlo approximations 
in order to obtain cheap estimates
for the error covariance matrices that appear in the Kalman filter method.

It is well-known that Kalman smoothing is equivalent to the solution of 
the corresponding {weighted} least squares problem,
while the Kalman filter end result agrees with that of the Kalman smoother.
However, what is not taken directly into account in~\cite{iglast,caerso} are the additional regularization terms
in the MAP merit function~\eqref{x6}.
In fact, we generalize~\eqref{x7} by considering
\begin{eqnarray}
\min_{a_h} \phi_R(a_h) &=& \big(d-Pu_h(a_h)\big)^T\Gamma^{-1}\big(d-Pu_h(a_h)\big)+
(a_h-a_0)^TD^{-1}(a_h-a_0)+\label{enkf2}\\
 &&(L_{\tau} a_0-L_{\tau}a_h)^TD_\tau^{-1}(L_\tau a_0-L_\tau a_h)+
 (L_ya_0-L_ya_h)^TD_y^{-1}(L_ya_0-L_ya_h), \nonumber
 \end{eqnarray}
where $L_\tau\in \Re^{(M_{\tau}-1)\times M_y}$  and $L_y\in \Re^{M_{\tau}\times (M_y-1)}$ are the scaled discrete sums
that multiply $\alpha_2$ and $\alpha_3$ there, respectively.
The matrices $D$, $D_\tau$ and $D_y$ are  corresponding error covariance matrices, 
with $D$ as yet unknown and to be determined in the EnKF process.

Next, we modify the data misfit part in~\eqref{enkf2} using the state augmentation approach~\eqref{enkf1}.
For this we redefine the matrices $L_\tau$ and $L_y$ by
\[ L_\tau \leftarrow \begin{pmatrix} L_\tau & 0 \end{pmatrix}, \ \  L_y \leftarrow \begin{pmatrix} L_y & 0 \end{pmatrix} .\]
The matrix $D$ is also modified properly, and we set 
\[ {\hat a}_0 = \begin{pmatrix} {a}_0​ \\ Pu_h ({a}_0​) \end{pmatrix}​.\] 
We can then rewrite~\eqref{enkf2} for the augmented variable ${\hat a}_h$ of~\eqref{enkf1} as
\begin{eqnarray}
\min_{\hat a_h} \phi_R(\hat a_h) &=& \big(d- H{\hat a}_h\big)^T\Gamma^{-1}\big(d-H {\hat a}_h\big)+
({\hat a}_h-{\hat a}_0)^TD^{-1}({\hat a}_h-{\hat a}_0)+\label{enkf4}\\
 &&(L_{\tau} {\hat a}_0-L_{\tau}{\hat a}_h)^TD_\tau^{-1}(L_\tau {\hat a}_0-L_\tau {\hat a}_h)+
 (L_y{\hat a}_0-L_y{\hat a}_h)^TD_y^{-1}(L_y{\hat a}_0-L_y{\hat a}_h). \nonumber
 \end{eqnarray}

Since~\eqref{enkf4} is just an enlarged weigthed least squares problem, it also corresponds to a Kalman smoother/filter
process, and we subsequently build an approximation for it as in~\cite{johns2008two}
by defining and solving a ``three-stage'' EnKF.

The details are somewhat tedious but straightforward, so we have gathered them in Appendix~\ref{sec:enkf_details}. 
We are led to the following algorithm, where we
let $\hat{a}_h^{(n,j)}$ denote the augmented state vector of the  $j$-th sample in the $n$-th iteration. 
\begin{enumerate}
\item {\tt Initialization}
\begin{enumerate}
\item
Generate $J$ samples of $a_h^{(0)}$, denoted
$\{a_h^{(0,j)}\}_{j=1}^{J}$, where $a_h^{(0,j)}\sim \mathcal{N}(a_0,D_0)$ and $D_0$ is an initial covariance matrix.  
\item
Compute $Pu_h(a^{(0,j)}_h), \ j = 1,\ldots,J$, thus defining $\{\hat{a}_h^{(1,j)}\}_{j=1}^J$. 
\end{enumerate}
\item  For $n=0,1,2,\ldots,$ until convergence criterion is satisfied, do:
\begin{enumerate}
\item {\tt Prediction step}
\begin{enumerate}
\item For $j = 1, \ldots , J$, set
$$\hat{a}_h^{(n+1,j)}=\Psi(a_h^{(n,j)})=
\begin{pmatrix}
a_h^{(n,j)}\\
Pu_h(a_h^{(n,j)})
\end{pmatrix}.$$
\item Define sample mean and covariance matrix as
\begin{eqnarray}
\bar{a}_h^{(n+1)}&=&\frac{1}{J}\sum_{j=1}^J \hat{a}_h^{(n+1,j)},\nonumber\\
D_{n+1}&=&\frac{1}{J}
\sum_{j=1}^J\hat{a}_h^{(n+1,j)}(\hat{a}_h^{(n+1,j)})^T- \bar{a}_h^{(n+1)}(\bar{a}_h^{(n+1)})^T.\nonumber
\end{eqnarray}
\end{enumerate}
\item {\tt Analysis step (three-stage Ensemble Kalman filter)} 

\begin{enumerate}
\item
Calculate
\begin{itemize}
\item $A^{(1)}_{n+1} = D_{n+1}H^T(HD_{n+1}H^T+
\Gamma)^{-1}$
\item $B^{(1)}_{n+1} = (I-A^{(1)}_{n+1}H)D_{n+1}$
\item $A^{(2)}_{n+1} = B^{(1)}_{n+1}L_{\tau}^T(L_{\tau}
B^{(1)}_{n+1}L_{\tau}^T
+D_{\tau})^{-1}$
\item $B^{(2)}_{n+1} = (I-A^{(2)}_{n+1}L_{\tau})B^{(1)}_{n+1}$
\item $A^{(3)}_{n+1} = B^{(2)}_{n+1}L_y^T(L_yB^{(2)}_{n+1}
L_y^T+D_y)^{-1}$
\item $B^{(3)}_{n+1} = (I-A^{(3)}_{n+1}L_y)B^{(2)}_{n+1}$
\end{itemize}
(Here $I$ is an identity matrix of appropriate size.)
\item
For $j = 1, \ldots , J$, update
\begin{eqnarray}
\tilde{a}^{(n+1,j)}_{h}&=&\hat{a}^{(n+1,j)}
_{h}+B^{(3)}_{n+1}\Big(H^T\Gamma^{-1}(d^{(j)}_{n+1}-H\hat{a}^{(n+1,j)}_{h}) \nonumber\\
&+&L_{\tau}^TD_{\tau}^{-1}(r_{\tau}^{(n+1,j)}-L_{\tau} \hat{a}^{(n+1,j)}_{h}) \nonumber\\
&+& L_y^TD_y^{-1}(r_y^{(n+1,j)}-L_y\hat{a}^{(n+1,j)}_{h})\Big) , \label{13n}
\end{eqnarray}
where $d^{(j)}_{n+1}=d+\eta^{(j)}_{n+1}$, $\eta^{(j)}_{n+1}\sim \mathcal{N}(0,\Gamma)$;
$r_{\tau}^{(n+1,j)}$ and $r_y^{(n+1,j)}$ are sampled from 
$\mathcal{N}(L_{\tau}\bar{a}_h^{(n+1)}, D_{\tau})$ and $\mathcal{N}(L_y\bar{a}_h^{(n+1)}, D_y)$, respectively.
\item
For $j = 1, \ldots , J$, set
$$a_h^{(n+1,j)}=\begin{pmatrix} I &  0\end{pmatrix} \tilde{a}_{h}^{(n+1,j)} .$$
\end{enumerate}
\item {\tt Convergence test}

Compute
$$a_h^{(n+1)}=\frac{1}{J}\sum_{j=1}^Ja_h^{(n+1,j)} $$
and check for convergence.
\end{enumerate} 
\end{enumerate}


Assuming that this algorithm stops after $N$ iterations, it requires $(N+1)J$
solutions of the forward problem.


\section{Uncertainty in the asset price}
\label{sec:S0}

In a typical image processing application of deblurring, denoising or inpainting, data values are
prescribed at pixels, so the measurement locations are known. However, when denoising a point cloud
or a surface mesh, for instance, there is no such distinction between a datum value and location: 
the unknowns are nodal mesh points in 3D, and as such they live in a higher dimensional space
(namely, $3$ rather than $1$).
This difference affects portability of image processing algorithms to similar problems on
surfaces~\cite{huas}.

A more subtle instance of location uncertainty arises in our volatility surface context.
The observed value for $S_0$ is actually an average of a day's trading (say), 
and as such contains uncertainty
whose variance can fortunately be directly estimated. In this section we thus make $S_0$ an unknown
that can be adjusted but should not stray too much from a measured 
(and thus observed) average value $\widehat S_0$.

This affects the regularization prior, which now depends also on $S_0$, so $R = R(a_h,S_0)$,
and has a term added to~\eqref{x7}:
\begin{subequations}
\begin{eqnarray}
R(a_h,S_0) &=&  \alpha_1 \sum_i\sum_j (a_{i,j} - a_0)^2 + \alpha_5 (S_0 - \widehat S_0)^2 \label{S0vara} \\
&+& \frac{\alpha_2}{\Delta \tau^2} \sum_j \sum_{i=1}^{M_\tau} (a_{i,j}-a_{i-1,j})^2
+ \frac{\alpha_3}{\Delta y^2} \sum_i \sum_{j=1}^{M_y} (a_{i,j}-a_{i,j-1})^2 . \nonumber
\end{eqnarray}
The additional parameter $\alpha_5$ is determined by the daily price variance. 

Furthermore, in the dimension-less form~\eqref{x2} of the Dupire PDE
problem, which to recall allows bounded PDE coefficients and uniform discretization step sizes in $(\tau,y)$,
the independent log moneyness variable $y = \log (K / S_0)$ now contains uncertainty as well!
We therefore update also the misfit function~\eqref{x4} to read
\begin{eqnarray} 
		\phi (a_h ,S_0) &=& \| P(S_0) u_h (a_h) - d \|_{\Gamma^{-1}}^{2} \label{S0varb} \\
	&+& \alpha_4  \sum_{j=1}^{M_y}\left((1 - \exp (y_j(S_0)))^+ -  (1 - \exp (y_j(\widehat S_0)))^+\right)^2 . \nonumber
\end{eqnarray}
\label{S0var}
\end{subequations}%

The parameter $\alpha_4$, like $\alpha_5$, is determined from the variance in $S_0$.
In fact, in our experiments we have found that it is safe to set $\alpha_5 = 0$ and
control the uncertainty penalty through $\alpha_4$ alone.
Note that in the data projection matrix, $P = P(S_0)$.
Our optimization problem replacing~\eqref{x6} is now
\begin{equation} 
		\min_{a_h,S_0}~\phi_{R}(a_h,S_0) =  \phi(a_h,S_0) +  R(a_h,S_0) . \label{x6S0}
		\end{equation}

To solve the extended optimization problem~\eqref{x6S0}
we apply a splitting method, 
alternately minimizing \eqref{x6S0}
for $a_h$ and for $S_0$. When $S_0$ is held fixed, the problem returns to that considered in 
Sections~\ref{sec:scarce} and~\ref{sec:enkf}.
When $a_h$ is frozen in turn, the remaining minimization problem for $S_0$ is in 1D
and causes no difficulty.
Furthermore, fortunately the coupling between these variables
is weak, so fast convergence of this splitting method is observed in all our experiments.


\section{Testing our methods on real and synthetic data}
\label{sec:test}

In this section we present a number of tests that were performed in order to illustrate the points made in the previous
sections as well as to compare the different methodologies. 
Specifically,
we compare results related to the introduction of uncertainty in the value of the 
reference price $S_0$, the effect of choosing between the EnKF and the Tikhonov-type 
regularization approach, and the introduction of a penalty associated to the mean value of the local variance surface $a_h$. 

The computations for all the examples described in this paper
were performed using garden-variety personal computers,
with typical runtimes clocking within a few minutes. Thus, being concerned in this work
with addressing several fundamental questions,   
no special effort was made to optimize the runtime performance of our codes.

The basic setting of our numerical experiments consists of solving the inverse problem described in Section~\ref{sec:volatility},  
first using synthetic data (tainted by multiplicative Gaussian noise) and then using real data. 
In the  synthetic data examples, we first assume a ground truth local variance surface $a_{true}$ on a fine grid, from which we solve the discretized PDE for $u_h$ as in Section~\ref{sec:scarce}.
From that we sample data on a coarser mesh (so as to avoid a so-called inverse crime) and then subject the data to noise. 
In the real data examples we selected publicly available option data from the NY stock exchange. 
The reconstructions are performed by using the techniques put forth in Sections~\ref{sec:scarce} and \ref{sec:enkf}.
%

\subsection{Results for uncertain $S_0$ using synthetic data}
\label{sec:adjusts0scar}

In order to test our techniques in situations where we have full control of the unknown volatility surface, 
we postulate a local volatility surface with known form . 
We concentrate on the issue of underlying price uncertainty, and as we shall see, this can be well addressed by the method 
discussed in Section~\ref{sec:S0}.


\subsubsection{First experiment with underlying price uncertainty}
We start by experimenting with the uncertainty about the value of $S_0$ for a ground truth
volatility given by
\begin{equation}
\sigma(\tau,y) = \left\{
\begin{array}{ll}
\displaystyle\frac{2}{5}-\frac{4}{25}\text{e}^{-\tau/2}\cos\left(\displaystyle\frac{4\pi y}{5}\right),& \text{ if } -2/5 \leq y \leq 2/5\\
\\
2/5,& \text{ otherwise.}
\end{array} \right.
\label{vol}
\end{equation}
From \eqref{vol} we produced call option prices by solving the PDE system~\eqref{x2}--\eqref{x21} discretized  
over a  mesh with step sizes $\Delta \tau = 0.005$ and $\Delta y = 0.025$. 
The maximum time to maturity is $\tau_{\max} = 0.5$ year, and the minimum and maximum log-moneyness strikes are $y_{\min}= -5$ and $y_{\max} = 5$, respectively. 
We also chose, for simplicity, the true underlying asset price as $S_0 = 1.0$ and the risk-free interest rate as $r = 0$. 

The computed call prices $u_{i,j}$ were then polluted by a relative noise,
\begin{equation}
 u_{i,j}^{\text{noisy}} = u_{i,j}(1  + 0.01 \eta_{i,j}),
 \label{eq:noise}
\end{equation}
with $\eta_{i,j}$ drawn from the standard normal distribution $\mathcal N (0,1)$.
The data is then composed on a coarse mesh by the noisy prices $u^{\text{noisy}}_{i,j}$ with 
$\tau_i = i\delta\tau$, $i=1,...,5$, $\delta\tau = 0.1$, and $y_j = -0.75 + j\delta y$, $\delta y = 0.05$, $j=0,1,...,30$. 
The observed underlying price is $\hat{S_0} = 0.95$.

As mentioned in Section~\ref{sec:S0}, the minimization of the Tikhonov-type functional \eqref{x6S0} is 
achieved by alternating minimizations, namely, 
the minimization w.r.t. $a_h$, which is performed by a gradient descent method as in~\cite{AlbAscZub2015}, 
and w.r.t. $S_0$, which is performed with the {\sc Matlab} function {\tt lsqnonlin}. 
In both stages, the mesh width used are $\Delta \tau = 0.01$ and $\Delta y = 0.05$, 
the minimum and maximum log-moneyness values are taken as $-5$ and $5$, respectively, and the maximum time to maturity is $0.5$.

We started the algorithm with $a^0_h \equiv 0.45^2/2$ and $S^0_0 = \hat{S_0} = 0.95$.
During the minimization w.r.t. $a_h$, we took the parameters in the penalty function as  
$\alpha_0 = 10^{5}$, $\alpha_1 = 10^2$, $\alpha_2 = 10^4$, $\alpha_3 = 1$, and $\alpha_4=\alpha_5=0$, 
whereas in the minimization w.r.t. $S_0$, we used
$\alpha_1=\alpha_2=\alpha_3=0$, and $\alpha_4=\alpha_5=10^5$. The {\em a priori} surface was taken as $a_0 \equiv 0.45^2/2$.

\begin{figure}[!ht]
    \centering
           \includegraphics[width=0.2\textwidth]{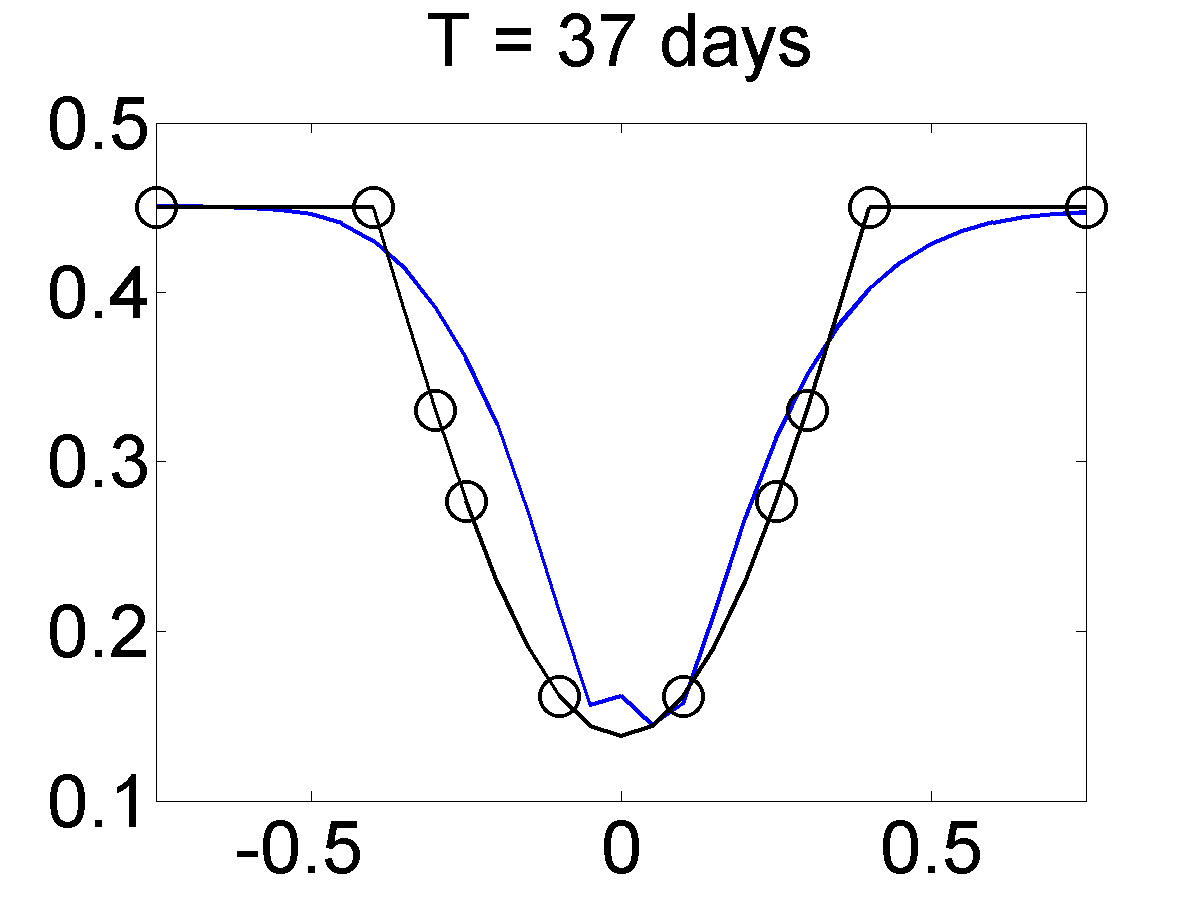}\hfill
           \includegraphics[width=0.2\textwidth]{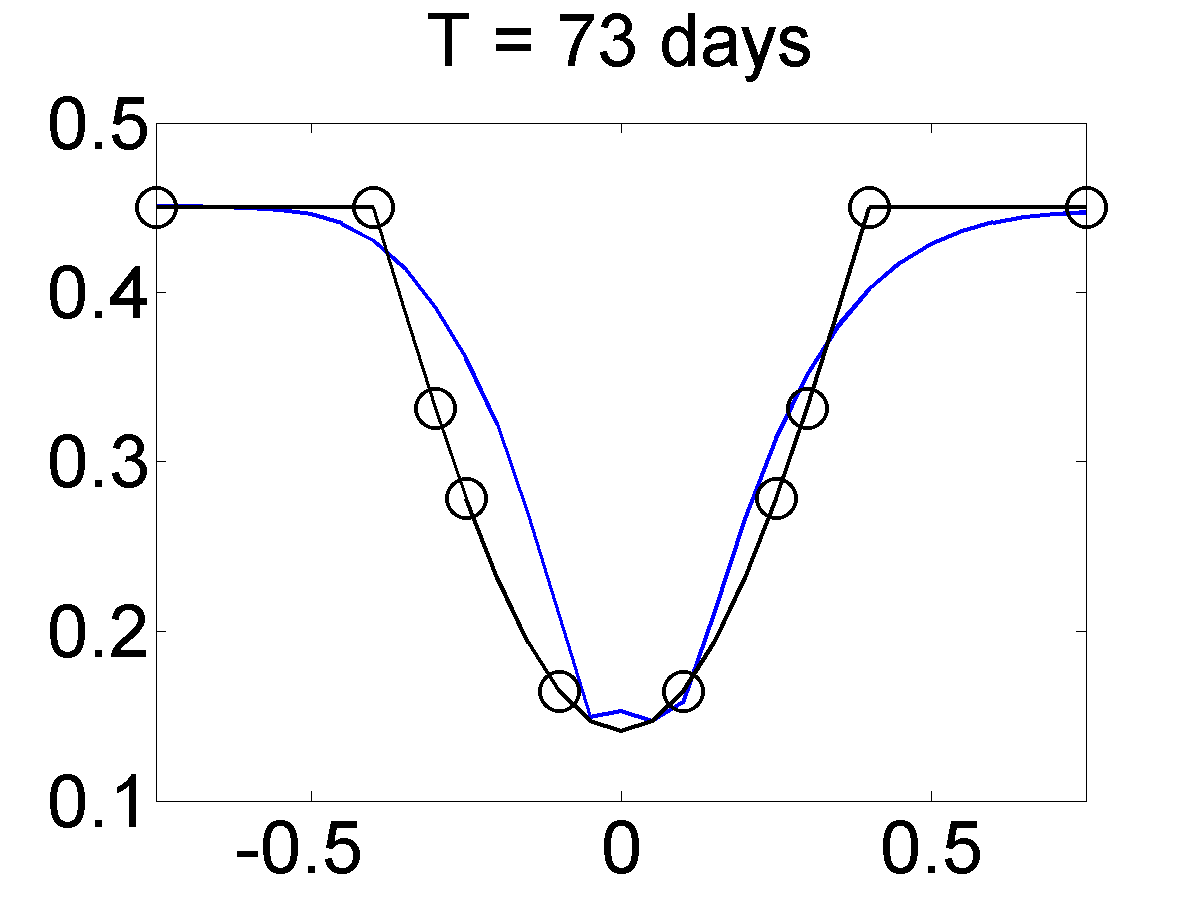}\hfill
 	\includegraphics[width=0.2\textwidth]{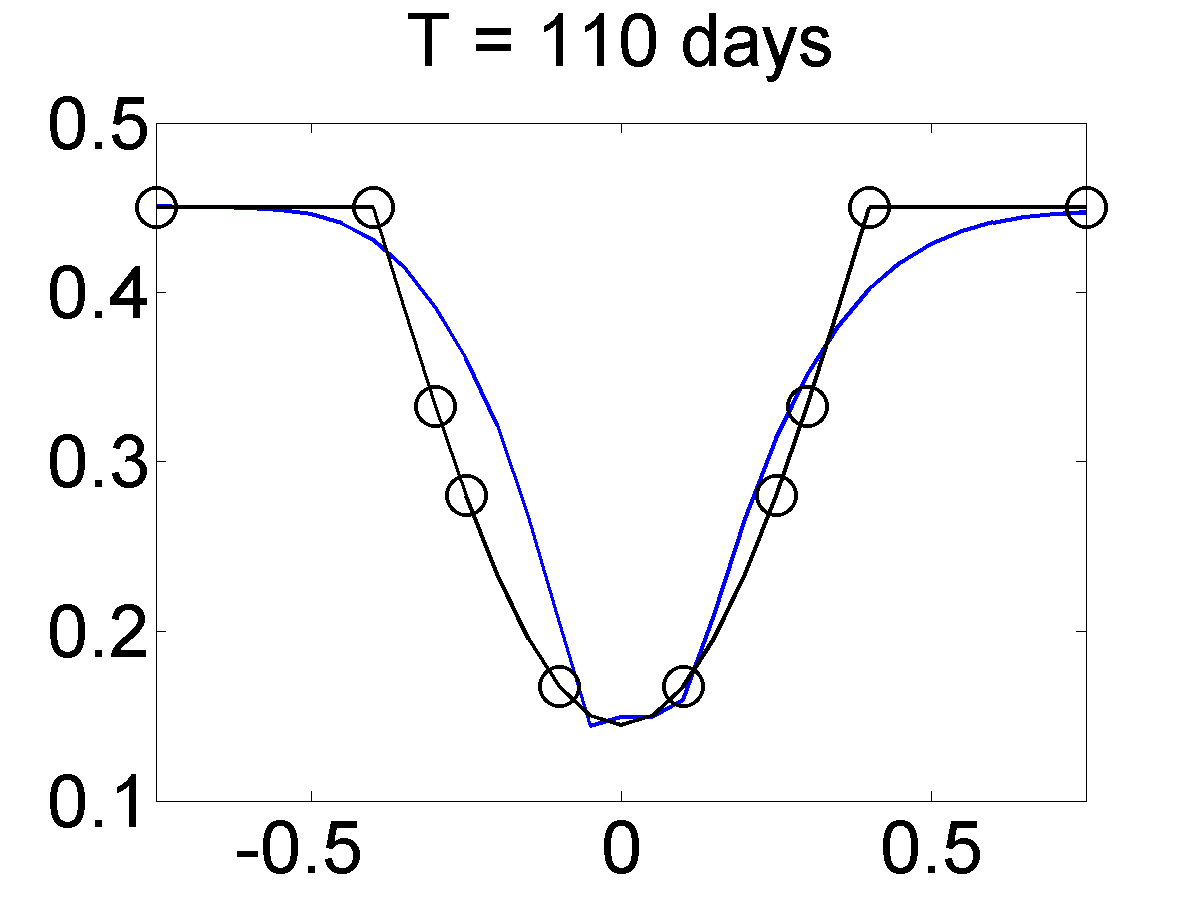}\hfill
 	\includegraphics[width=0.2\textwidth]{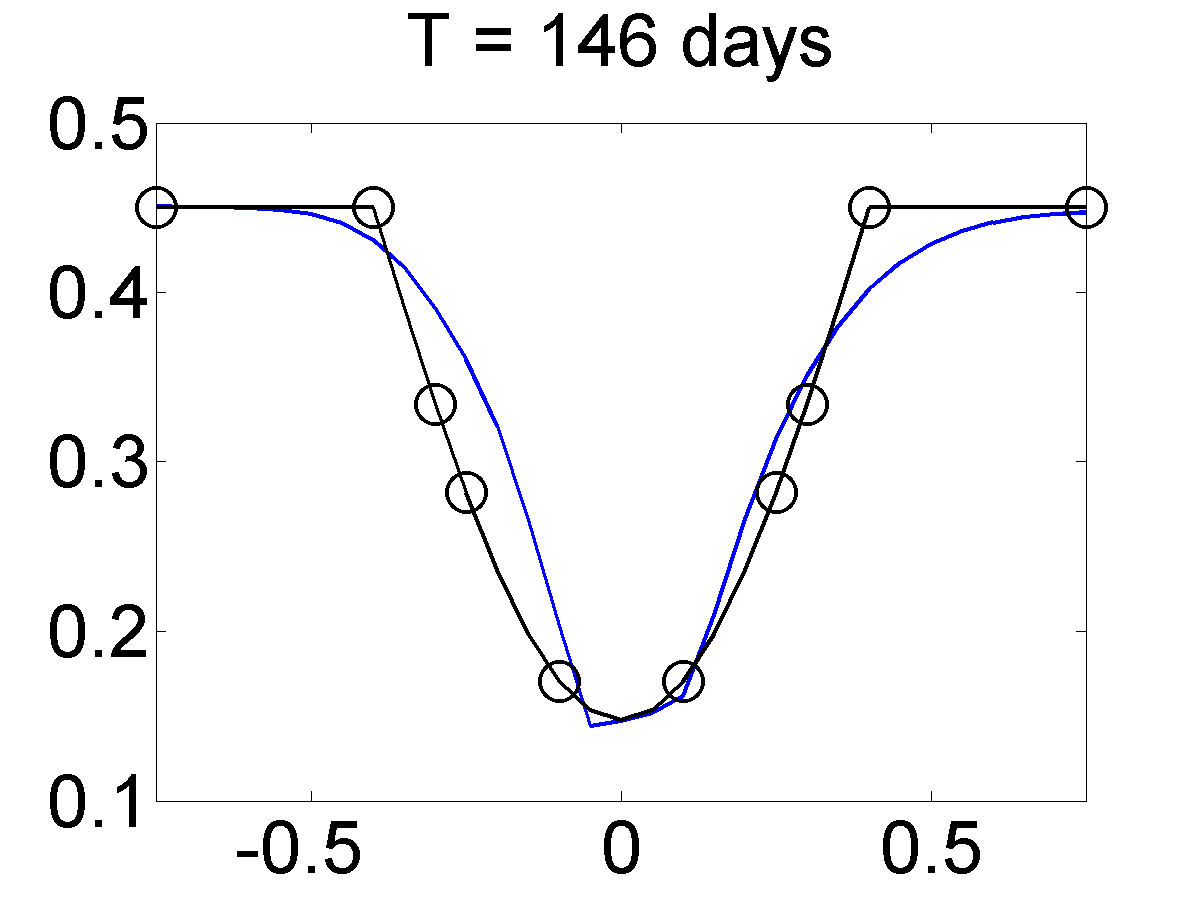}\hfill
 	\includegraphics[width=0.2\textwidth]{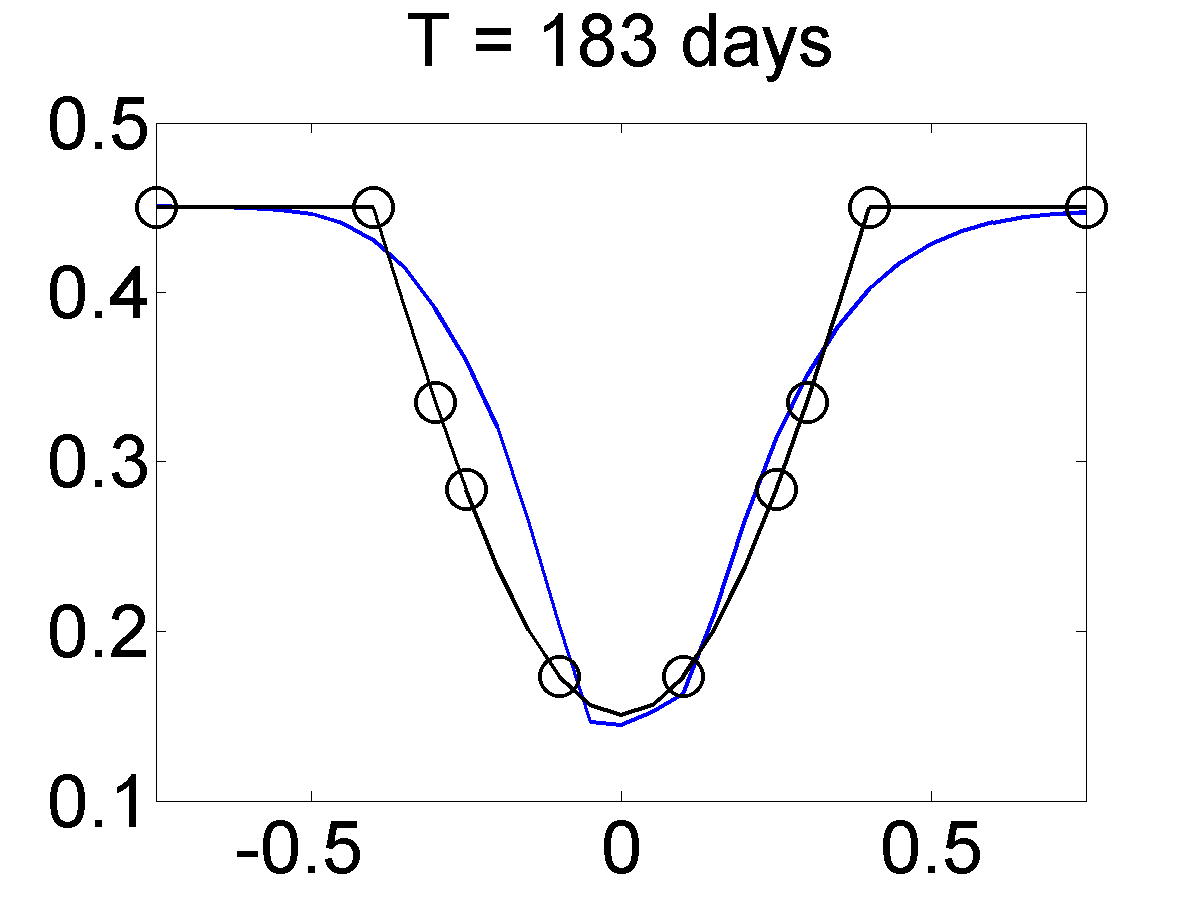}
  \caption{Reconstructed (continuous line) and true (line with circles) local volatility surfaces at the five different maturities. The reconstructed local volatility surface corresponds to the last one obtained with the adjustment algorithm of the underlying asset $S_0$.}
   \label{fig:scars01}
 \end{figure}

After $8$ iterations, the resulting underlying asset price was $0.999$, 
and the normalized $\ell_2$-distance between the true and the reconstructed local volatility surfaces was $0.13$. 
At the beginning, with $S_0 = 0.95$, the normalized distance was $5.55$. 
Figure~\ref{fig:scars01} compares between the reconstructed local volatility surface and the ground truth at each maturity.  
With the same scarce data set, as the underlying asset price was adjusted, 
we found a local volatility surface much closer to the original one. 
Table~\ref{tab:scars01} presents the evolution of the normalized $\ell_2$-distance between the true and the reconstructed local volatility surfaces, illustrating the latter observation.
 \begin{table}[!ht]
 \centering
 \caption{Normalized $\ell_2$-distance between  the true and the reconstructed local volatility surfaces and the value of $S_0$ at each step of the algorithm for adjusting $S_0$.}\label{tab:scars01}
 \begin{tabular}{|c|cccccccc|}
 \hline
 Iteration & 1 & 2 & 3 & 4 & 5 & 6 & 7 & 8\\
 \hline
 Normalized Distance & $5.55$ & $3.41$ & $2.39$ & $1.22$ & $0.78$ & $0.47$ & $0.21$ & $0.13$\\
 $S_0$ & $0.950$ & $0.963$ & $0.977$ & $0.985$ & $0.989$ & $0.994$ & $0.007$ & $0.999$\\
 \hline
 \end{tabular}
 \end{table}

\subsubsection{Another experiment with underlying price uncertainty} 
\begin{figure}[H]
   \centering
       \includegraphics[width=0.475\textwidth]{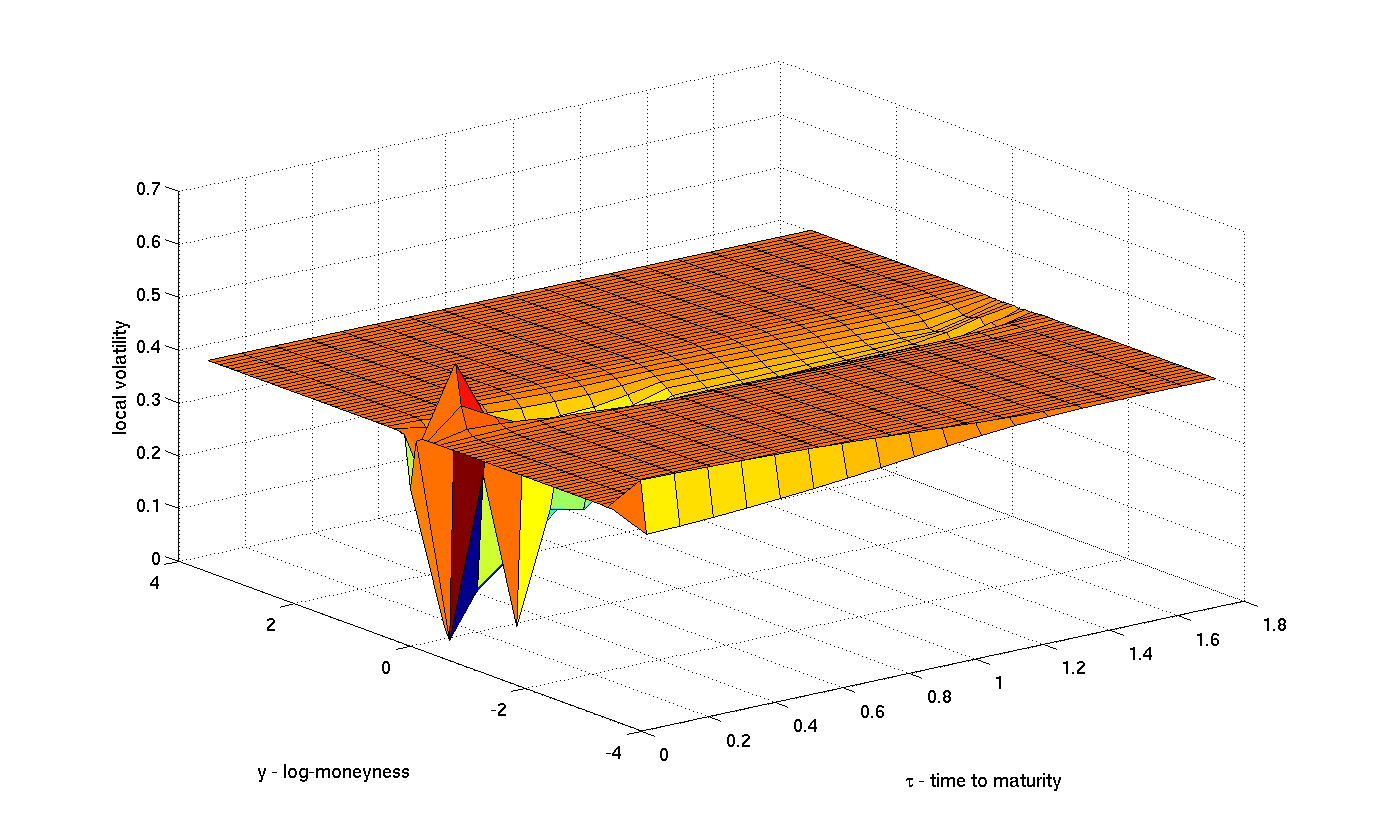}\hfill
       \includegraphics[width=0.475\textwidth]{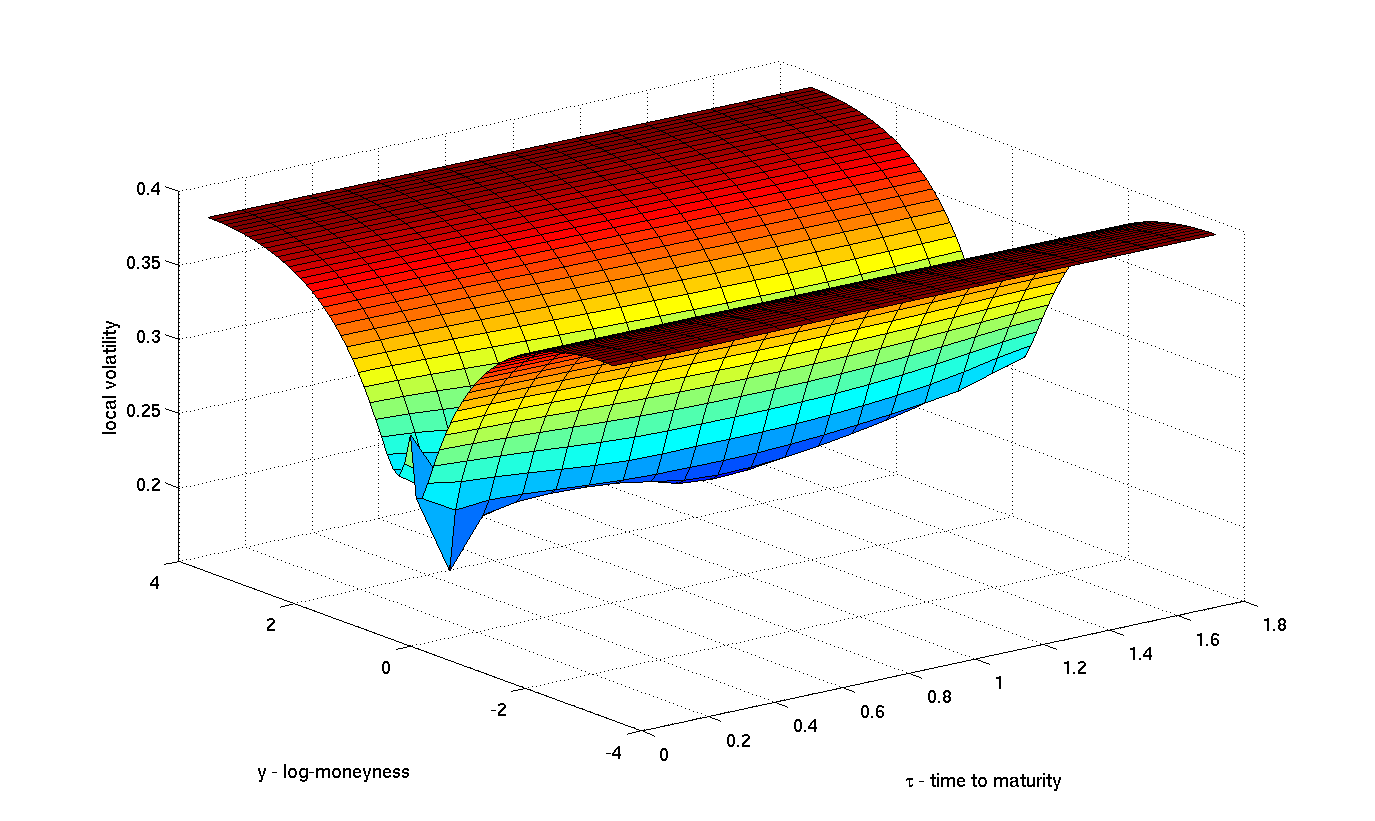}\hfill
       \includegraphics[width=0.475\textwidth]{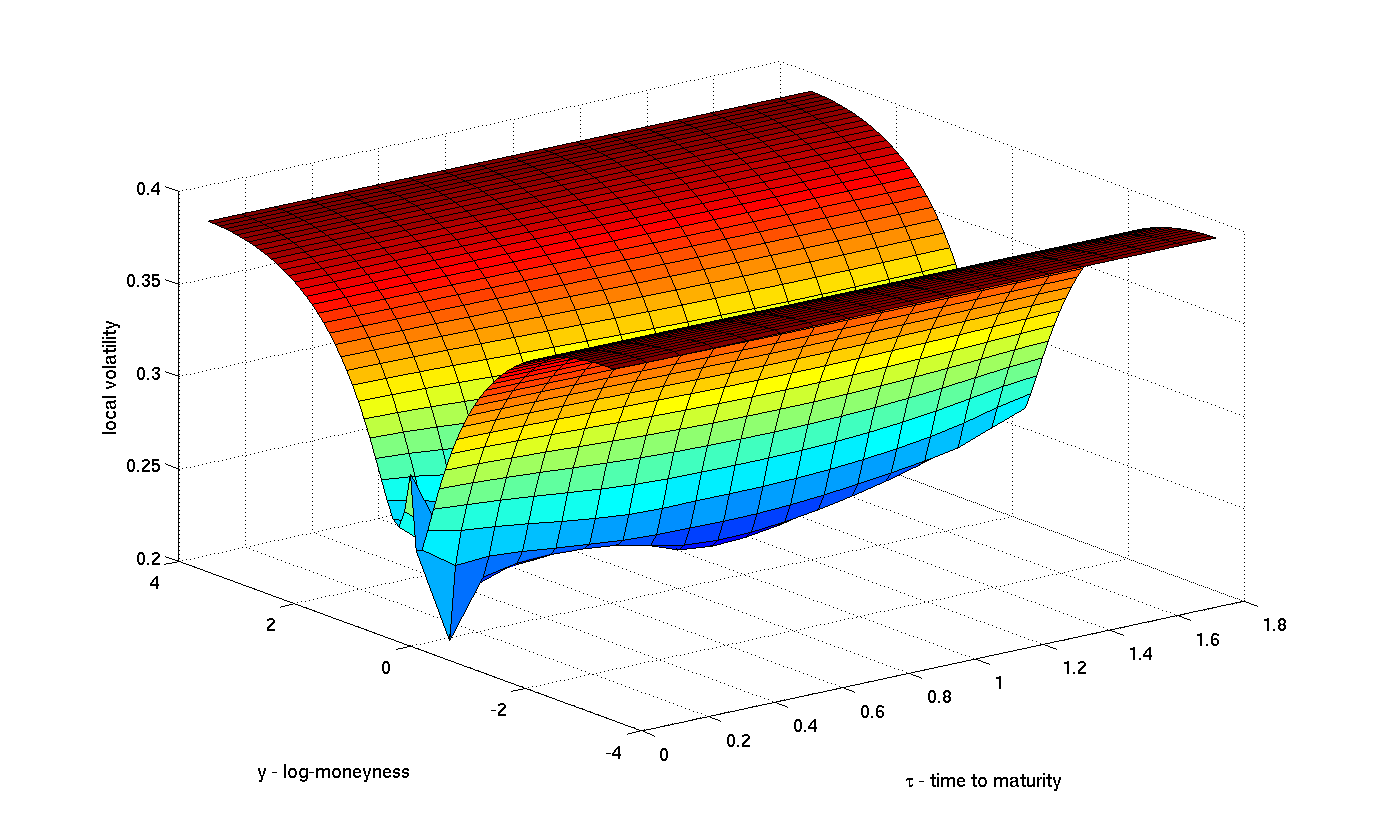}\hfill
       \includegraphics[width=0.475\textwidth]{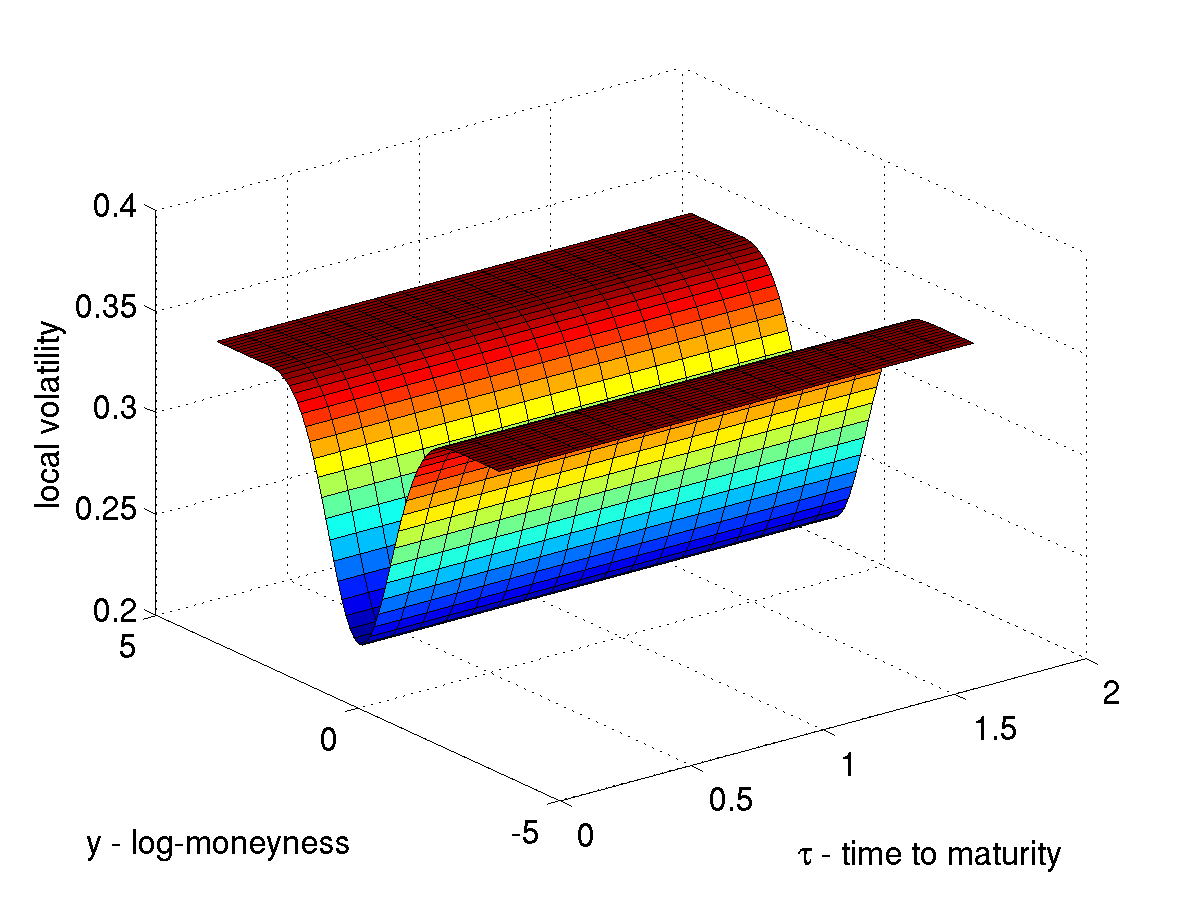}\hfill
   \caption{Calibration of the local volatility in 5 iterations. 
    Shown from the upper left, clockwise, are the 
   1st iteration, 3rd iteration, 5th iteration and the ground truth. 
   }
  \label{syn_1}
\end{figure}

We have carried out another experiment which deals with the underlying price uncertainty. 
The same general form as given in~\eqref{vol} is used for 
the volatility surface, except for a little bit of smoothing: see the bottm right subfigure in Figure~\ref{syn_1}. 
However, we make sure to employ 
parameters that are close to the ones generally encountered with real data, 
and we use data that follows a practical grid. The parameters for the example are given in Table~\ref{tabsyn1}.

\begin{table}[H]
\centering
\caption{{\small Parameters for the example of Figure~\ref{syn_1}. }\label{tabsyn1}}
\begin{tabular}{|c|c|}
\hline
$\widehat{S}_0$ initial spot price & 2500 \\
$S_{true}$ optimal spot price & 2200 \\
$r$ interest rate & 0.25\%  \\
the maximum maturity & 1.8  \\
Minimum $y$& -3.5 \\
Maximum $y$& 3.5  \\
$\Delta \tau$  & 0.1 \\
$\Delta y$ & 0.1\\
{\em a priori} surface $a_0$ & $0.4^2/2$  \\
\hline
\end{tabular}
\end{table} 

\begin{figure}[H]
\begin{center}
\includegraphics[height=0.2\textheight]{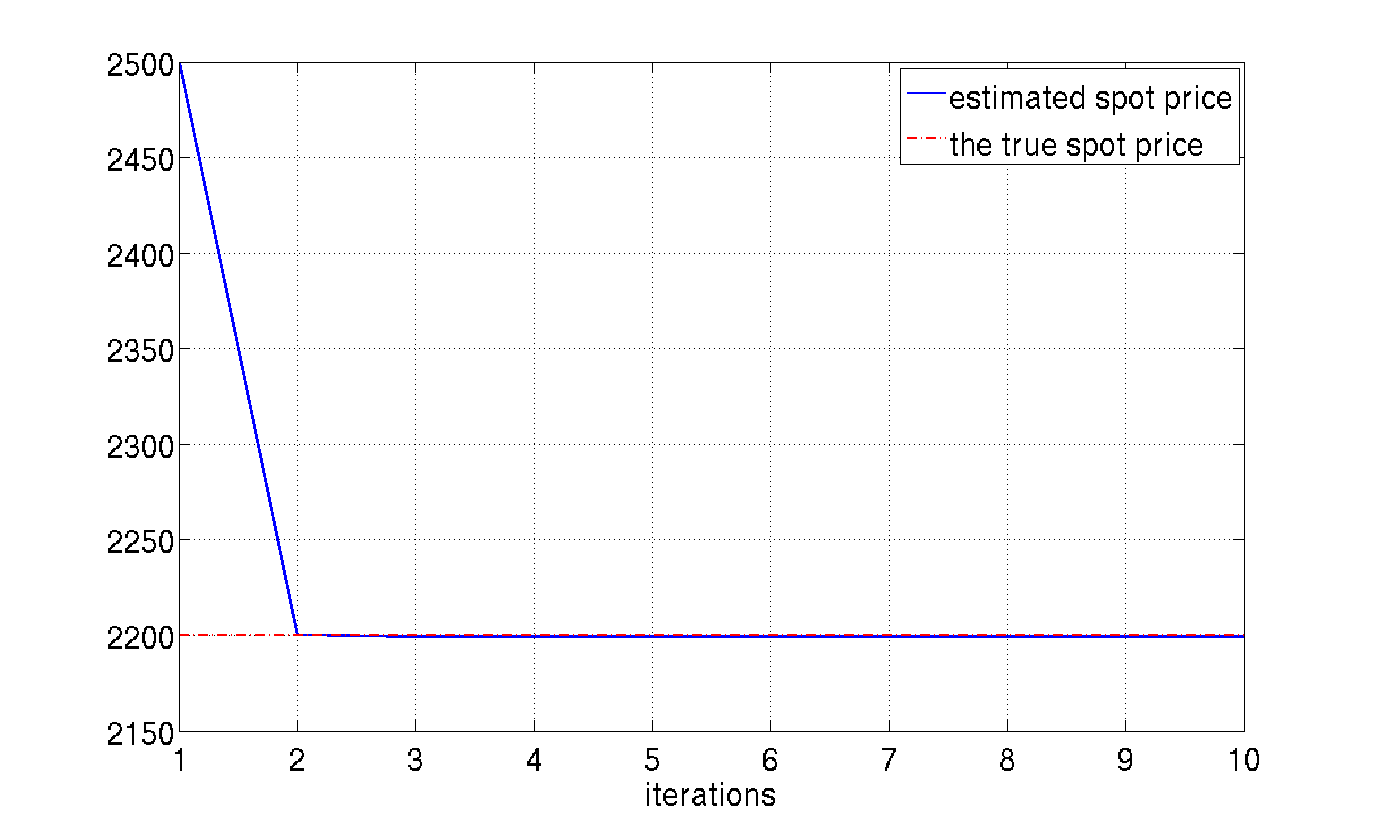}\\
\caption{The estimated spot price converges to the true price. \label{spotconvsyn_1}
}
\end{center}
\end{figure}

Following the same algorithm as before, we obtain that $S_0$ approximates
$S_{true}$ well, and furthermore, the local variance $a_h$ approximates $a_{true}$ on a coarse grid.  
This is illustrated in Figures~\ref{syn_1} and \ref{spotconvsyn_1}. 

\paragraph{Discussion of the synthetic data results}

The experiment with synthetic data indicates that including the underlying stock price as one of the unknowns 
(which corresponds to handling data location uncertainty in $y$)
leads to better 
results when there is uncertainty about such value. 
Indeed, Table~\ref{tab:scars01} and Figures~\ref{fig:scars01} and \ref{syn_1} 
show that the normalized distance between the reconstructions and the true surface decreases considerably 
when we combine local volatility calibration with the adjustment of the underlying asset price, 
upon using the algorithm presented in Section~\ref{sec:S0}. 
In the first experiment, the distance of the initial price and the price was relatively small, whereas in the second one, the two prices were significantly apart. In both cases, we can see that after only a few iterations, the prices were well-approximated. 

\subsection{Results for equity data}
\label{sec:equity}
\begin{figure}[htb]
\begin{center}
\includegraphics[scale=.45]{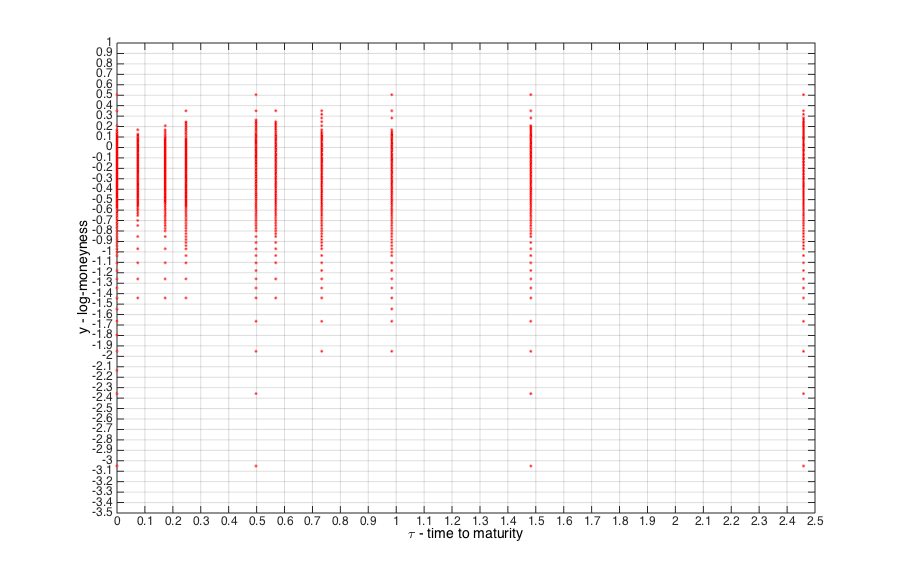}\\
\caption{Locations of the SPX data in the $(\tau,y)$ domain with our coarsest mesh
in the background. \label{fig_spxdata}}
\end{center}
\end{figure}

In this subsection and the next we consider real data from financial markets. 
Such markets are a tremendous source of data which can 
be promptly accessed by researchers. In our experiments we chose options on the Standard and Poors (SPX) index.\footnote{The Standard and Poors is a weighted index of actively traded large capitalization common stocks in the United States.}
Figure~\ref{fig_spxdata} depicts the locations at which our data set is given.

Such options are fairly liquid ones and thus amenable to the models introduced in Section~\ref{sec:volatility}. 
The data were collected on 19-Jun-2015 and contain prices for 9 different maturities ranging  from 1 day to 2 years. 
The parameters for all the models are given in Table~\ref{tabspx1}. 
Note that the optimal spot price 
in the table refers to the optimized spot price
using the method described in Section~\ref{sec:S0}. 
%

\begin{table}[H]
\centering
\caption{{\small Parameters for the equity data examples.}}\label{tabspx1}
\begin{tabular}{|c|c|}
\hline
$S_0$ initial spot price & 2112.7 \\
$S_0$ optimal spot price & 2095.6 \\
$r$ interest rate & 0.25\%  \\
the maximum maturity & 2.5  \\
Minimum $y$& -4.5 \\
Maximum $y$& 1.5  \\
$\Delta \tau$  & 0.05 \\
$\Delta y$ & 0.1\\
initial $a_0$ & $0.14^2/2$  \\
\hline
\end{tabular}
\end{table}

The parameters $\alpha_0$, $\alpha_1$, $\alpha_2$ and $\alpha_3$ in the penalty functional \eqref{x7} or \eqref{enkf2} used in this experiment can be found in Table~\ref{tabspx2}.

\begin{table}[htb]
\centering
\caption{{\small Parameters of the penalty functional \eqref{x7} or \eqref{enkf2} with SPX data.}\label{tabspx2}}
\begin{tabular}{|l|c|c|c|c|}
\hline
Parameter & $\alpha_0$ & $\alpha_1$ & $\alpha_2$ & $\alpha_3$ \\
\hline
Value     & 4.e+8 & 1.e+6 or 0 & 1.e+6 & 1.e+6 \\
\hline
\end{tabular}
\end{table}


\begin{table}[htb]
\centering
\caption{{\small Residuals of the 6 method variants.}}\label{residual_6}
\begin{tabular}{c|c|c|c|c|c|c|}
\cline{2-7}
 &\multicolumn{4}{c|}{Tikhonov-type} & \multicolumn{2}{c|}{EnKF}\\
\cline{2-7}
 &Scarce&Comp.&Scarce (no $a_0$)& Comp. (no $a_0$)&Scarce&Comp.\\
\hline
\multicolumn{1}{|l|}{Residual} & 0.0196 & 0.0314 & 0.0247 & 0.0289 & 0.0198 & 0.0294 \\
\hline
\end{tabular}
\end{table}

\begin{figure}[H]
\begin{center}
\includegraphics[width=.25\textwidth]{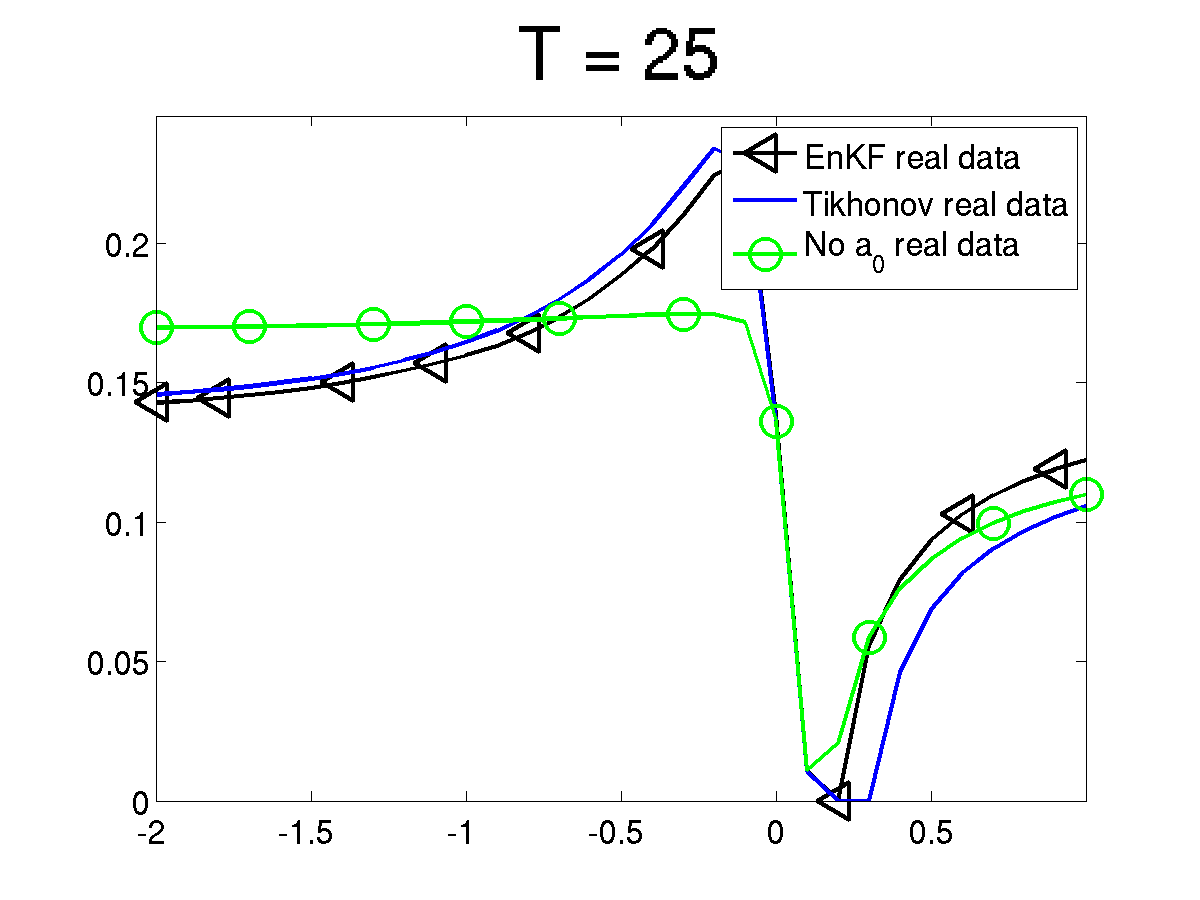}\hfill
\includegraphics[width=.25\textwidth]{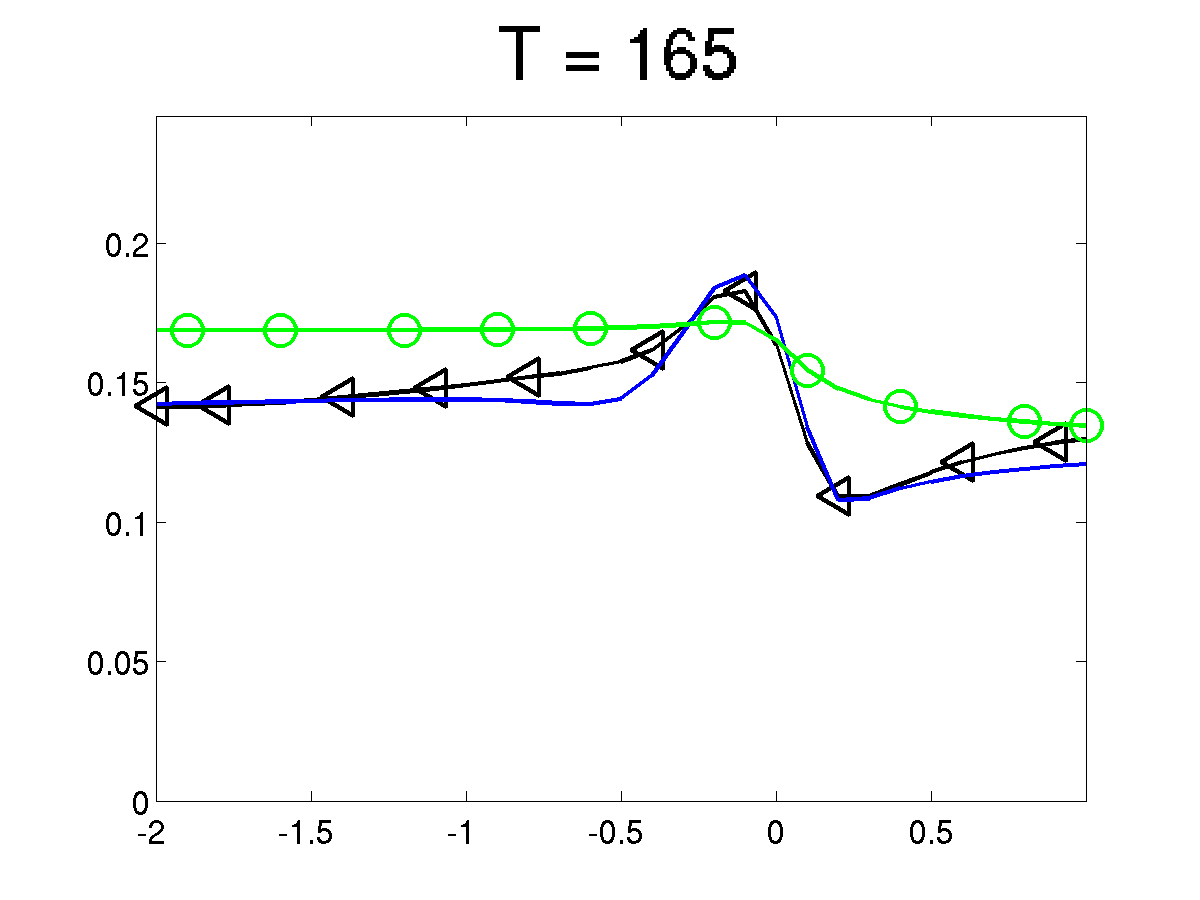}\hfill
\includegraphics[width=.25\textwidth]{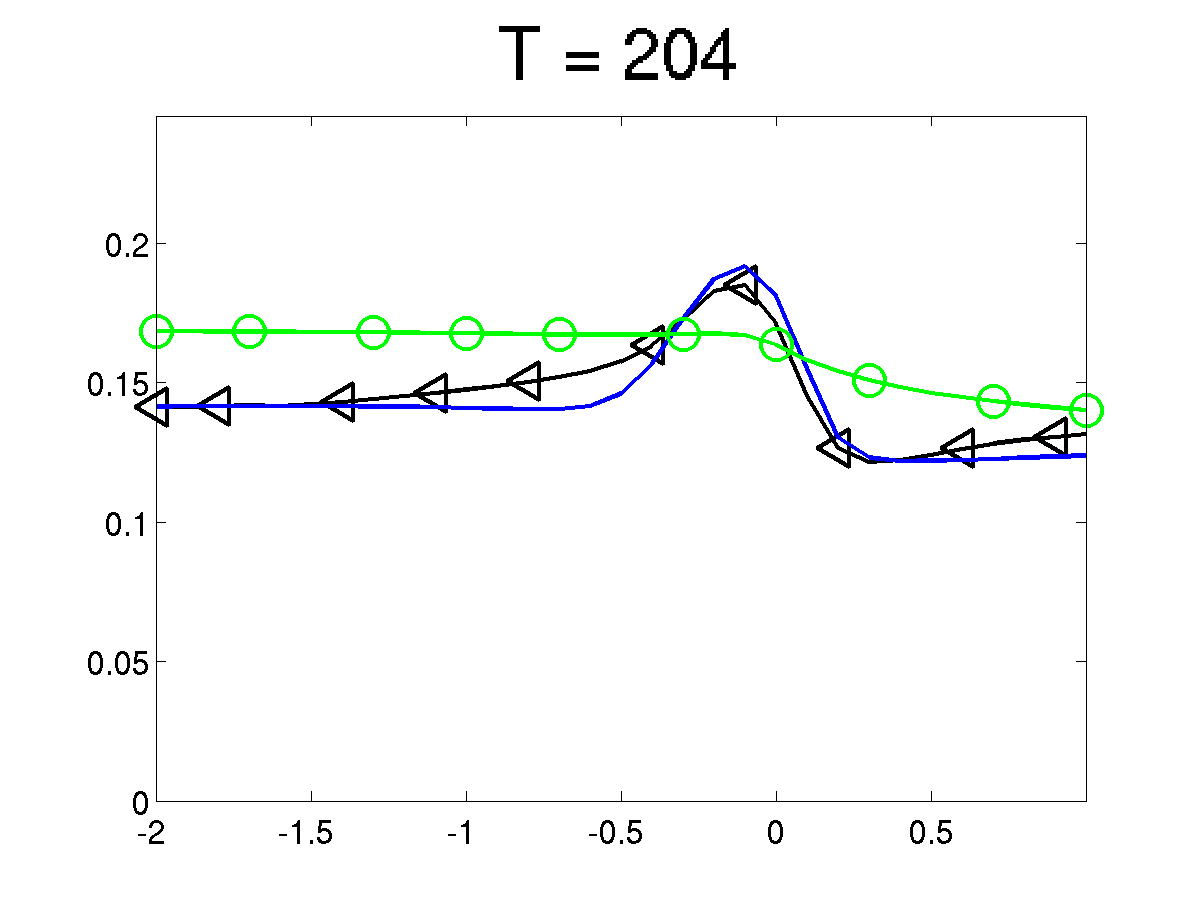}\hfill
\includegraphics[width=.25\textwidth]{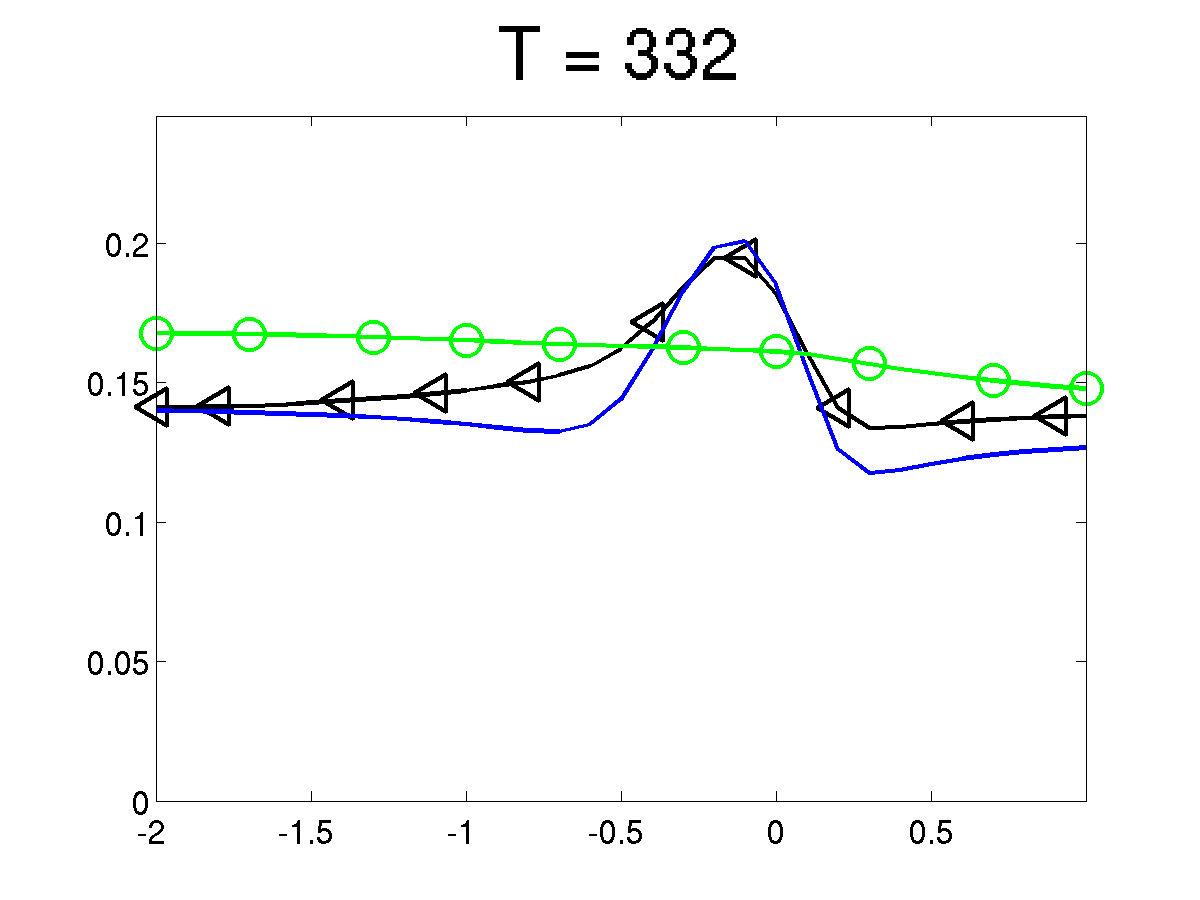}
\caption{\label{spx_r} Reconstructed SPX local volatility surfaces at different maturities obtained with three
method variants using scarce data.}
\end{center}
\end{figure}

Figure~\ref{spx_r} displays 
reconstructed SPX local volatility surfaces at different maturities obtained with three method variants using the original scarce data.
Note that the results generated by the Tikhonov-type method with $a_0$ penalty and EnKF are closer, whereas the
results without the $a_0$ penalty differ inexplicably.
  
\begin{figure}[H]
\begin{center}
\includegraphics[width=.25\textwidth]{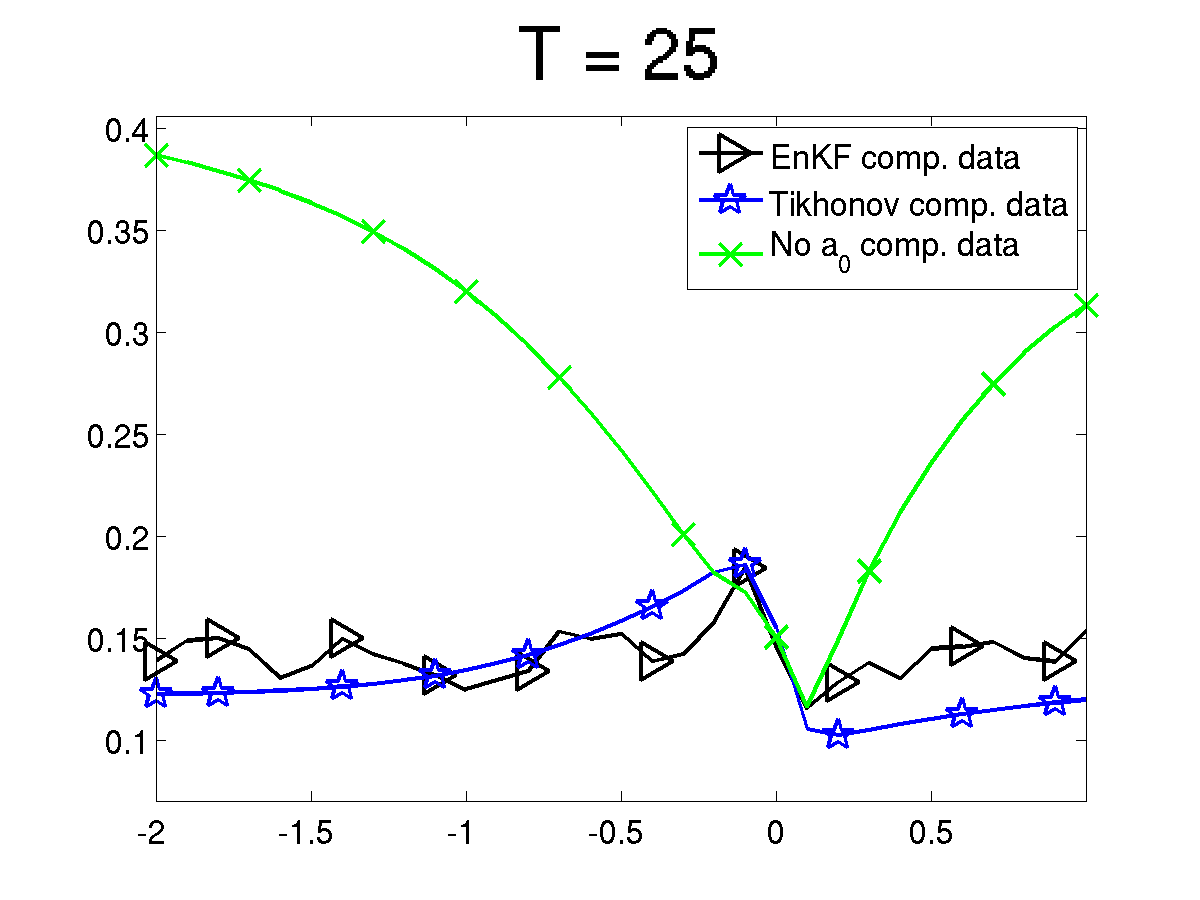}\hfill
\includegraphics[width=.25\textwidth]{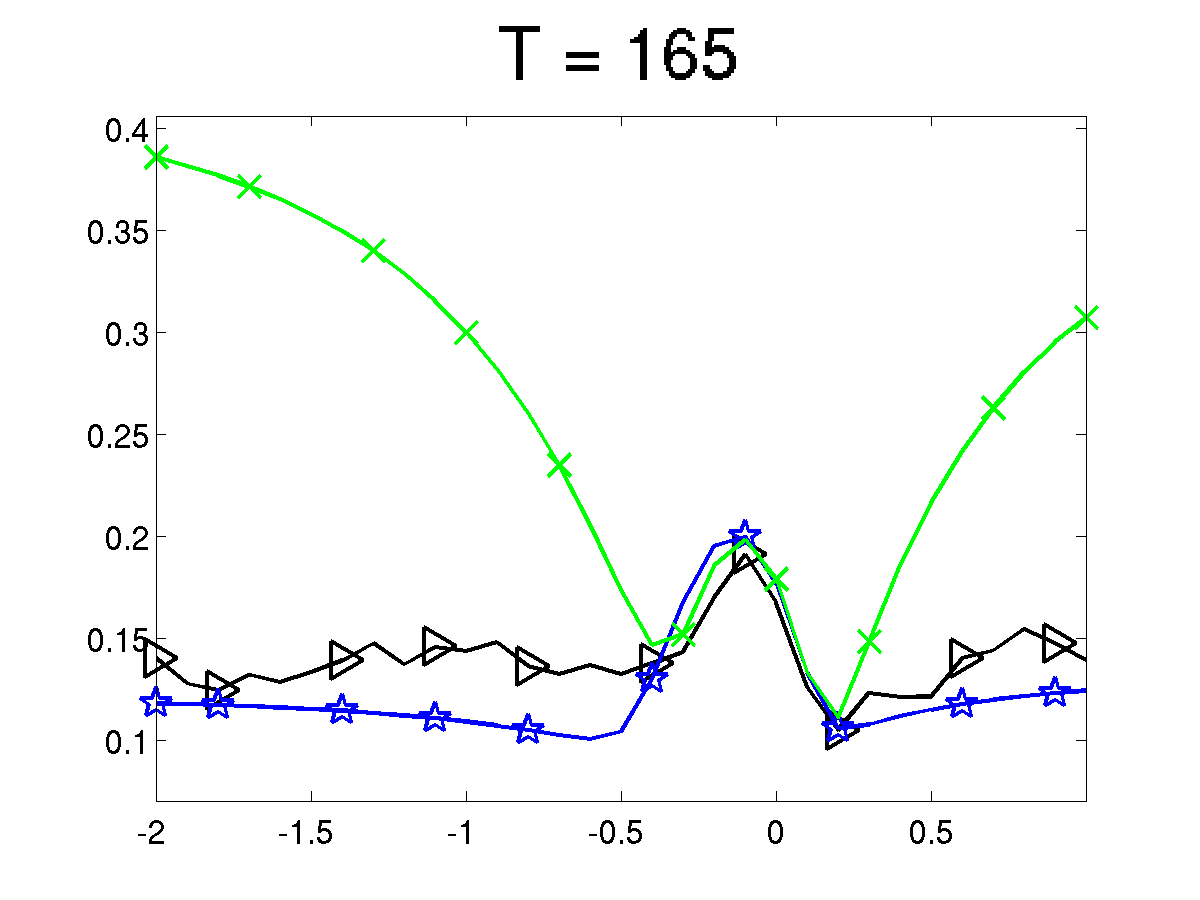}\hfill
\includegraphics[width=.25\textwidth]{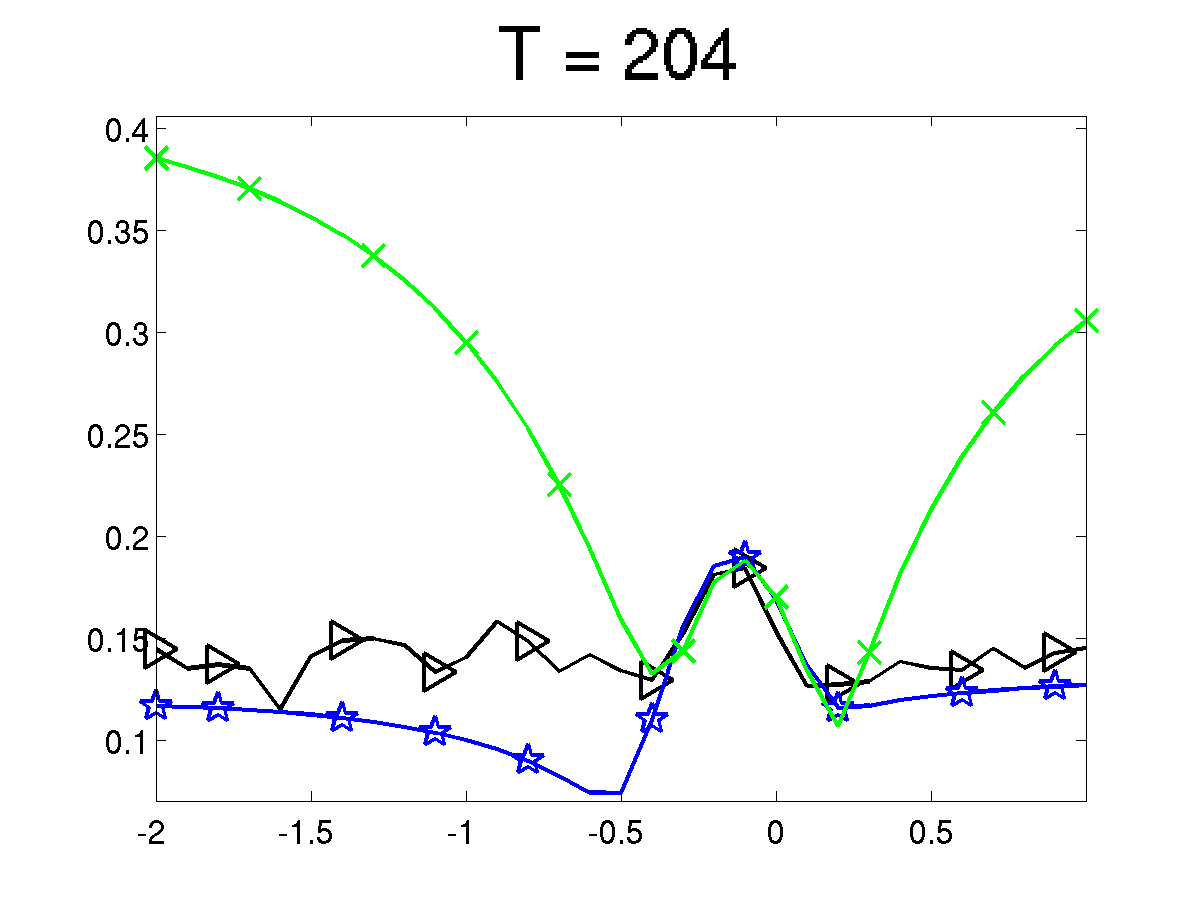}\hfill
\includegraphics[width=.25\textwidth]{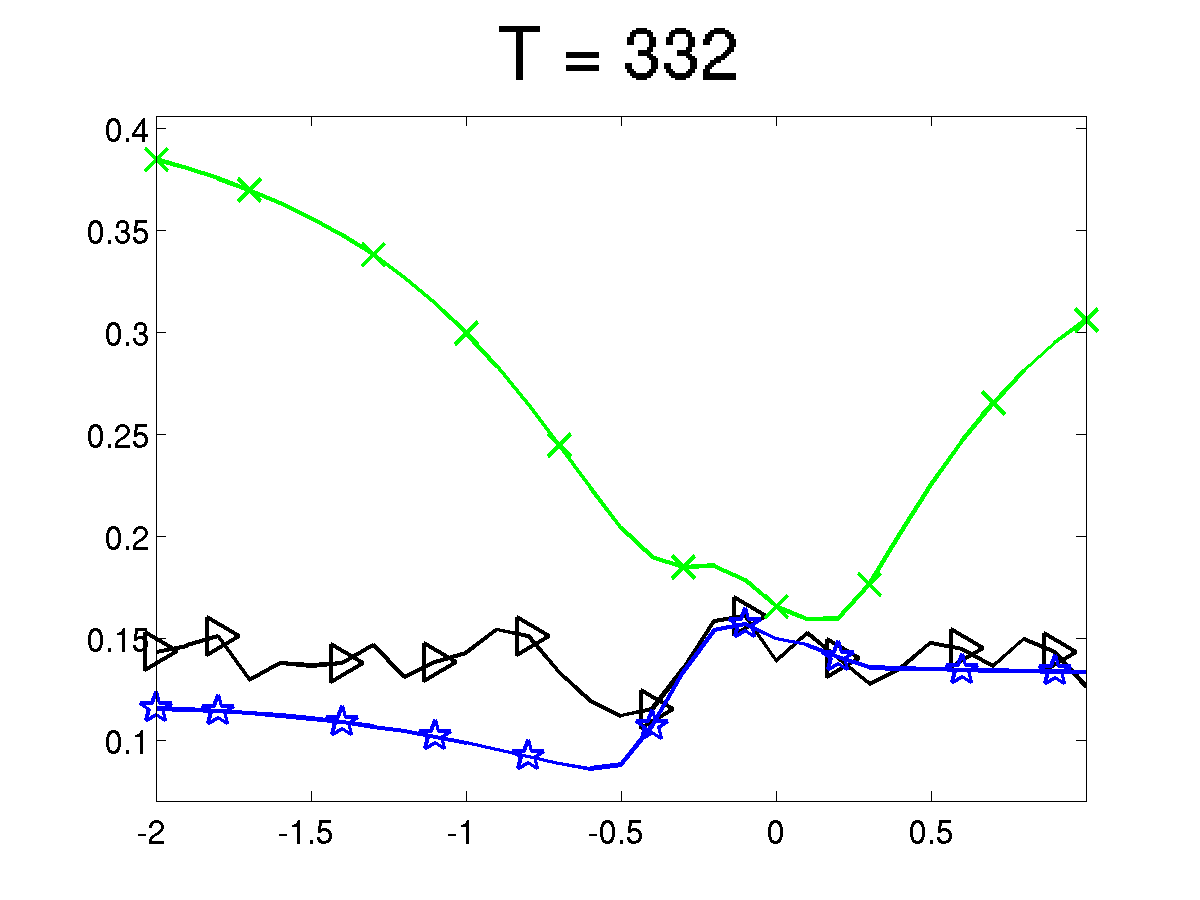}
\caption{\label{spx_c} Reconstructed SPX local volatility surfaces at different maturities obtained with Tikhonov-type and EnKF 
methods using completed data.}
\end{center}
\end{figure}

In Figure~\ref{spx_c} we consider SPX local volatility surfaces at different maturities obtained with the methods
of Sections~\ref{sec:scarce} and~\ref{sec:enkf}. 
They were computed using completed data as in~\cite{kahale}.
Note that the three method variants produce similar results around $y = 0$. (This is called {\em at-the-money}.) 
For $|y| > 0.5$ (the so-called {\em in-the-money} and {\em out-of-the-money} regions), 
if we do not add the $a_0$ penalty, the two wings blow up. 

\begin{figure}[H]
 \begin{center}
 \includegraphics[width=.25\textwidth]{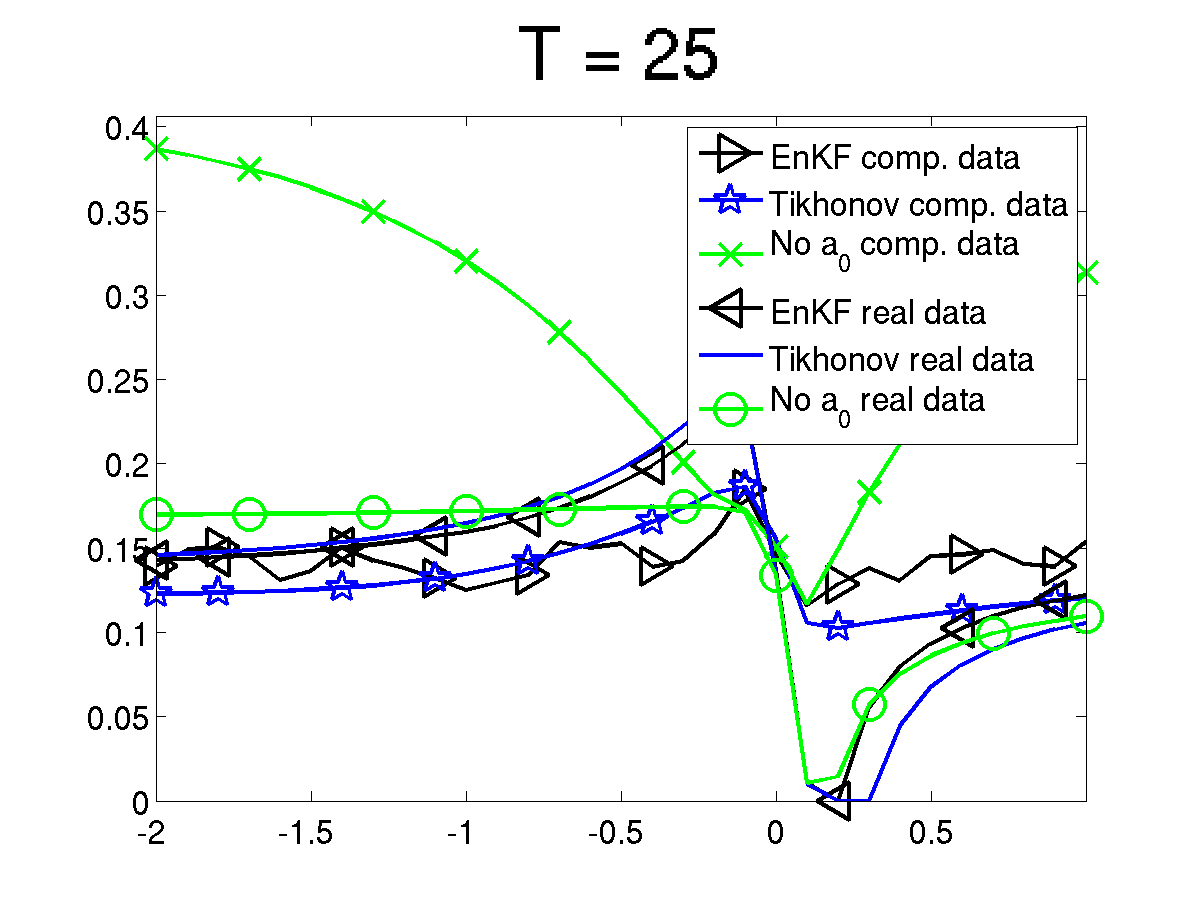}\hfill
 \includegraphics[width=.25\textwidth]{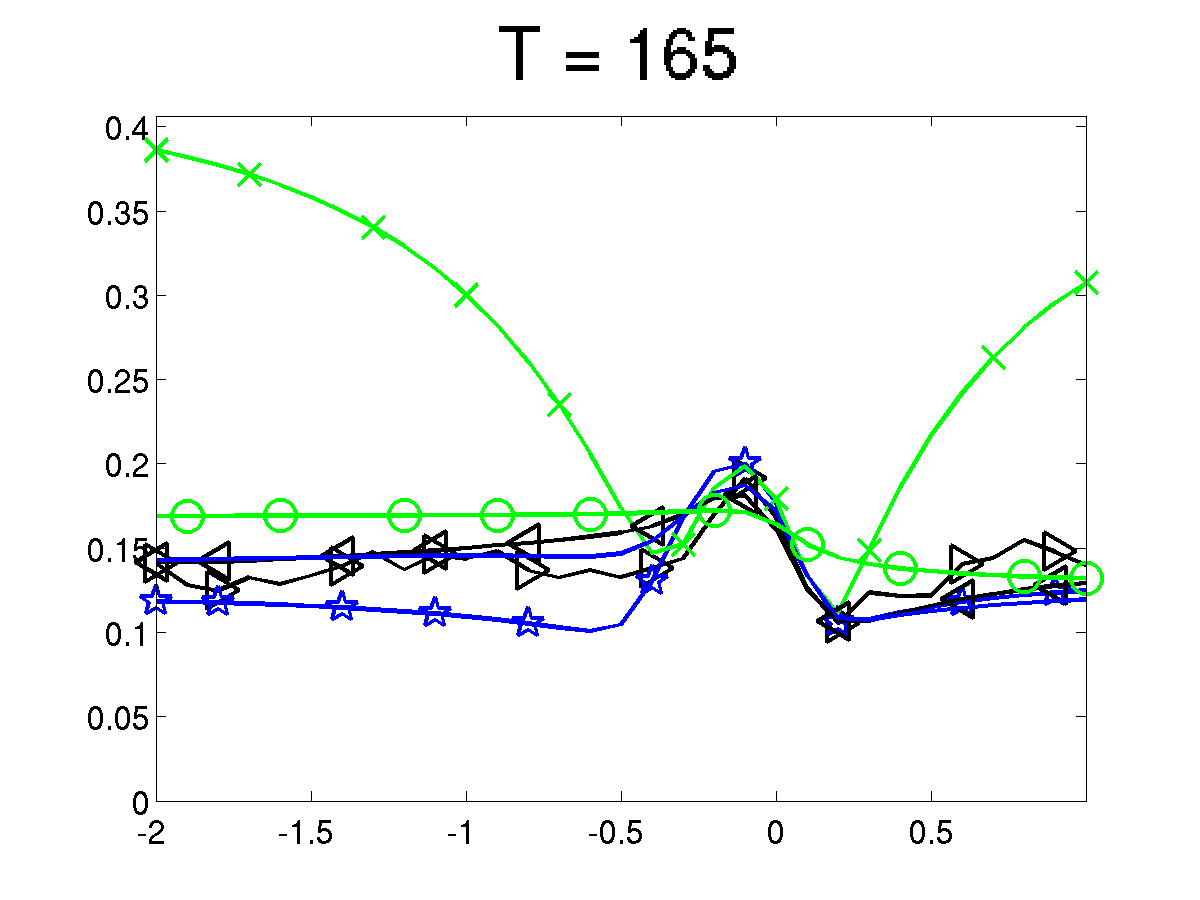}\hfill
 \includegraphics[width=.25\textwidth]{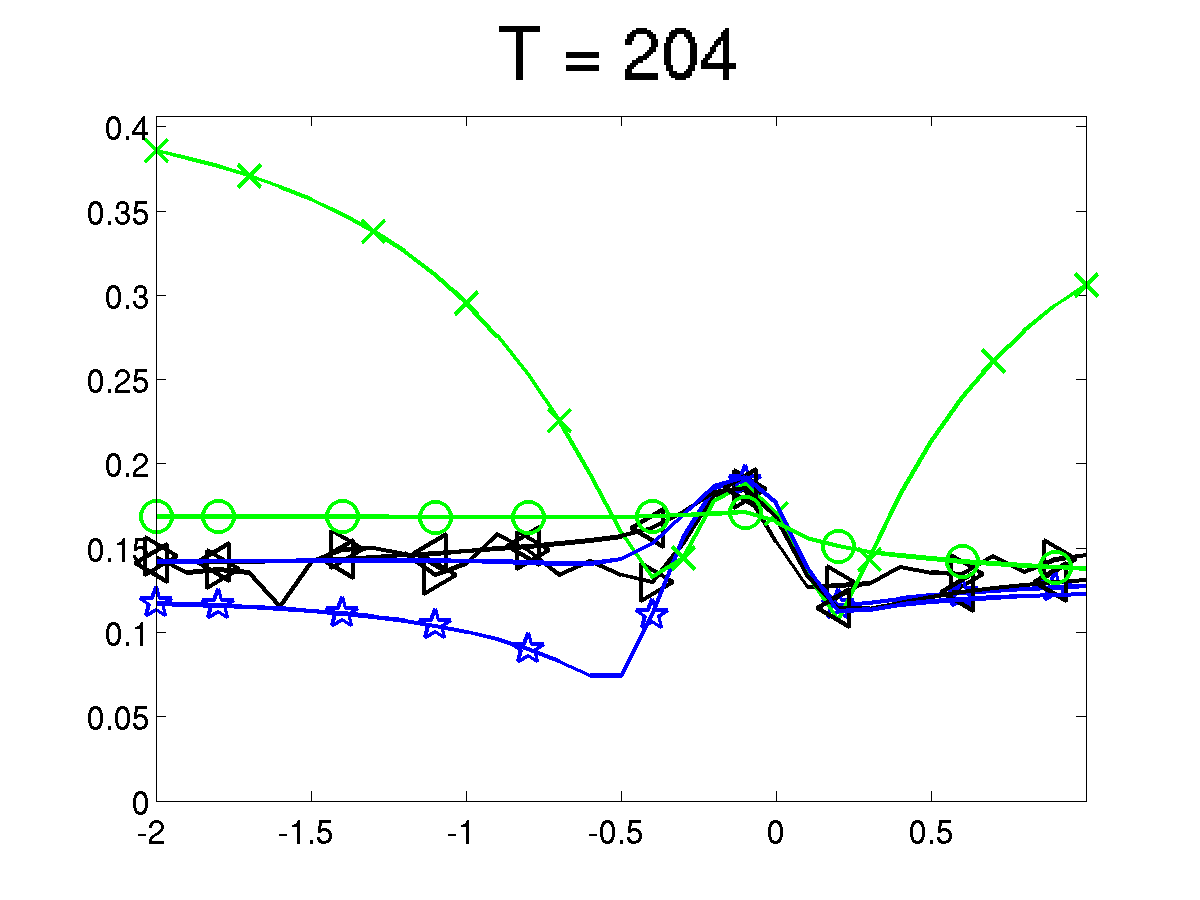}\hfill
 \includegraphics[width=.25\textwidth]{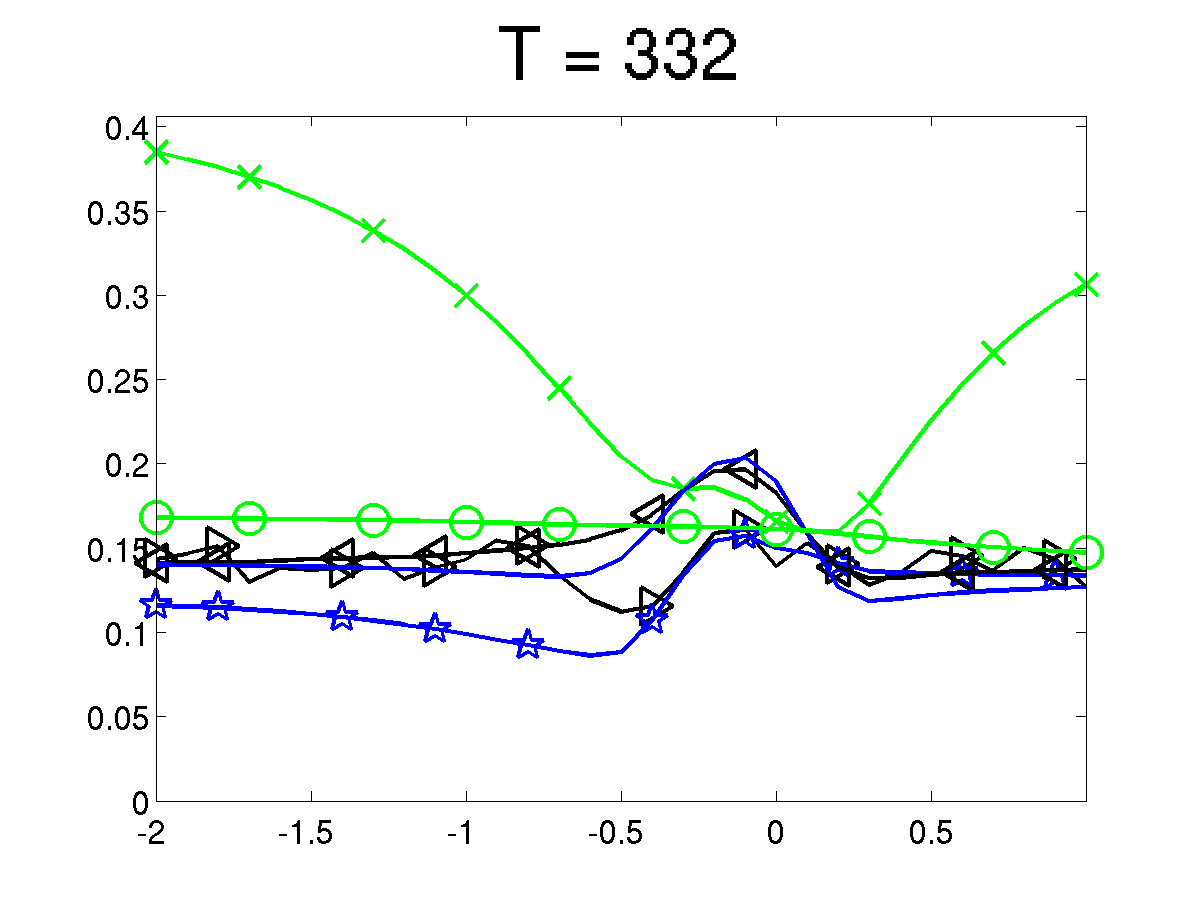}\
 \caption{\label{spx6} Reconstructed SPX local volatility surfaces obtained with six method variants.}
 \end{center}
 \end{figure}

Figure~\ref{spx6} presents reconstructed SPX local volatility surfaces obtained with all six method
variants. When put together, we can see that the different methods lead to results that 
get closer as the maturity gets longer around $y = 0$. 

\begin{figure}[H]
\begin{center}
\includegraphics[width=.25\textwidth]{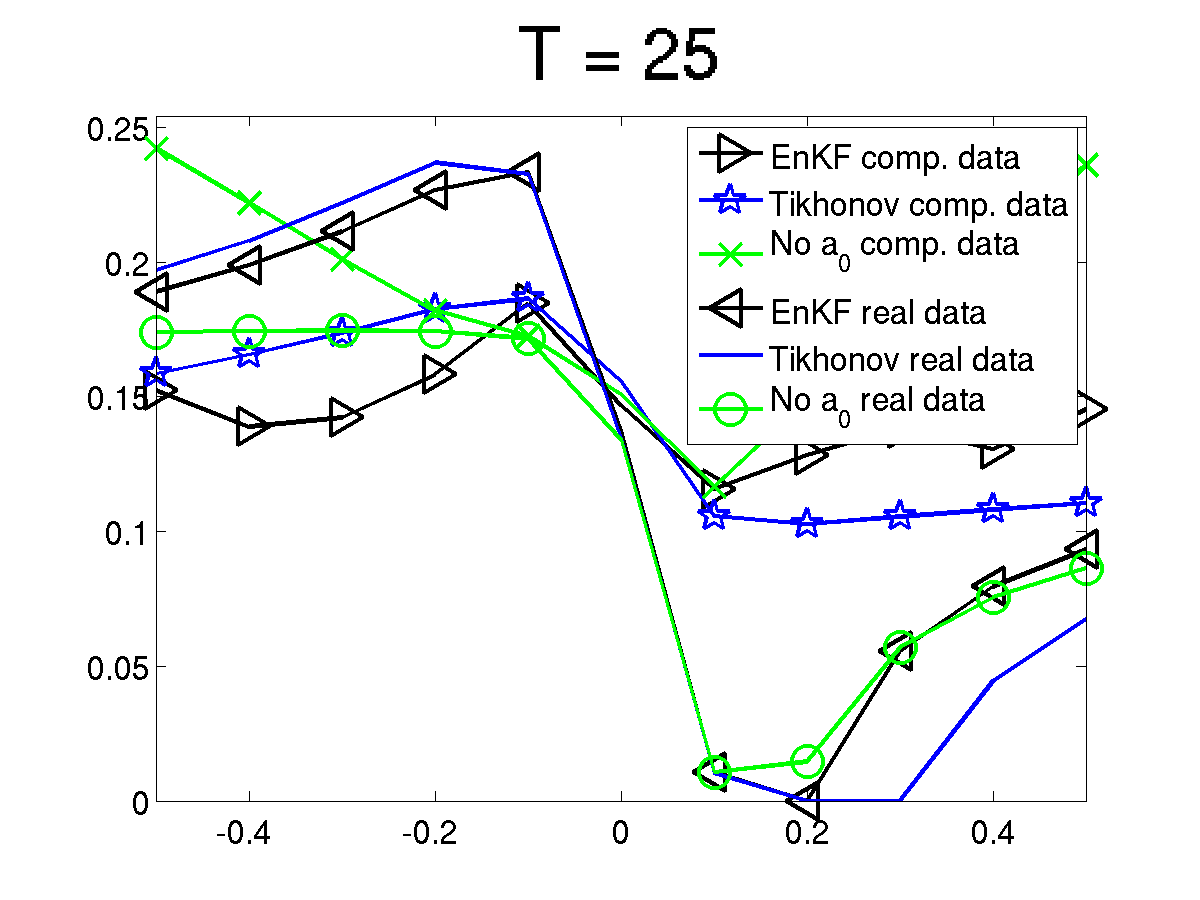}\hfill
\includegraphics[width=.25\textwidth]{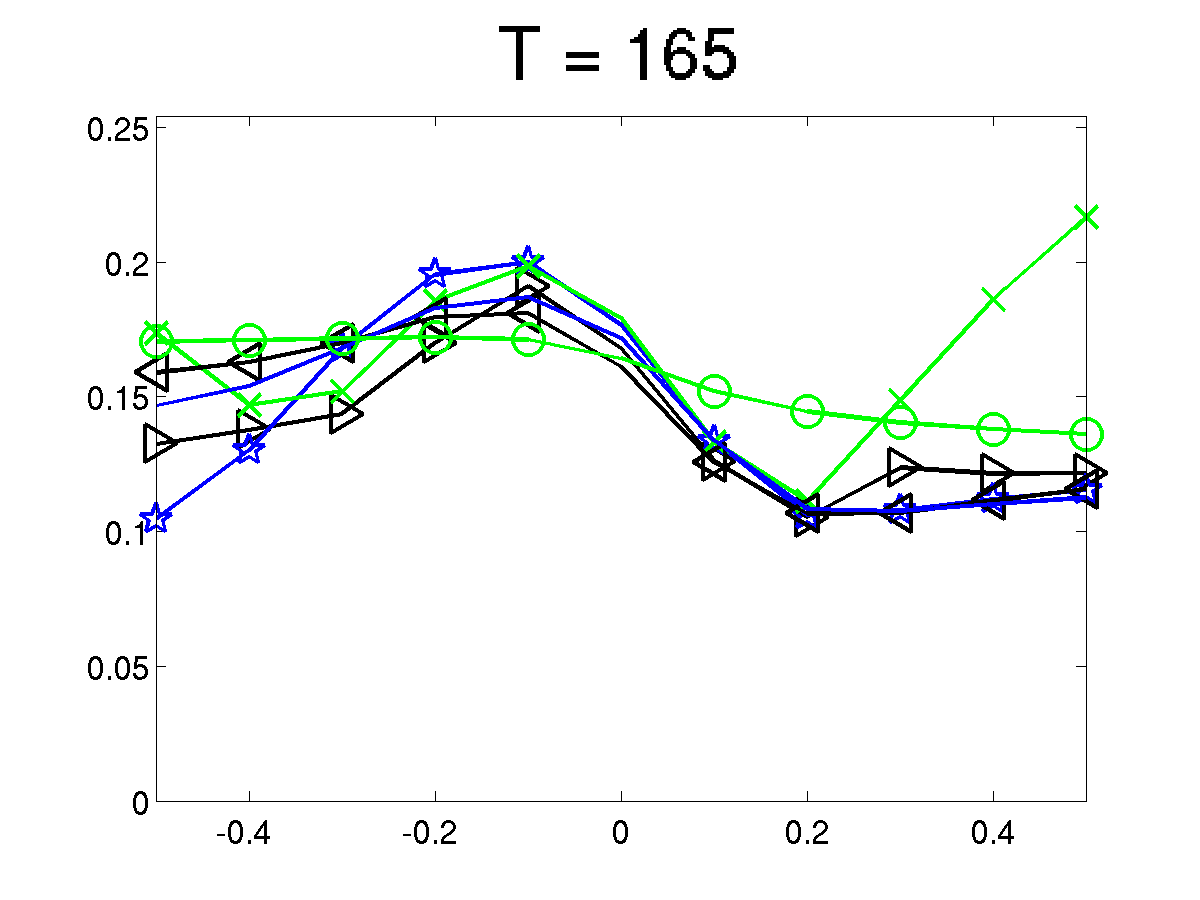}\hfill
\includegraphics[width=.25\textwidth]{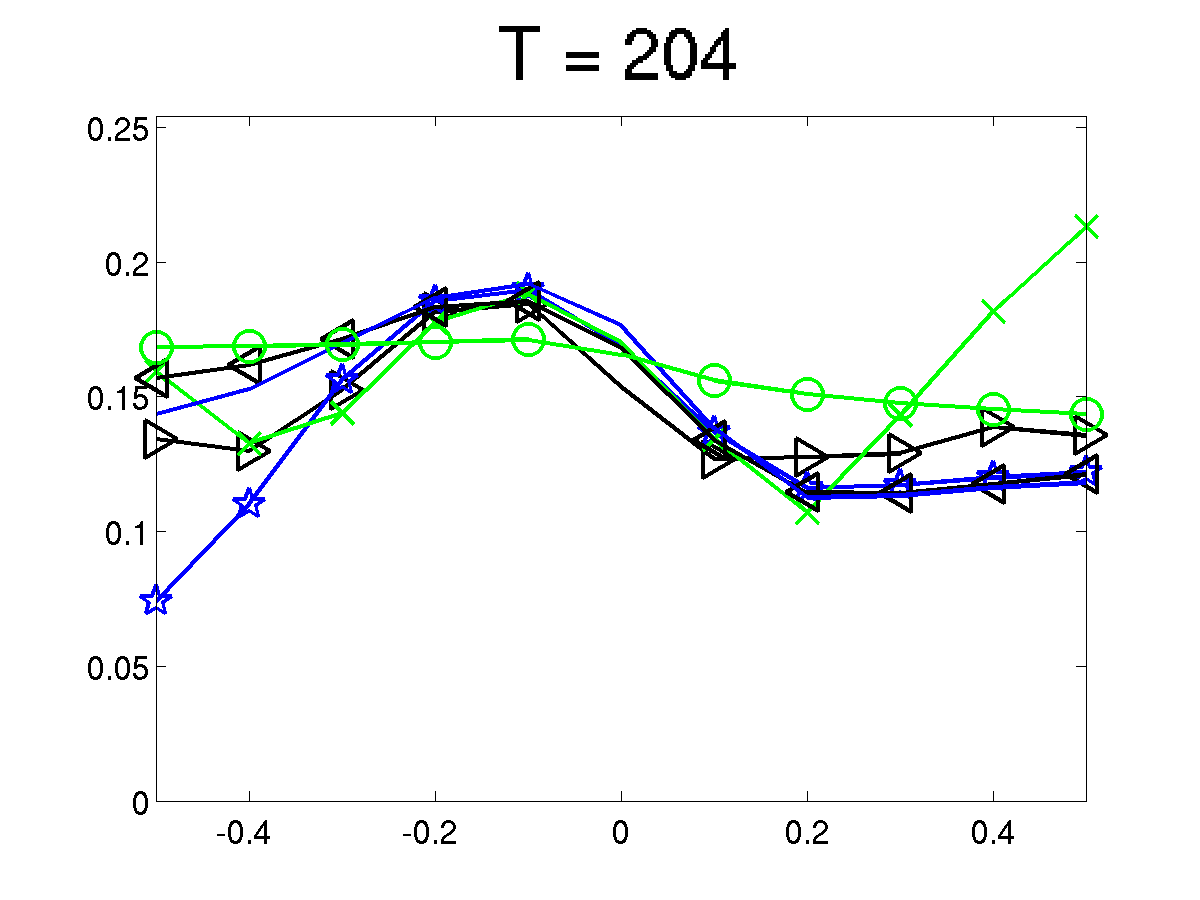}\hfill
\includegraphics[width=.25\textwidth]{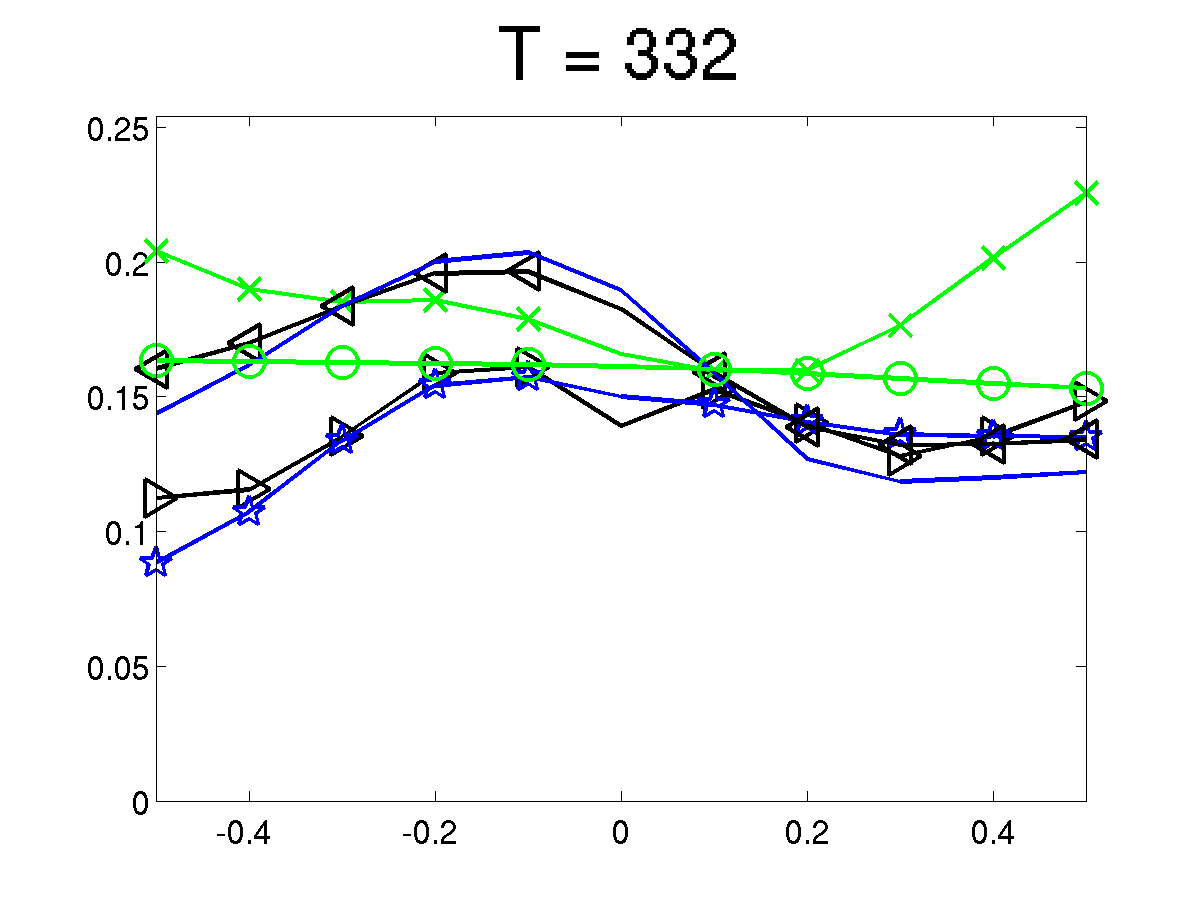}
\caption{\label{spx6a} Reconstructed SPX local volatility surfaces obtained with six method variants 
for different maturities in the at-the-money ($y=0$) neighborhood.}
\end{center}
\end{figure}

Figure~\ref{spx6a} is a zoom-in of Figure \ref{spx6} to the at-the-money region.  
Observe that the plotted curves  are generally divided into two groups: 
one is obtained from the original (scarce) data and the other from the completed 
data. This phenomenon is  clearer in the figures for the earliest and latest dates $T=25$ and $T=332$.   

\begin{figure}[H]
\begin{center}
\includegraphics[width=.25\textwidth]{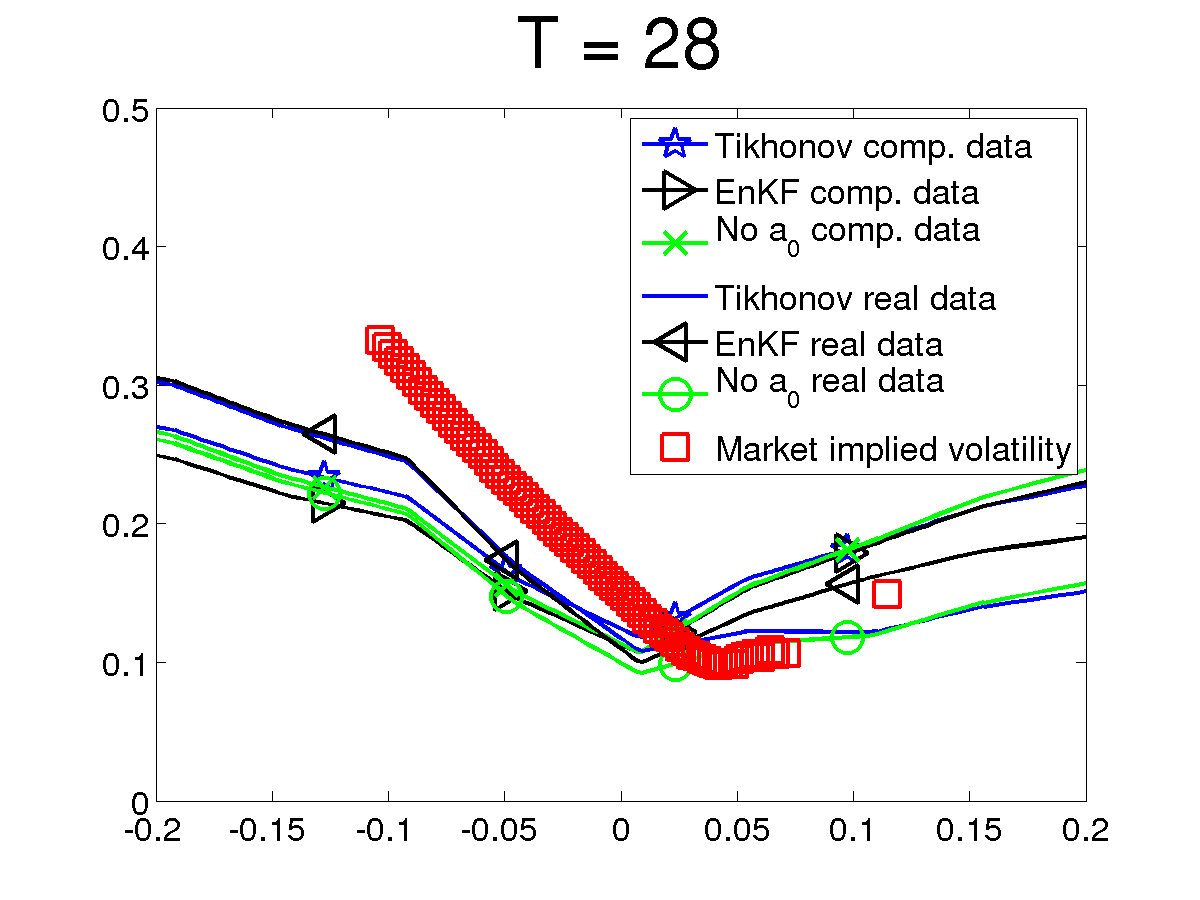}\hfill
\includegraphics[width=.25\textwidth]{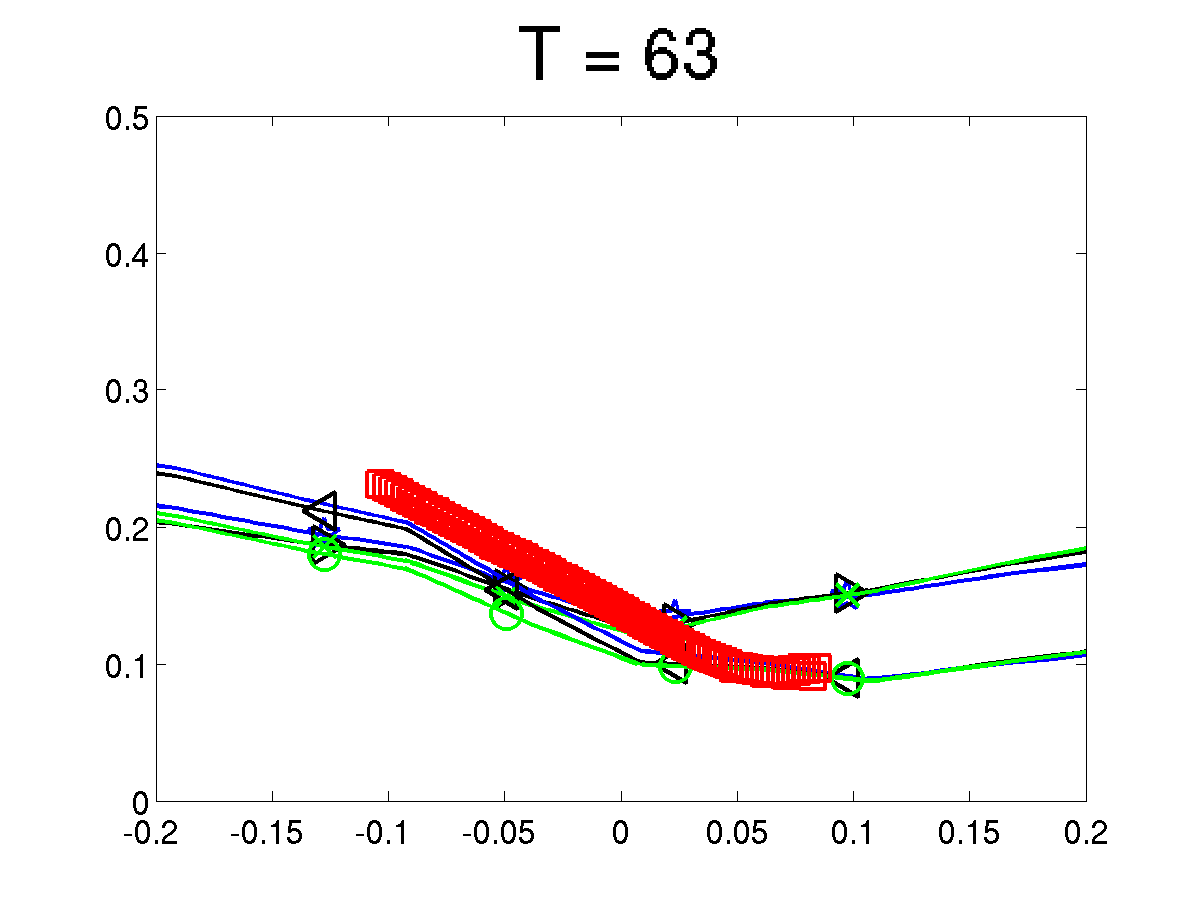}\hfill
\includegraphics[width=.25\textwidth]{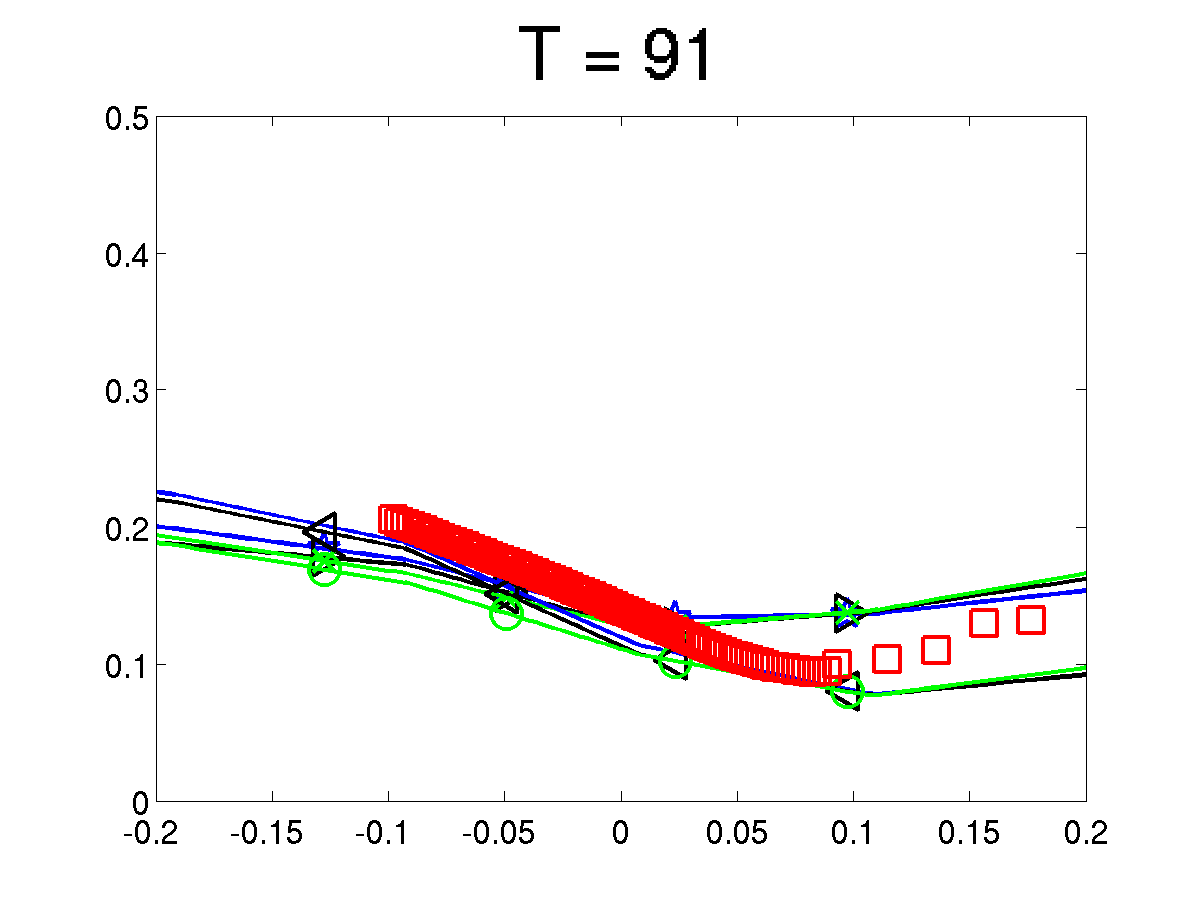}\hfill
\includegraphics[width=.25\textwidth]{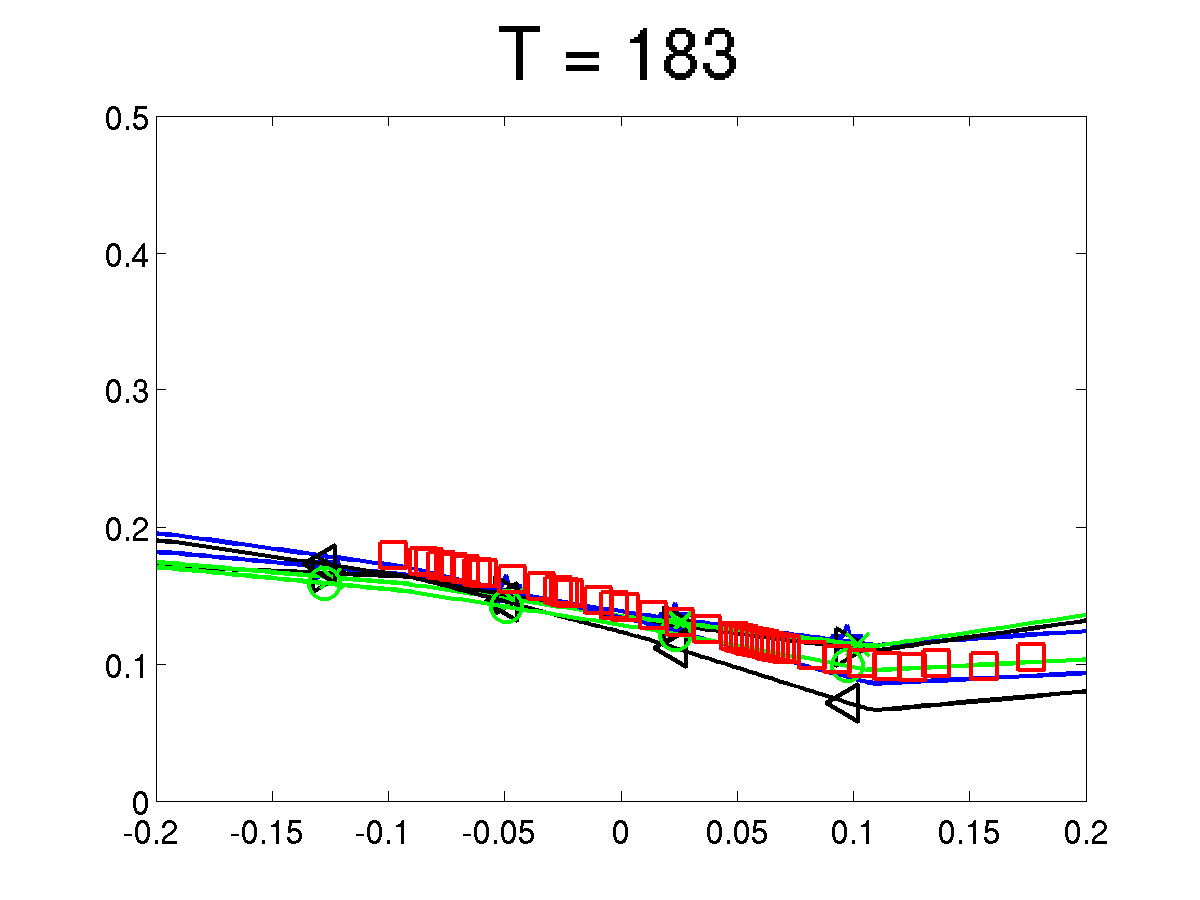}
\includegraphics[width=.25\textwidth]{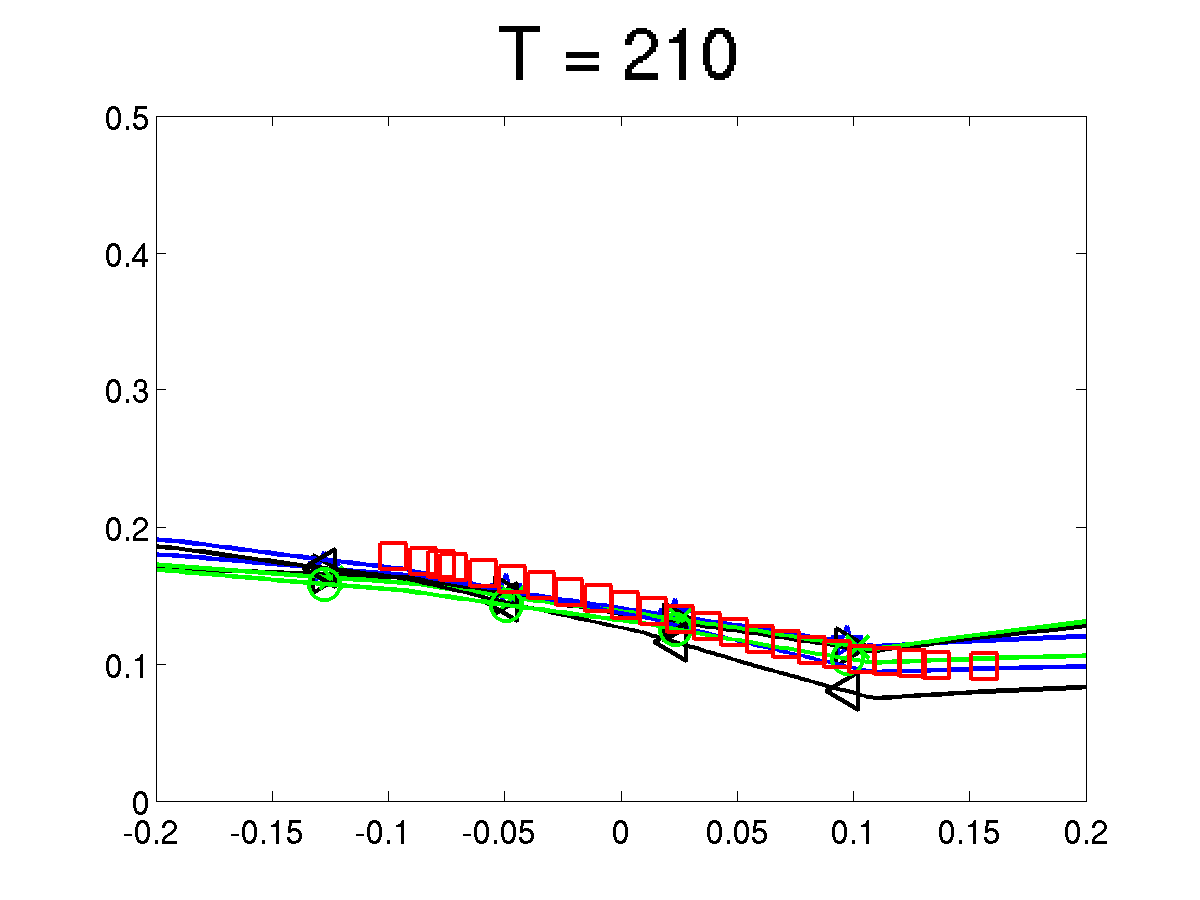}
\includegraphics[width=.25\textwidth]{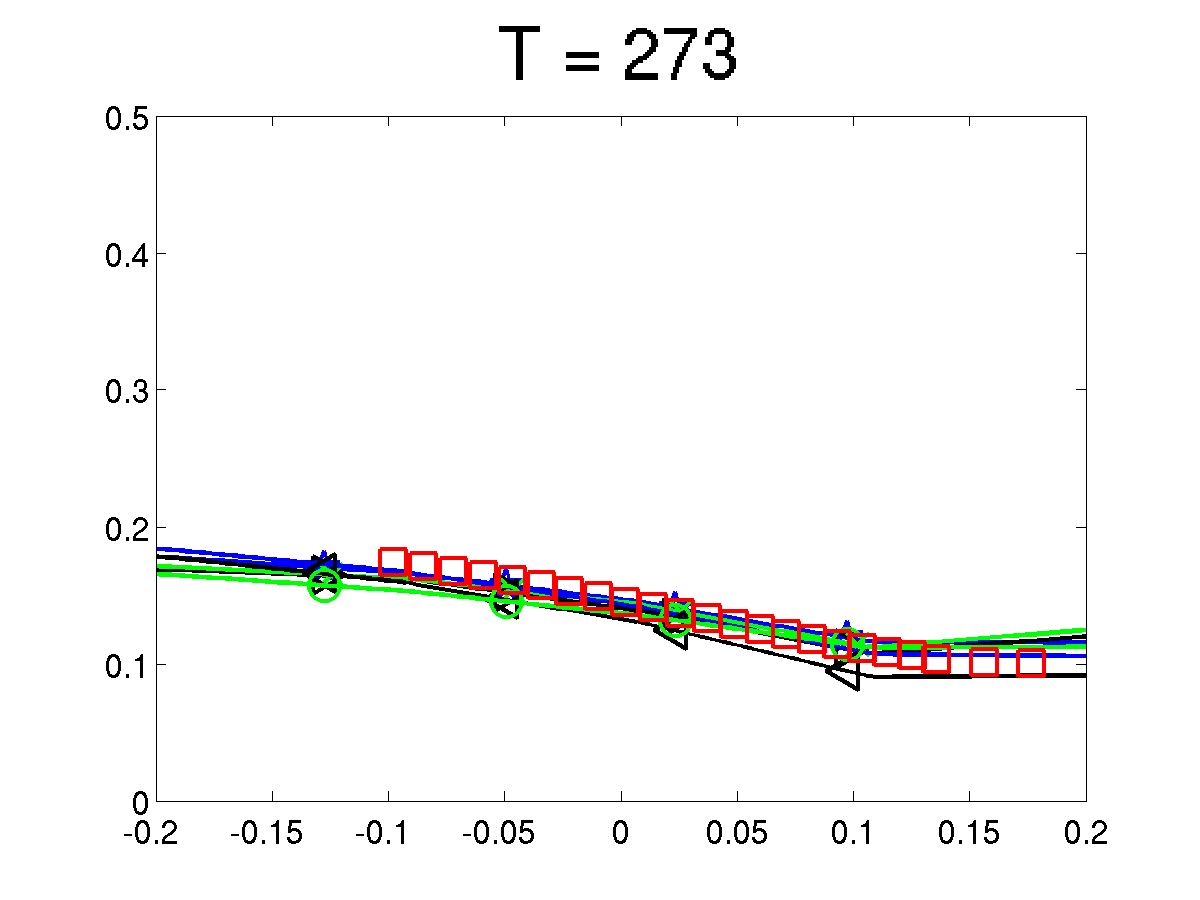}
\includegraphics[width=.25\textwidth]{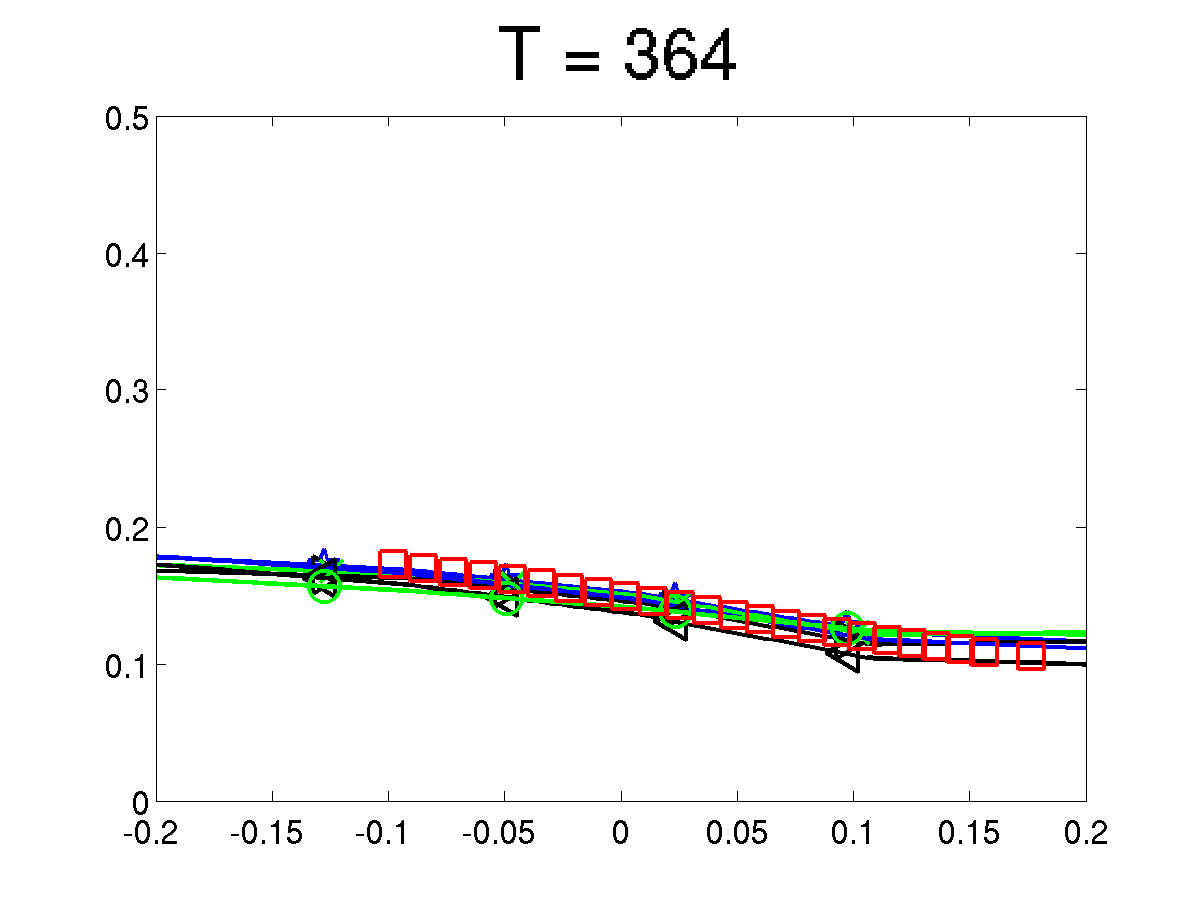}
\includegraphics[width=.25\textwidth]{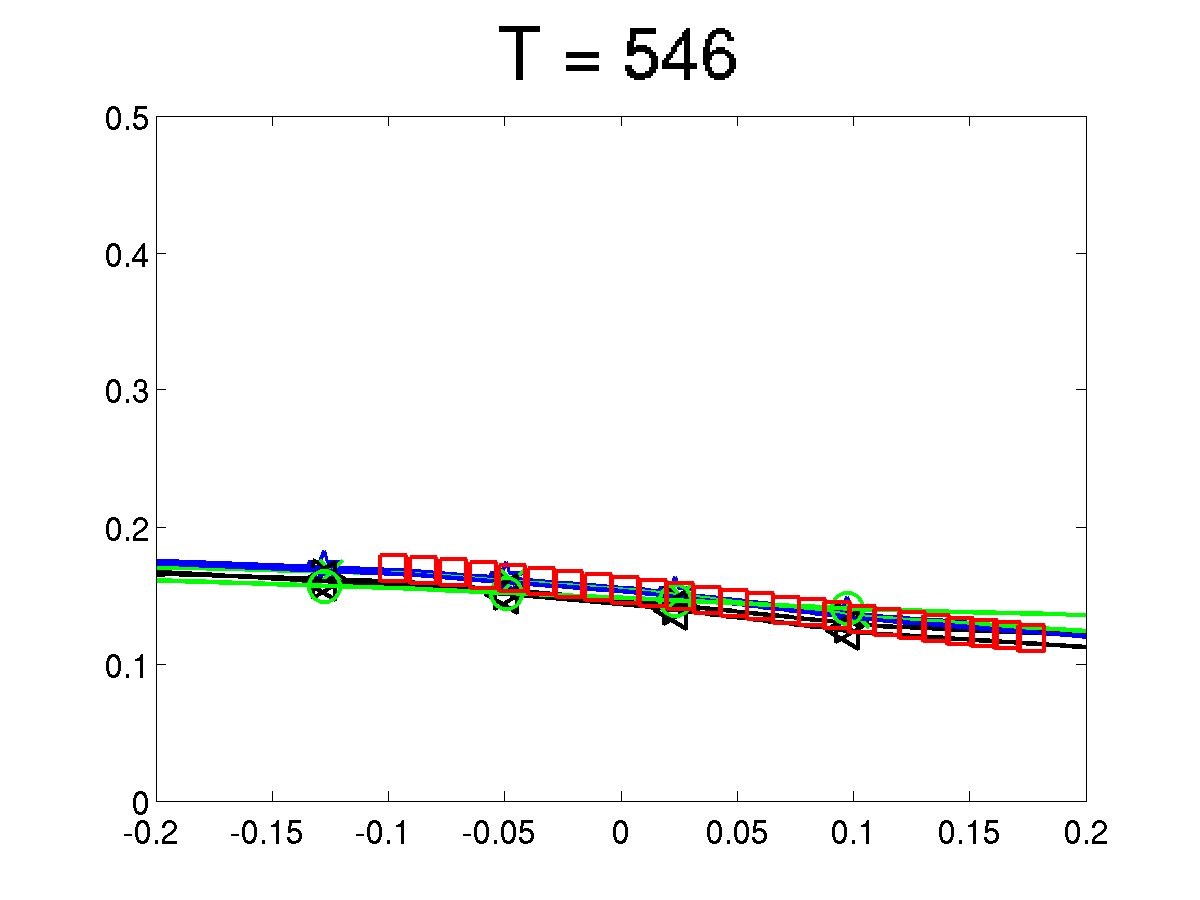}
\includegraphics[width=.25\textwidth]{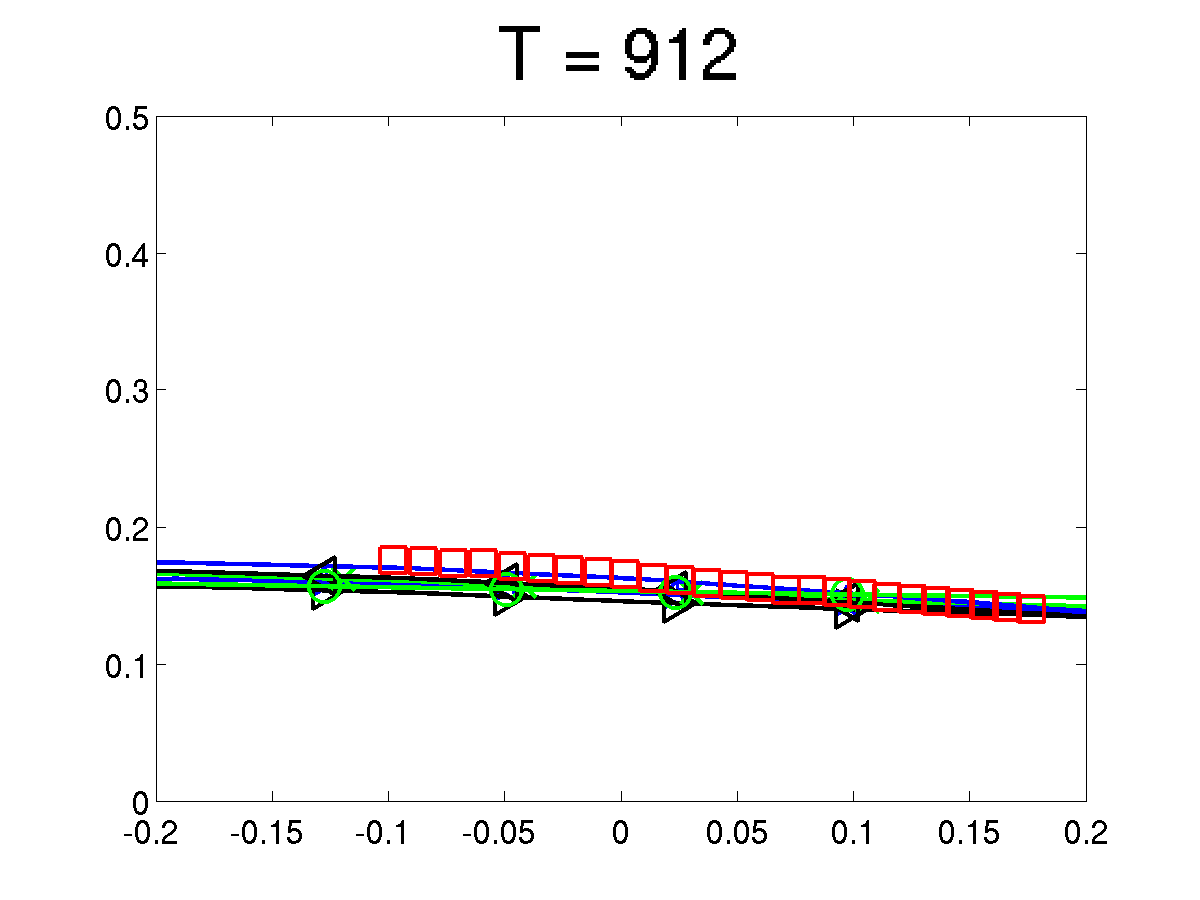}
\caption{\label{spx_im1} Implied (Black-Scholes) volatility corresponding to the local volatility surfaces obtained with the six 
method variants compared to the market one.  }
\end{center}
\end{figure}

Let us define {\em implied volatility} as the volatility 
that would be observed for a standard call (or put) contract to give the observed price if the classical Black-Scholes
formula were used.
This concept is prevalent in market related discussions.  
In Figure~\ref{spx_im1} we see that  
the reconstructions of implied volatility from local volatility do better for intermediate term  and long term maturities. However, when the maturities are short, the results split into two groups according to whether we use the original scarce data or the completed data.   

We measure the data misfit of all methods in three different ways. 
Let $I^L_{i,j}$ ($I^{ba}_{i,j}$) denote the implied volatility corresponding to the reconstructed local volatility 
(from the average of bid and ask option prices) with strike $K_i$ and maturity $\tau_j$. 
In the SPX example, we restrict the strikes to be between 1890 and 2500. 
Further, $V_{i,j}$ is the volume of the option price with strike $K_i$ and maturity $\tau_j$, 
and $N_{vol}$ is the number of contracts that have a nonzero volume. 
We define
\begin{eqnarray*}
RMSE &=& \sqrt{\sum_{i,j} (I^{L}_{i,j}-I^{ba}_{i,j})^2/N_{vol}}, \quad{\rm (root~mean~square~error)} \\
RWMSE &=& \sqrt{\sum_{i,j} (I^{L}_{i,j}-I^{ba}_{i,j})^2\times V_{i,j}/N_{vol}}, \quad{\rm (root~weighted~mean~square~error)} \\ 
RR&=&\sqrt{\sum_{i,j}(I_{i,j}^L-I_{i,j}^{ba})^2}/\sqrt{\sum_{i,j}(I_{i,j}^{ba})^2}, \quad{\rm (relative~residual)}.
\end{eqnarray*}
The resulting values are presented in Table~\ref{tabspx4}.
\begin{table}[htb]
\centering
\caption{{\small Measures of data misfit of the 6 models.
}}\label{tabspx4}
\begin{tabular}{c|c|c|c|c|c|c|}
\cline{2-7}
 &\multicolumn{4}{c|}{Tikhonov-type} & \multicolumn{2}{c|}{EnKF}\\
\cline{2-7}
 &Scarce&Comp.&Scarce (no $a_0$)& Comp. (no $a_0$)&Scarce&Comp.\\
\hline
\multicolumn{1}{|l|}{RMSE} & 0.0195 & 0.0321 & 0.0290 & 0.0325 & 0.0255 & 0.0324 \\
\multicolumn{1}{|l|}{RWMSE} & 0.0175 & 0.0241 & 0.0252 & 0.0242 & 0.0241 & 0.0242\\
\multicolumn{1}{|l|}{RR} & 0.1407 & 0.1987 & 0.2292 & 0.2186 & 0.1766 & 0.2186\\
\hline
\end{tabular}
\end{table}

\paragraph{Discussion of the real equity data results}

From the above experiment we conclude that using real or completed data sets we get 
quite different results. 
Within the completed data set, if we discard the $a_0$ penalty, 
the two wings of the local volatility surface are not stable,  
both for methods of Sections~\ref{sec:scarce} and~\ref{sec:enkf}. 
This is apparent in all the figures as well as Table~\ref{residual_6}. 
From the results involving reconstruction of the implied volatility, 
the Tikhonov-type method using the original (scarce) data has the best residuals in all three measures, 
and the EnKF method is the second best.

\subsection{Results for commodities and data completion}
\label{sec:datacompcom}
Commodities have been traded extensively in different markets throughout the world for centuries. In many of those markets, 
 a number of liquid options on such assets are also traded.  Here again, data from such markets are freely available,
 and modelling such data is very relevant for financial analysts and risk management applications.

In the present set of examples, we consider the adaption of Dupire's model to the context of option prices 
on commodity futures introduced in~\cite{AlbAscZub2015}. 
For the present purpose, it consists of essentially the model presented in 
Sections~\ref{sec:volatility} and~\ref{sec:scarce}.

We chose data from WTI\footnote{West Texas Intermediate (WTI) is a grade of crude oil used as a benchmark in oil pricing.}
futures and their options, as well as Henry Hub\footnote{Henry Hub (HH) natural gas futures are standardized 
contracts traded on the New York Mercantile Exchange (NYMEX).} contracts. 
The end-of-the-day WTI option and future prices were traded at 06-Sep-2013, 
with the maturity dates, 18-Oct-2013, 16-Nov-2013, 17-Dec-2013, 16-Jan-2014, 15-Feb-2014 and 20-Mar-2014. 
The end-of-the-day Henry Hub option and future prices were traded at 06-Sep-2013, with the maturity dates, 29-Oct-2013, 27-Nov-2013, 27-Dec-2013, 29-Jan-2014, 26-Feb-2014 and 27-Mar-2014. The option prices were divided by their underlying future prices, so $S_0 = 1$.

In what follows, by residual we mean the merit function \eqref{x4} with $\Gamma=\alpha_0^{-1} I$, $\alpha_0>0$. The penalty function~\eqref{x7} was used with $\alpha_1 = \alpha_2 = \alpha_3$.
A gradient descent method was applied in the minimization of the resulting Tikhonov-type regularization functional.
Inspired by the discrepancy principle, different values for $\alpha_0$ were tested, and the largest one which gave a residual below a fixed tolerance was chosen. 
See~\cite{AlbAscZub2015} for further details.

To complete the data (when this alternative was used in this subsection) 
we applied linear interpolation, taking into account the boundary and initial conditions in~\eqref{x21}.
The boundary conditions were applied at $y = \pm 5$. 

When using completed data, we evaluated the local volatility only in the maturity times given in the original data, 
and interpolated it linearly in time, at each step in the minimization.
Also, whenever $|y| > 0.5$ was encountered, we set $a(\tau,y) = a(\tau,0.5*{\rm sign}(y))$.
%
When using the given scarce data, at each iteration, for each maturity time $\tau$ in the data set, 
we interpolated $a(\tau,\cdot)$ linearly in the intervals $[-5,y_{\min})$, and  $(y_{\max},5]$, 
assuming that $a(\tau,-5) = \max\{0.08,a(\tau,y_{\min})\}$ 
and $a(\tau,5) = \max\{0.08,a(\tau,y_{\max})\}$, 
where $y_{\min}$ and $y_{\max}$ are the minimum and maximum log-moneyness strikes in the data set corresponding to $\tau$.

\begin{figure}[!ht]
\centering
\begin{minipage}{0.49\textwidth}
   \centering
       \includegraphics[width=0.475\textwidth]{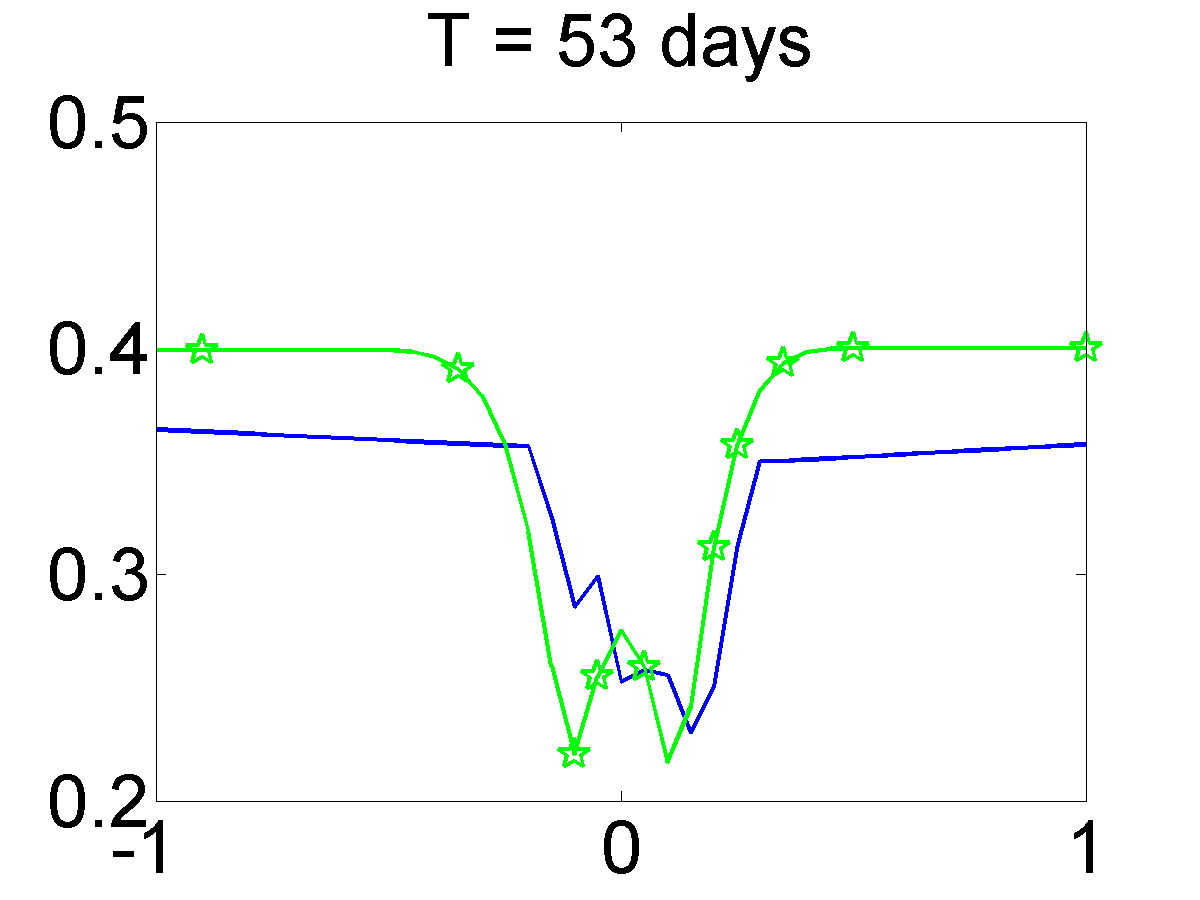}\hfill
       \includegraphics[width=0.475\textwidth]{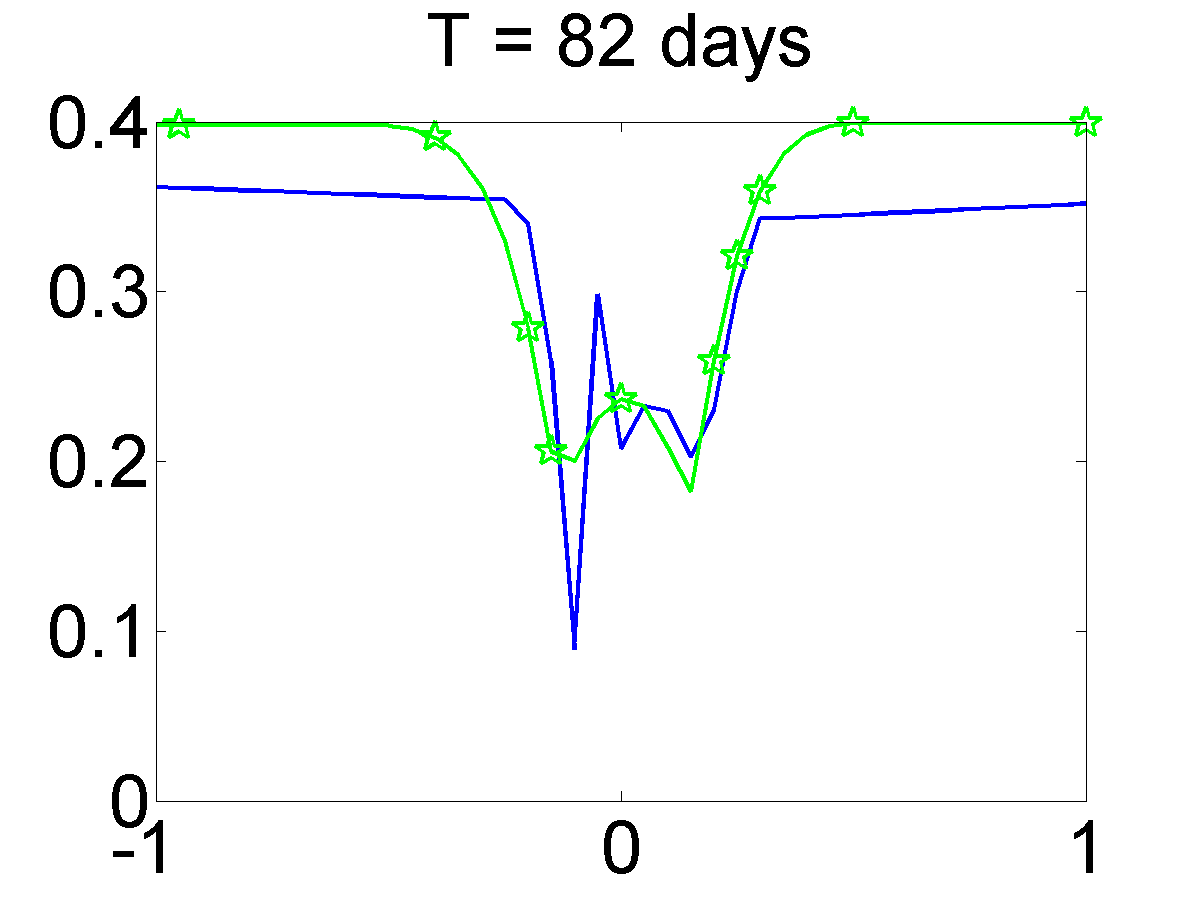}\hfill
       \includegraphics[width=0.475\textwidth]{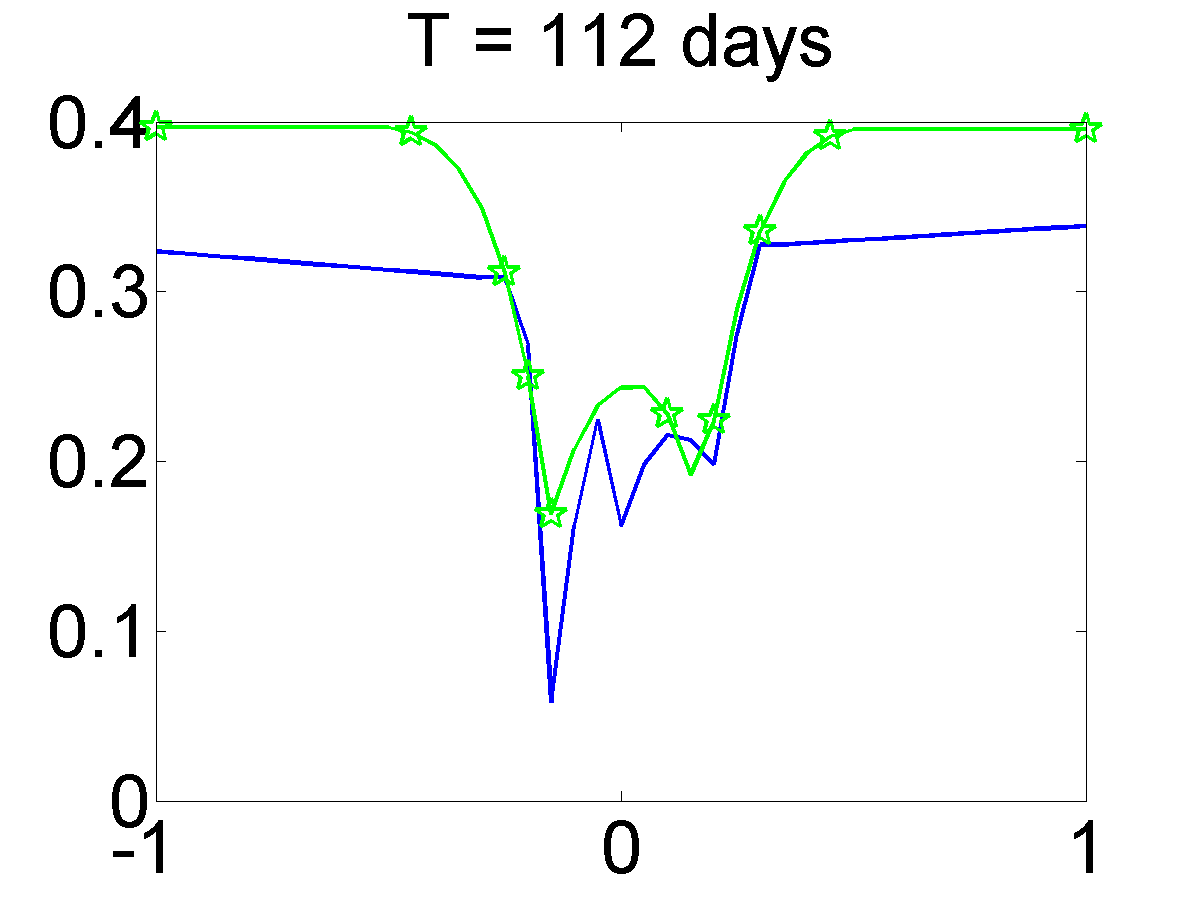}\hfill
       \includegraphics[width=0.475\textwidth]{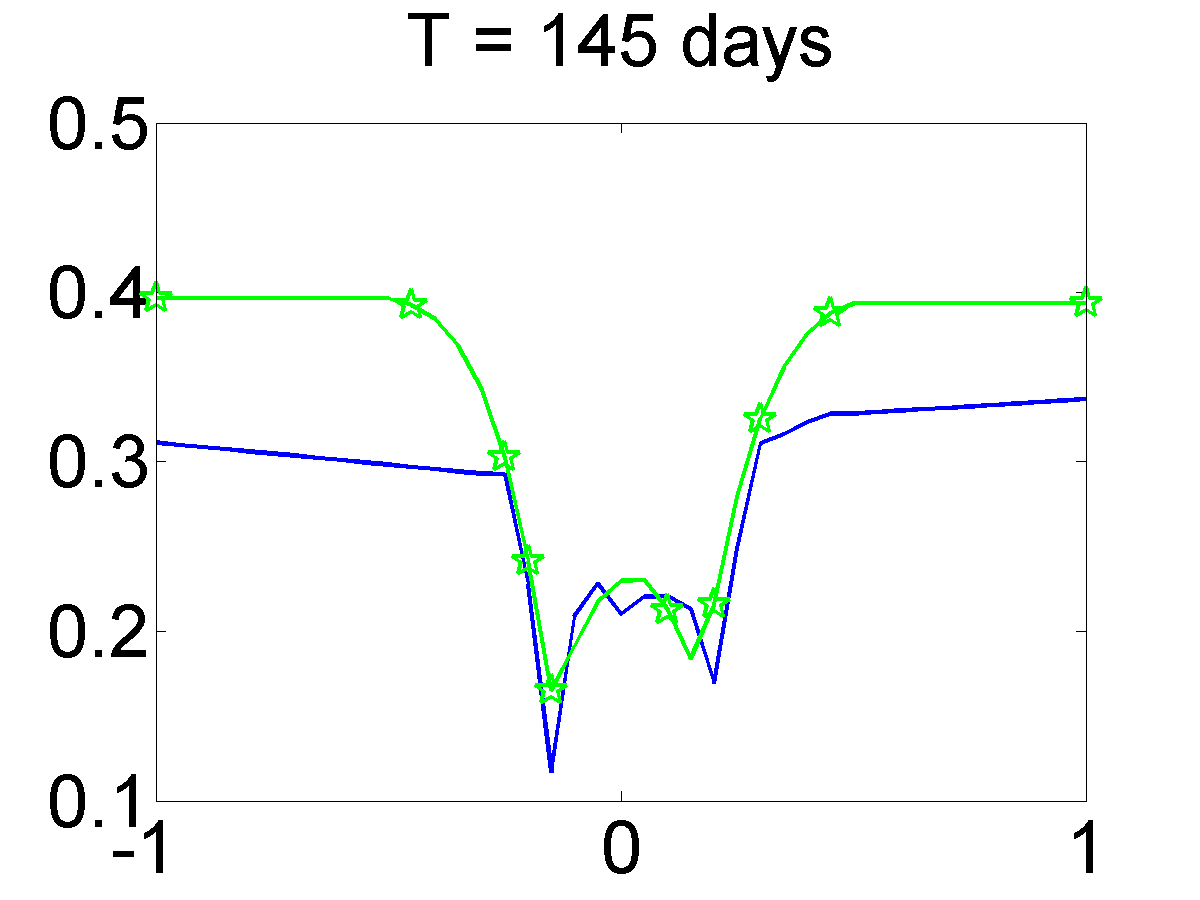}\hfill
       \includegraphics[width=0.475\textwidth]{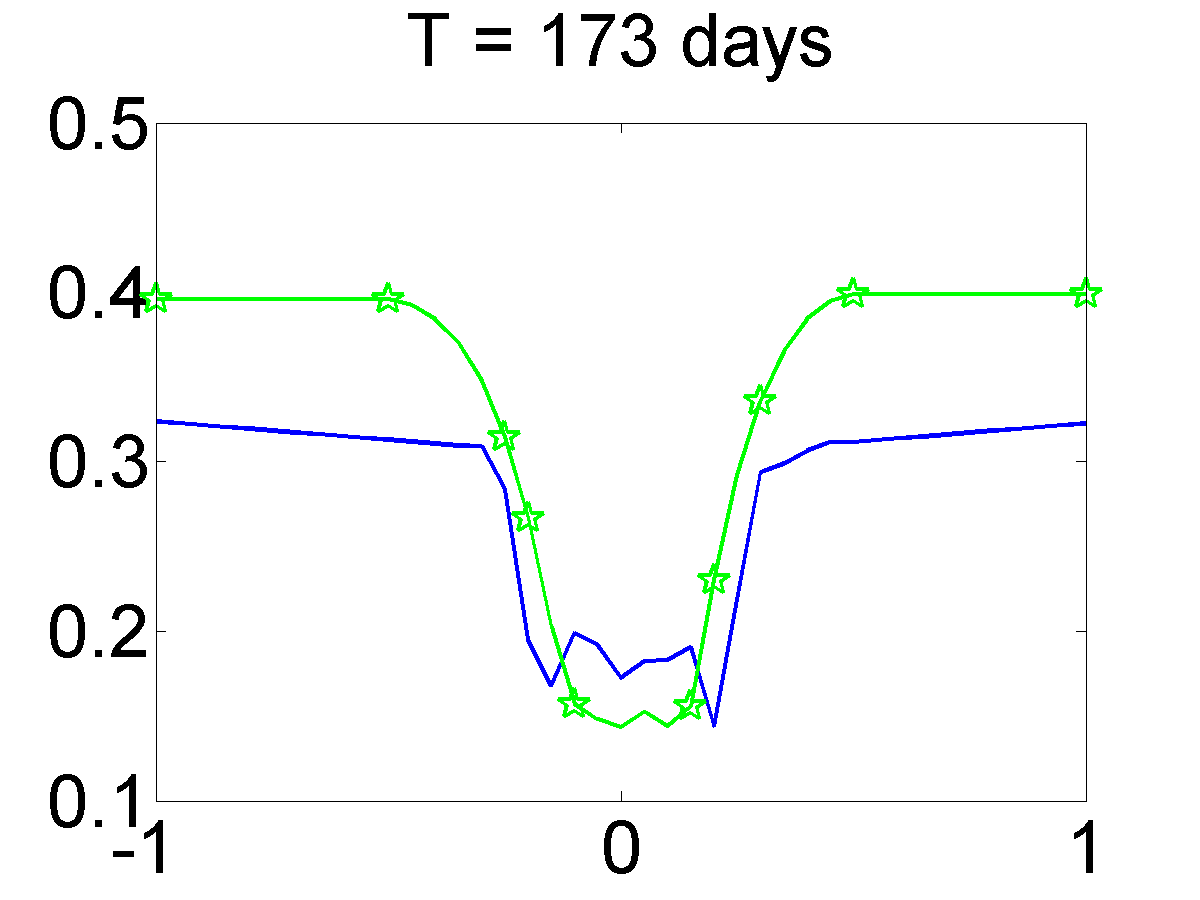}\hfill
       \includegraphics[width=0.475\textwidth]{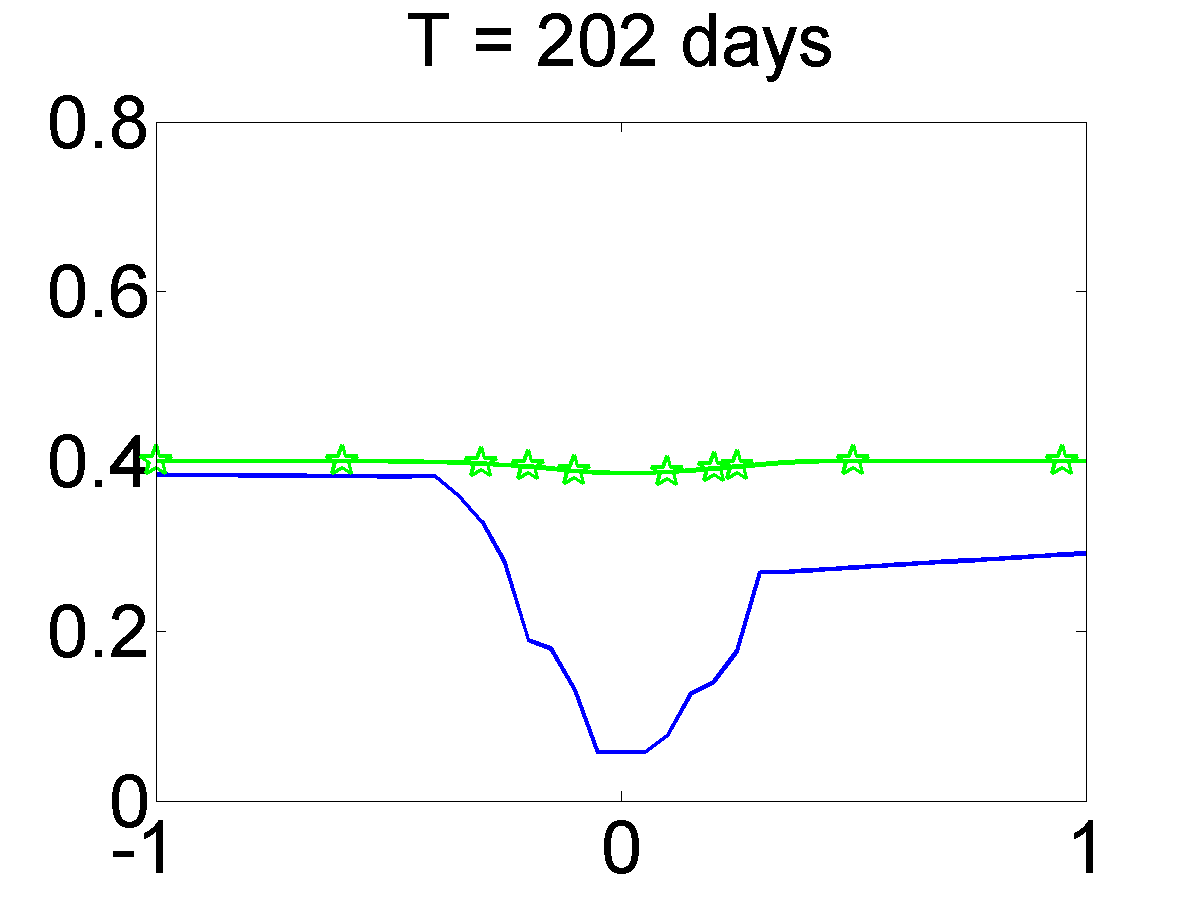}
 \caption{Reconstructed local volatility for different maturity dates for Henry Hub call option prices, 
comparing between completed data (green line with pentagram) and scarce data (blue line) results.}
  \label{fig:com_hh}
\end{minipage}\hfill
\begin{minipage}{0.49\textwidth}
  \centering
   \centering
       \includegraphics[width=0.475\textwidth]{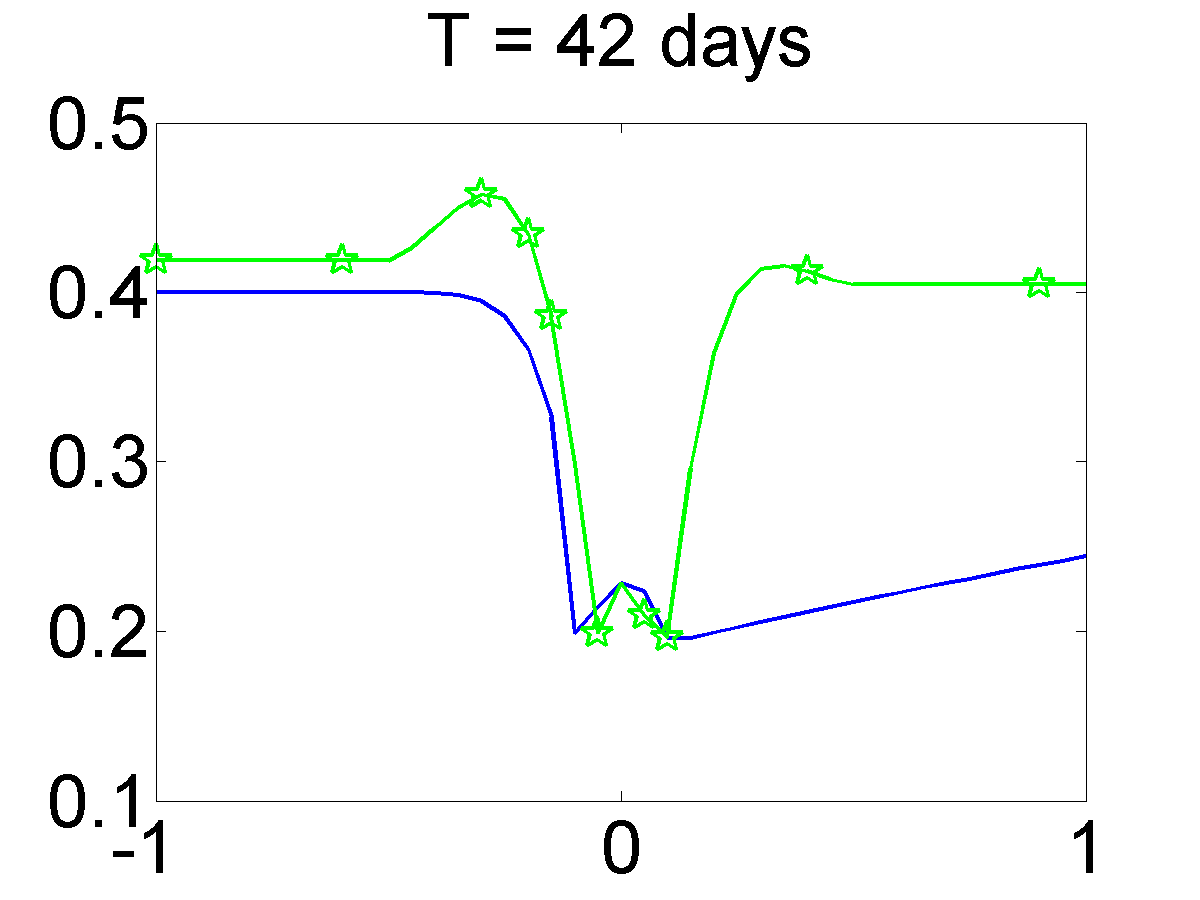}\hfill
       \includegraphics[width=0.475\textwidth]{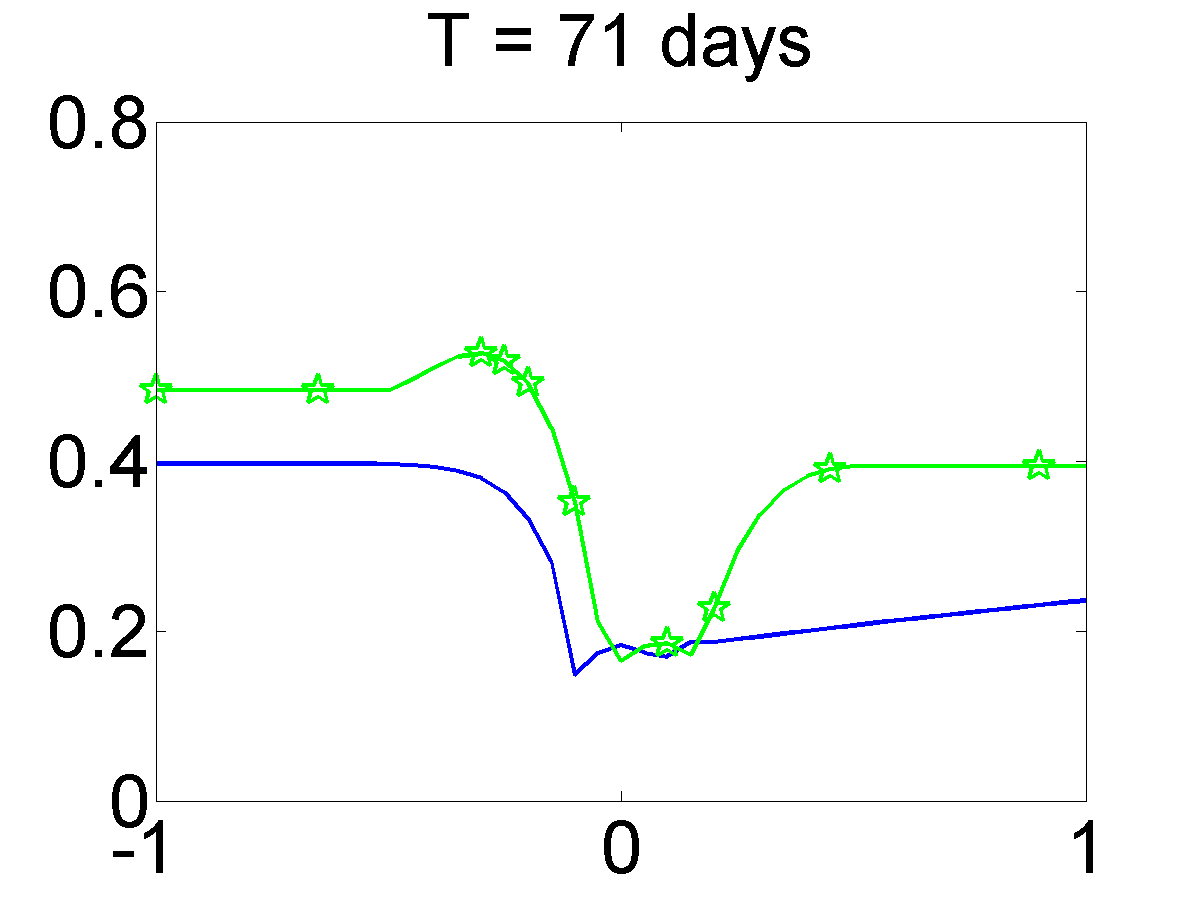}\hfill
       \includegraphics[width=0.475\textwidth]{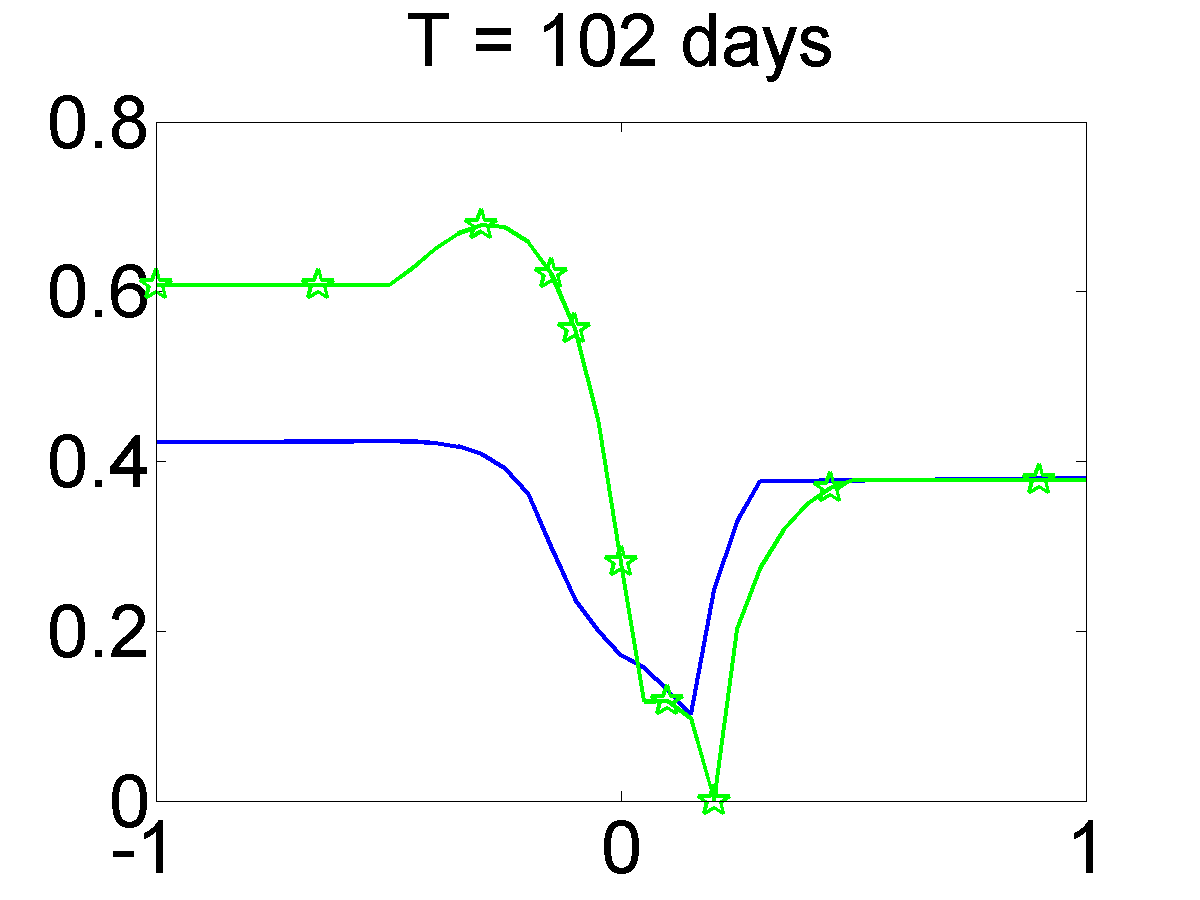}\hfill
       \includegraphics[width=0.475\textwidth]{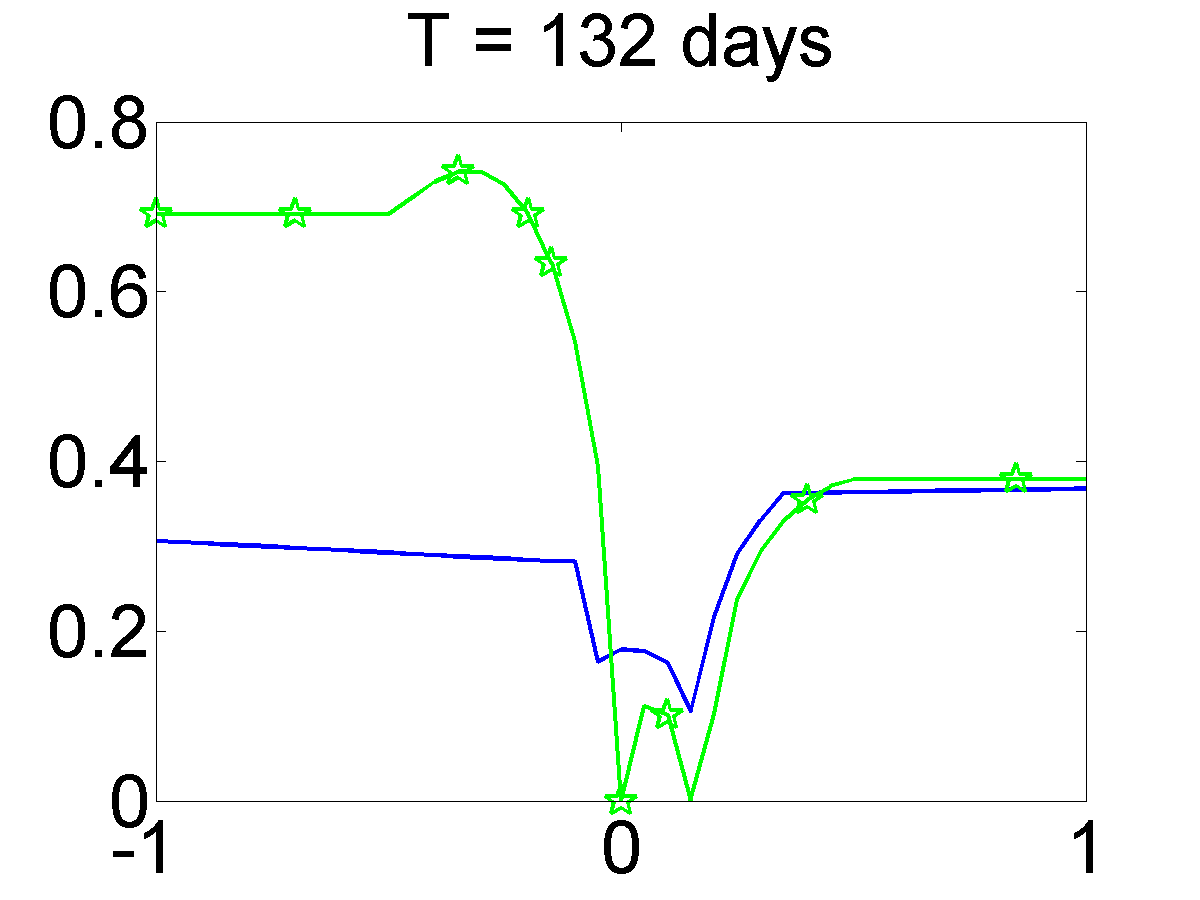}\hfill
       \includegraphics[width=0.475\textwidth]{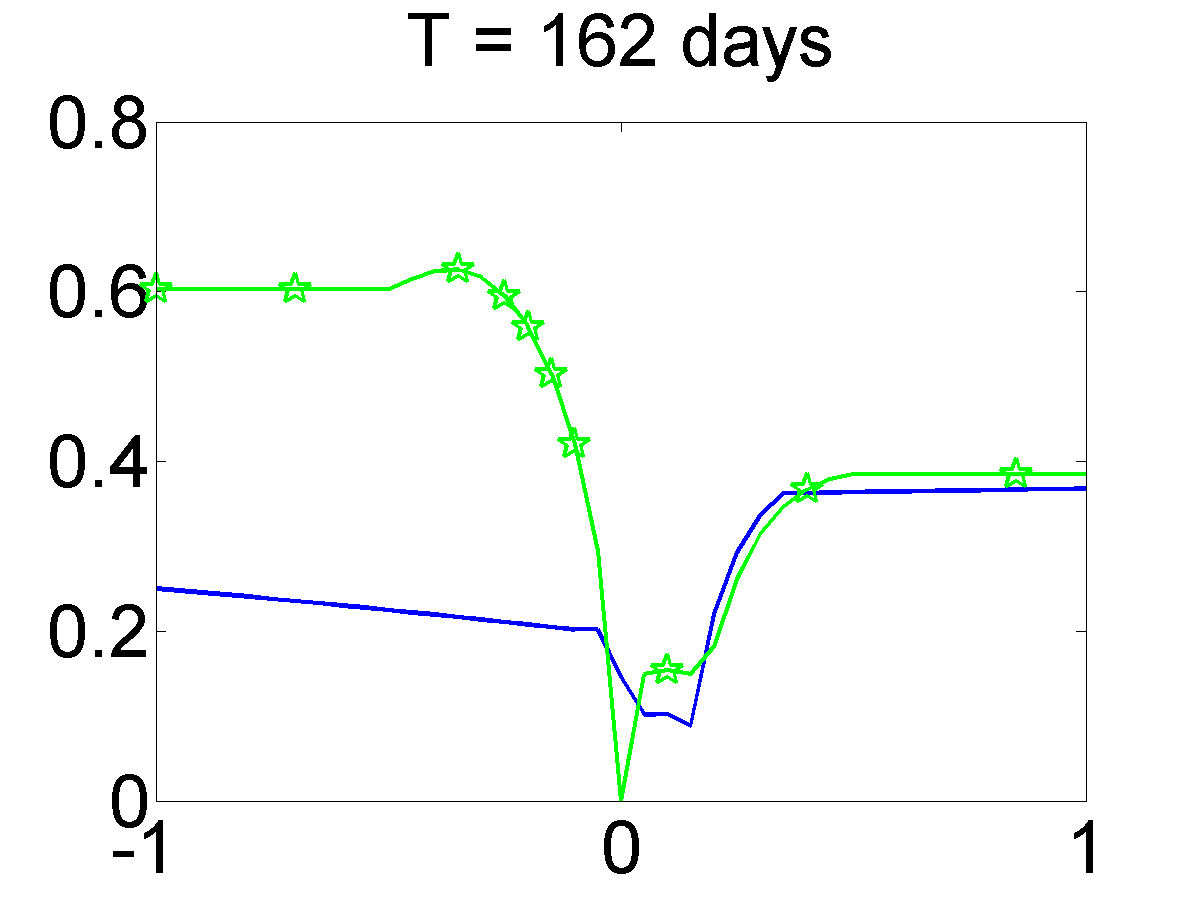}\hfill
       \includegraphics[width=0.475\textwidth]{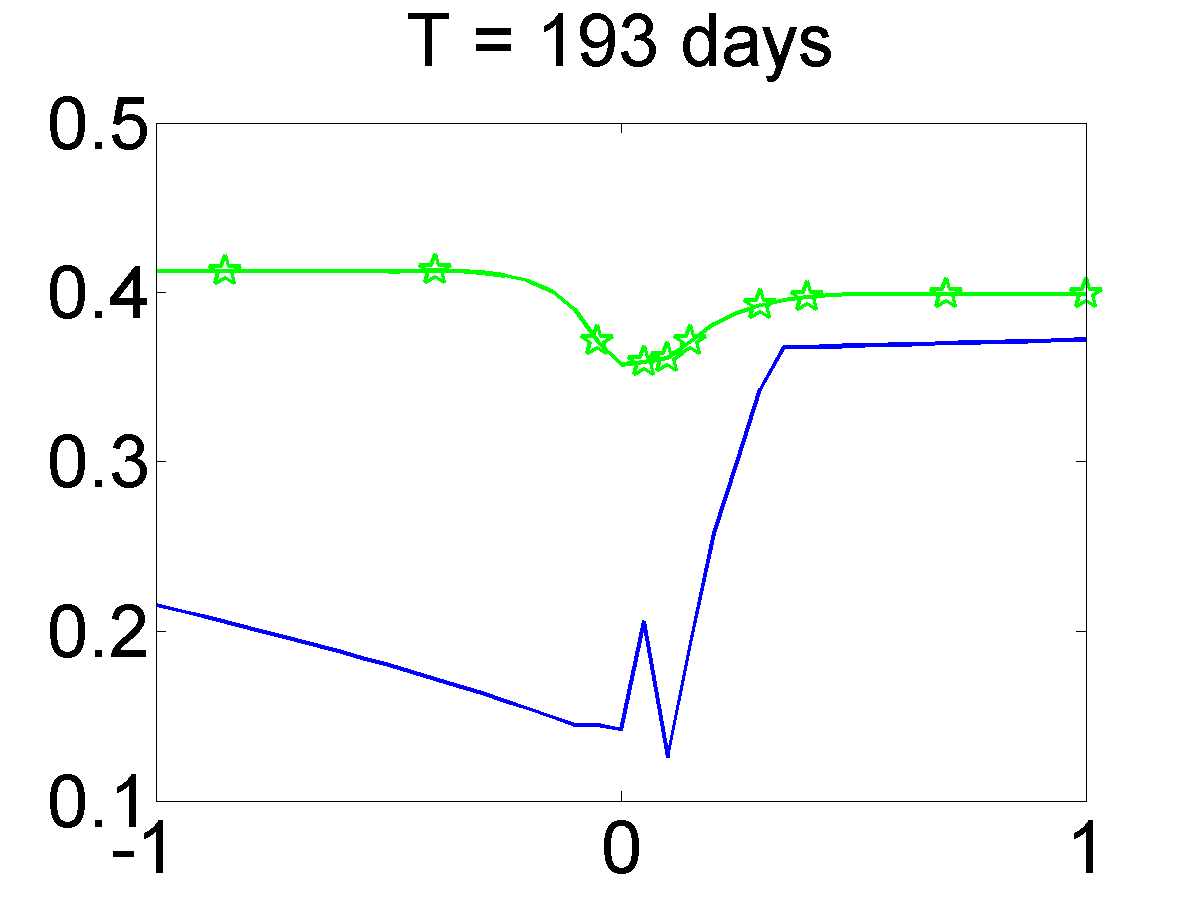}
 \caption{Reconstructed local volatility for different maturity dates for WTI call option prices, 
comparing between completed data (green line with pentagram) and scarce data (blue line) results.}
  \label{fig:com_wti}
\end{minipage}%
\end{figure}

Figures~\ref{fig:com_hh} and \ref{fig:com_wti} display reconstructed local volatility surfaces for the different maturities,
comparing between using the given data and the completed data.

\begin{figure}[!ht]
\centering
\begin{minipage}{0.49\textwidth}
   \centering
       \includegraphics[width=0.475\textwidth]{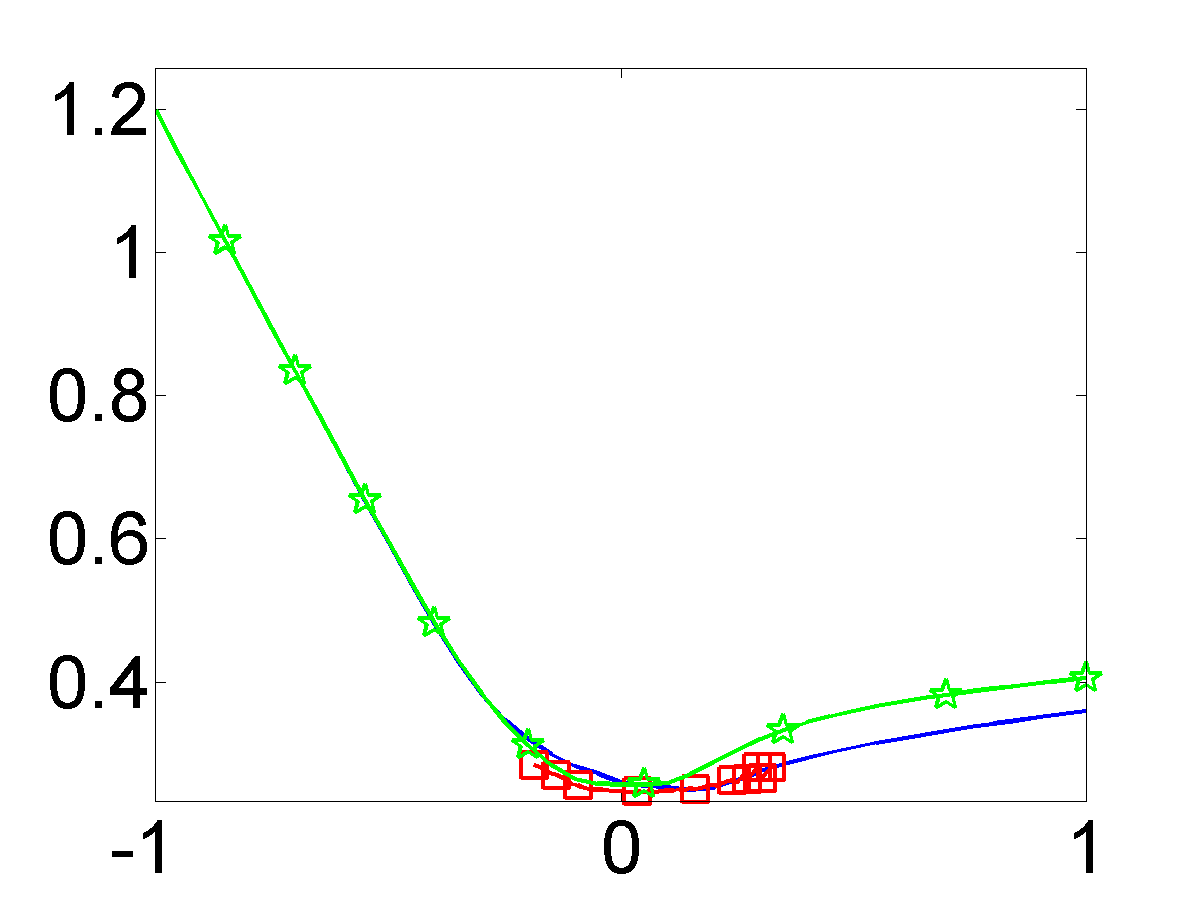}\hfill
       \includegraphics[width=0.475\textwidth]{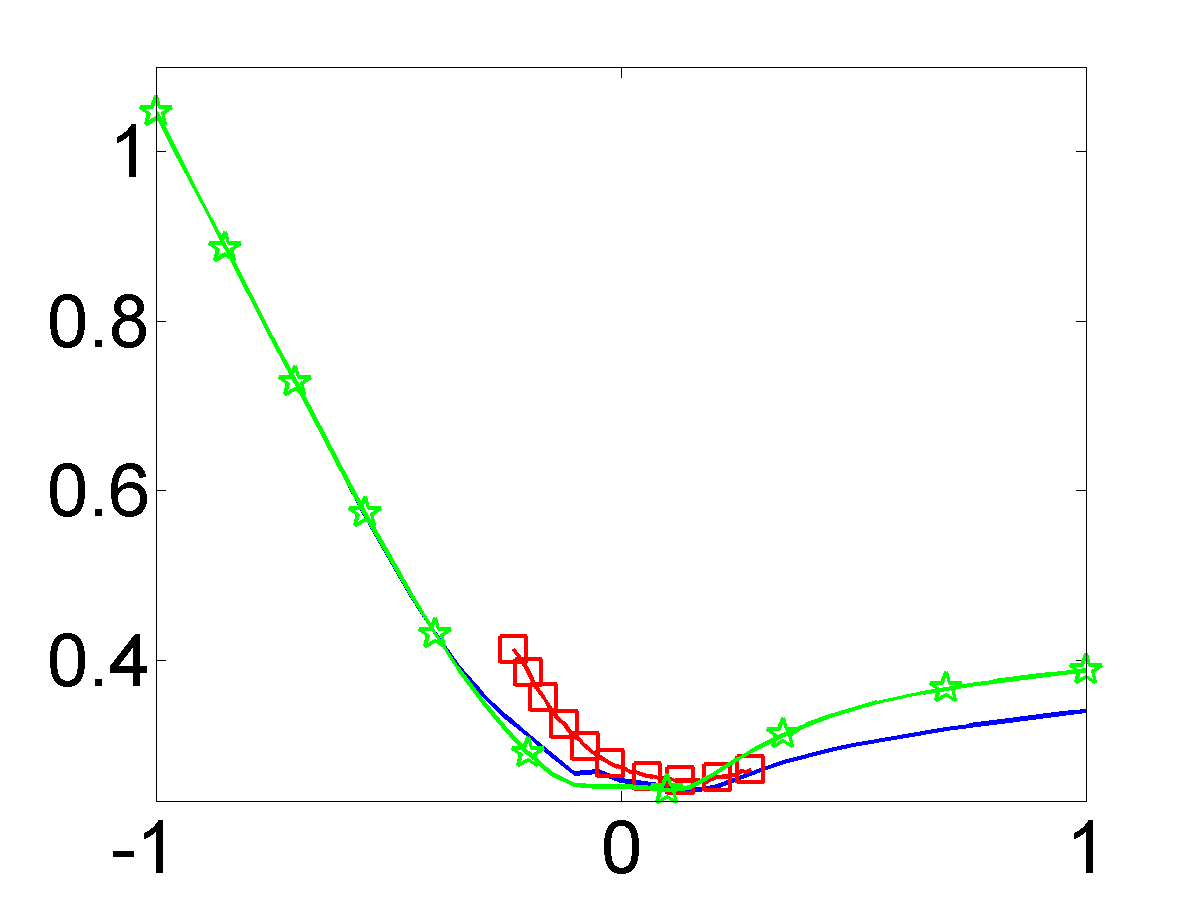}\hfill
       \includegraphics[width=0.475\textwidth]{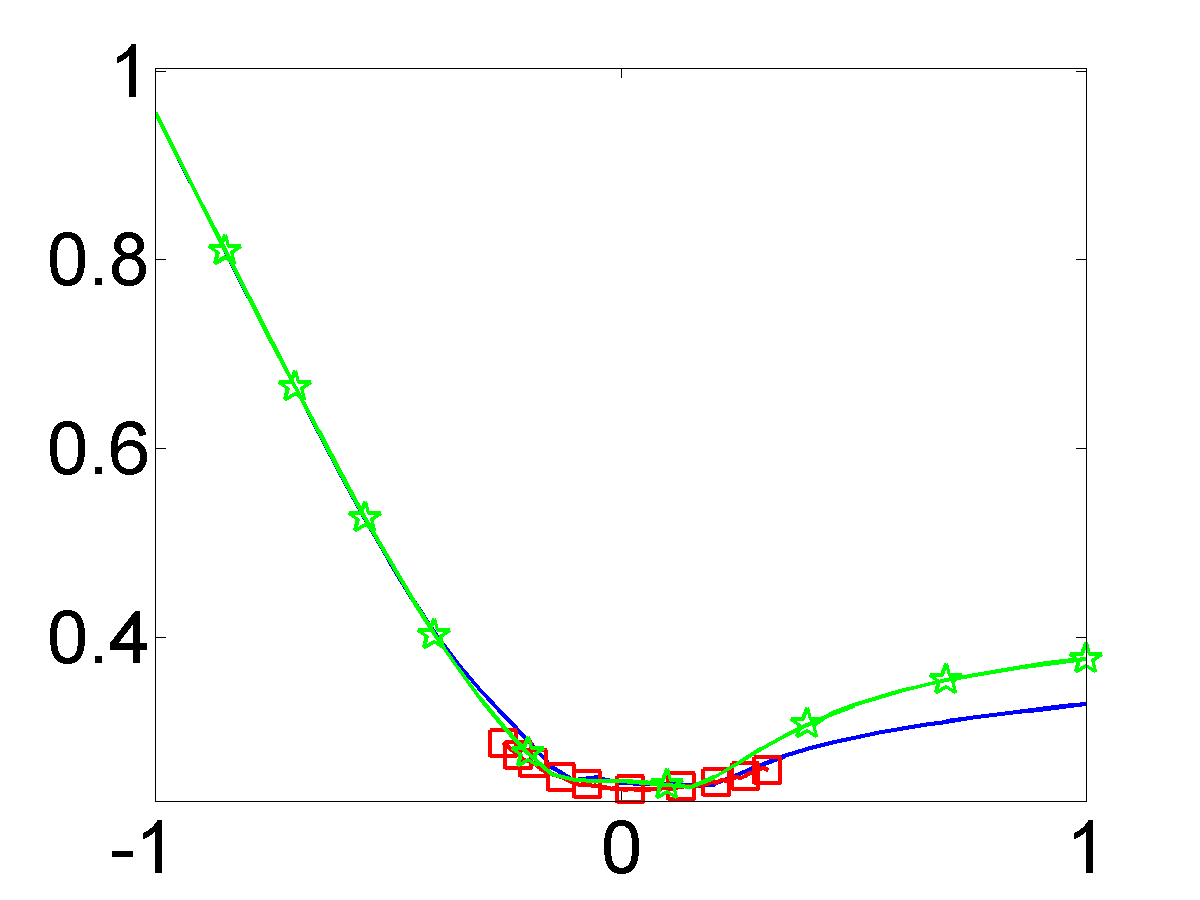}\hfill
       \includegraphics[width=0.475\textwidth]{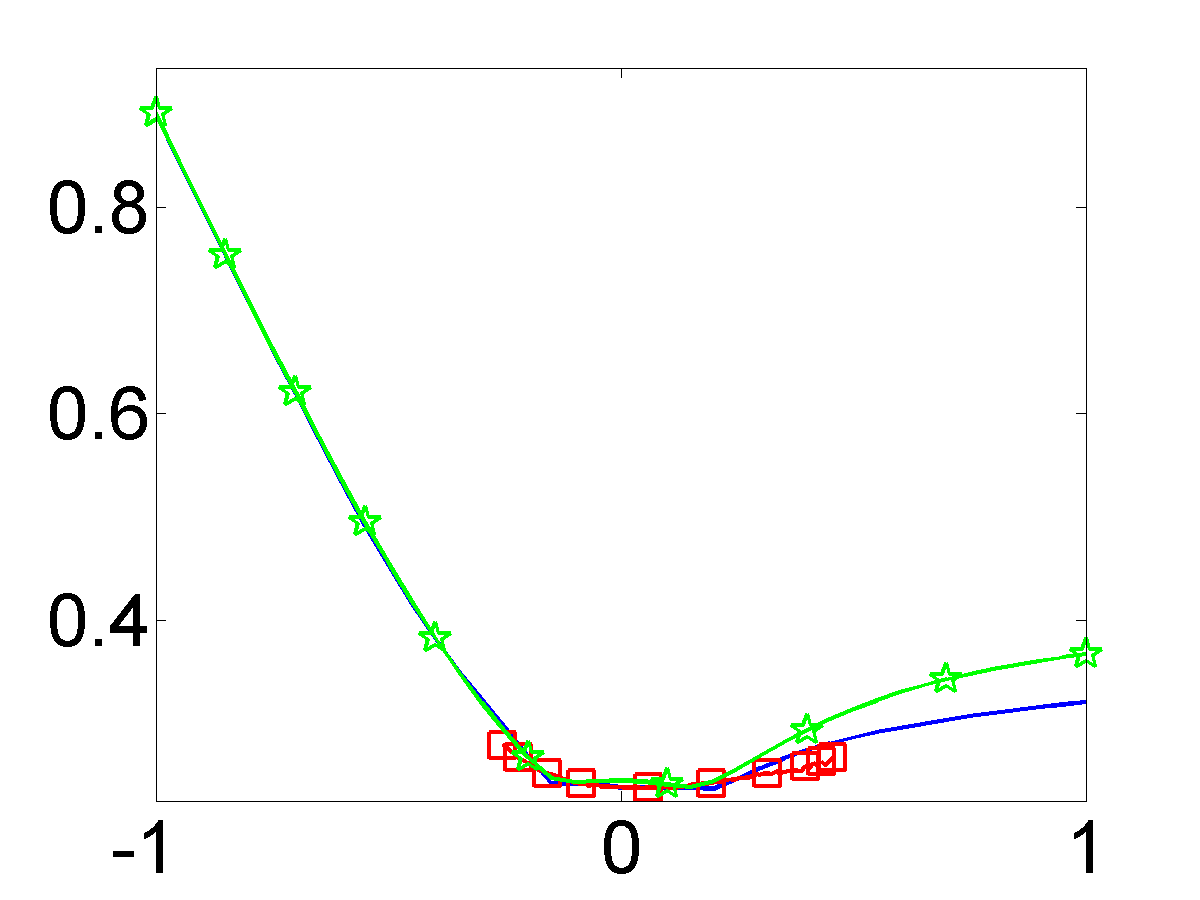}\hfill
       \includegraphics[width=0.475\textwidth]{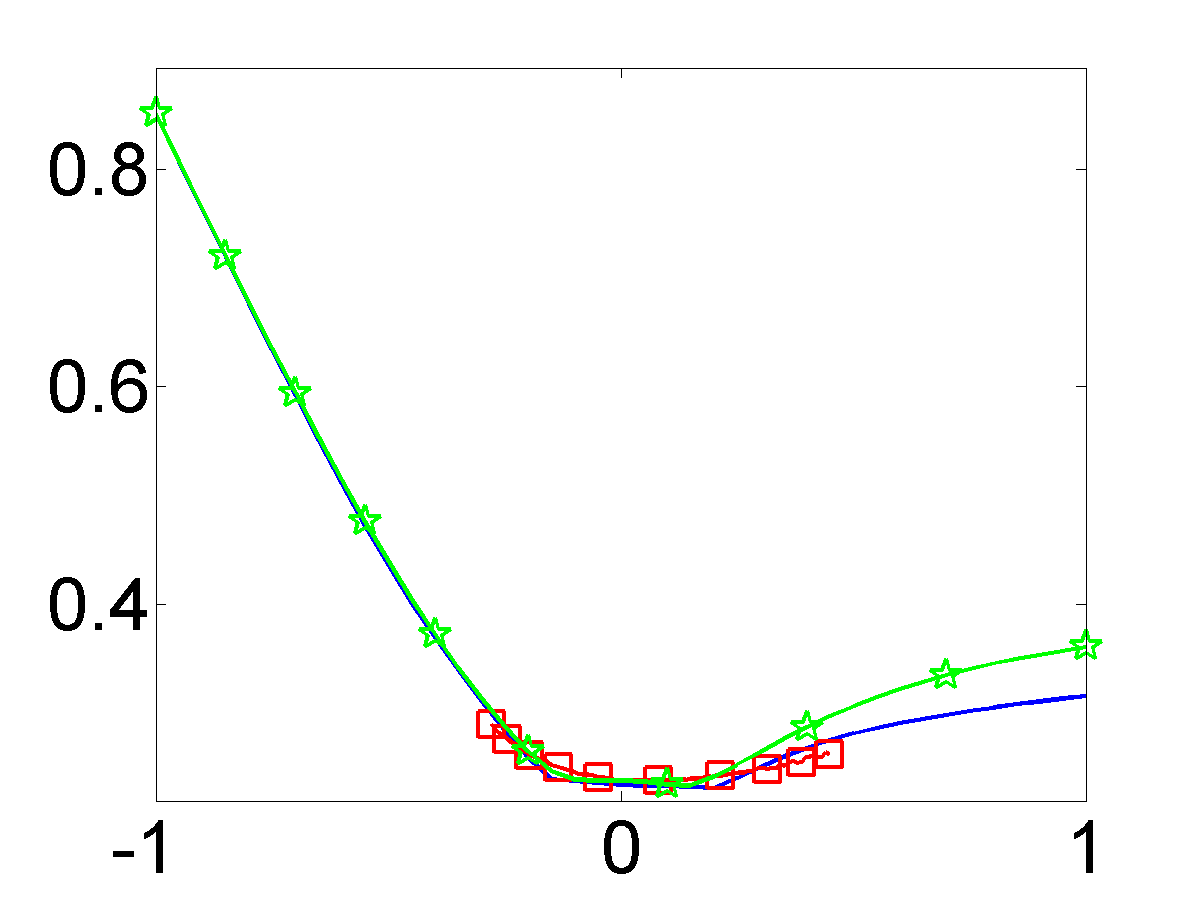}\hfill
       \includegraphics[width=0.475\textwidth]{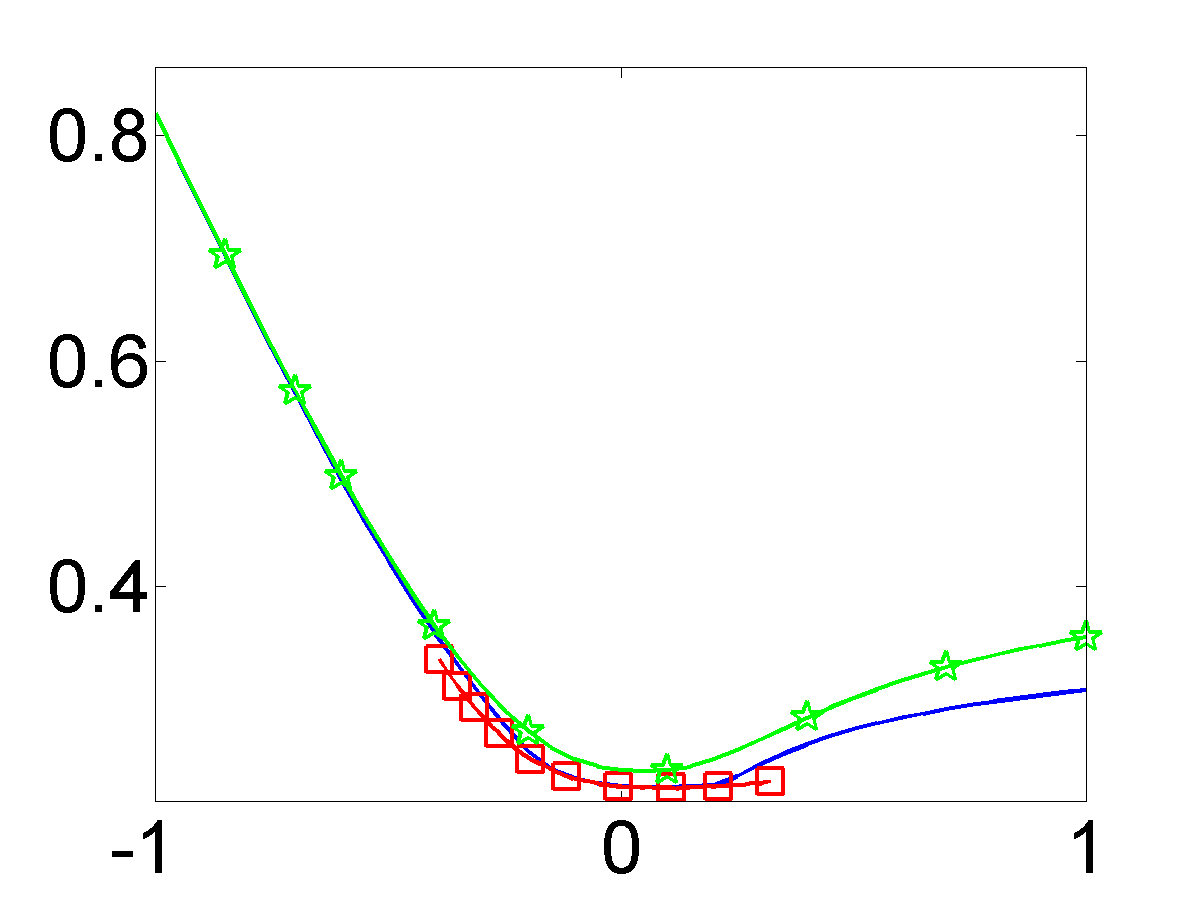}
 \caption{Henry Hub prices: completed data (green line with pentagram), scarce data (blue continuous line), 
and market (red squares) implied volatilities.}
  \label{fig:com2_hh}
\end{minipage}\hfill
\begin{minipage}{0.49\textwidth}
  \centering
   \centering
       \includegraphics[width=0.475\textwidth]{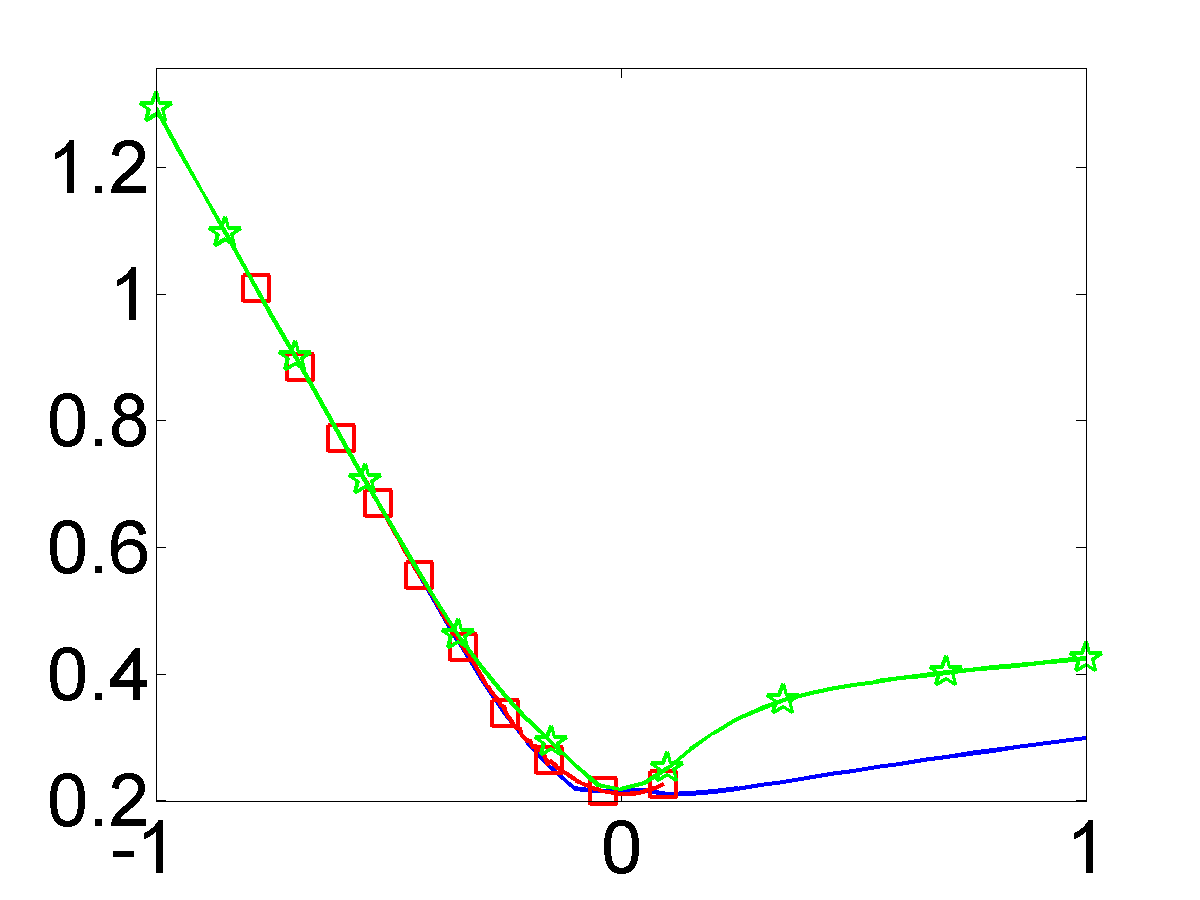}\hfill
       \includegraphics[width=0.475\textwidth]{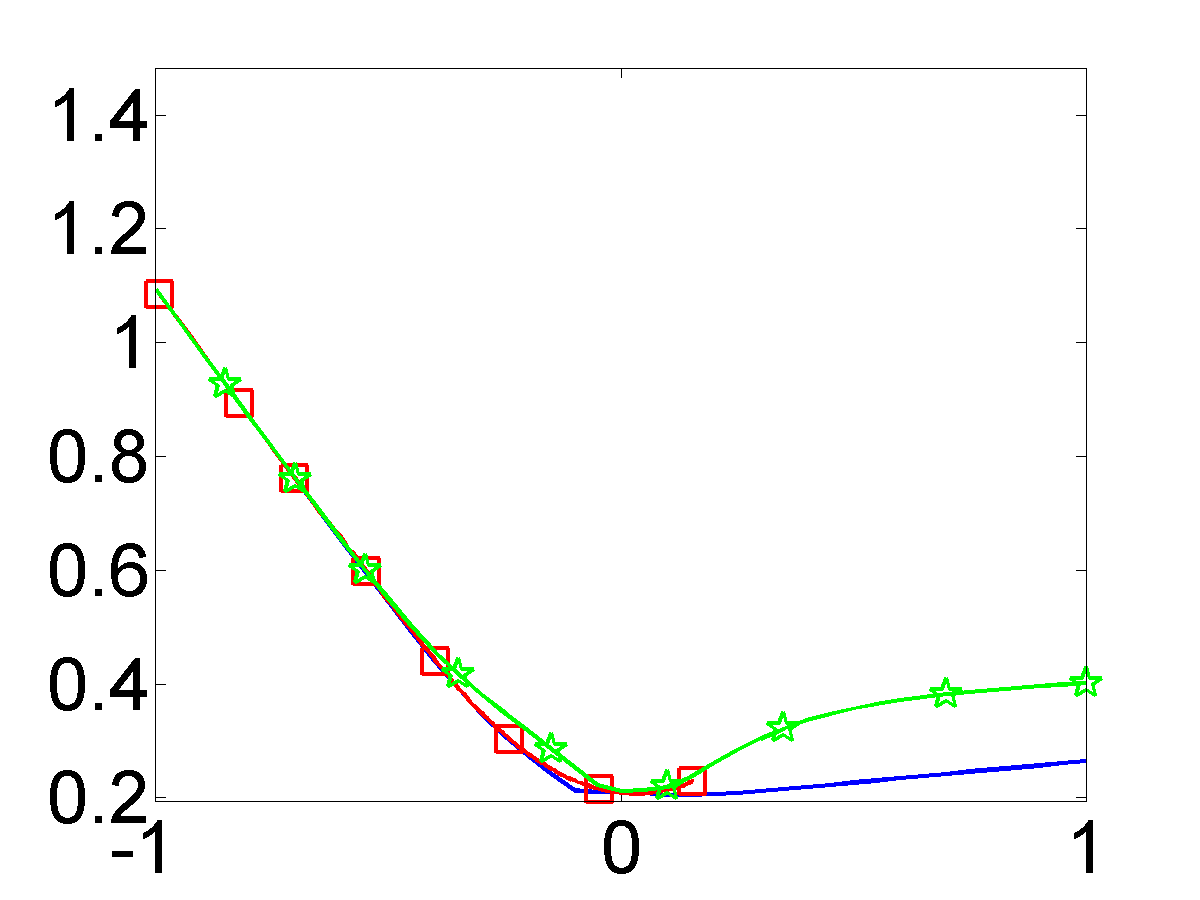}\hfill
       \includegraphics[width=0.475\textwidth]{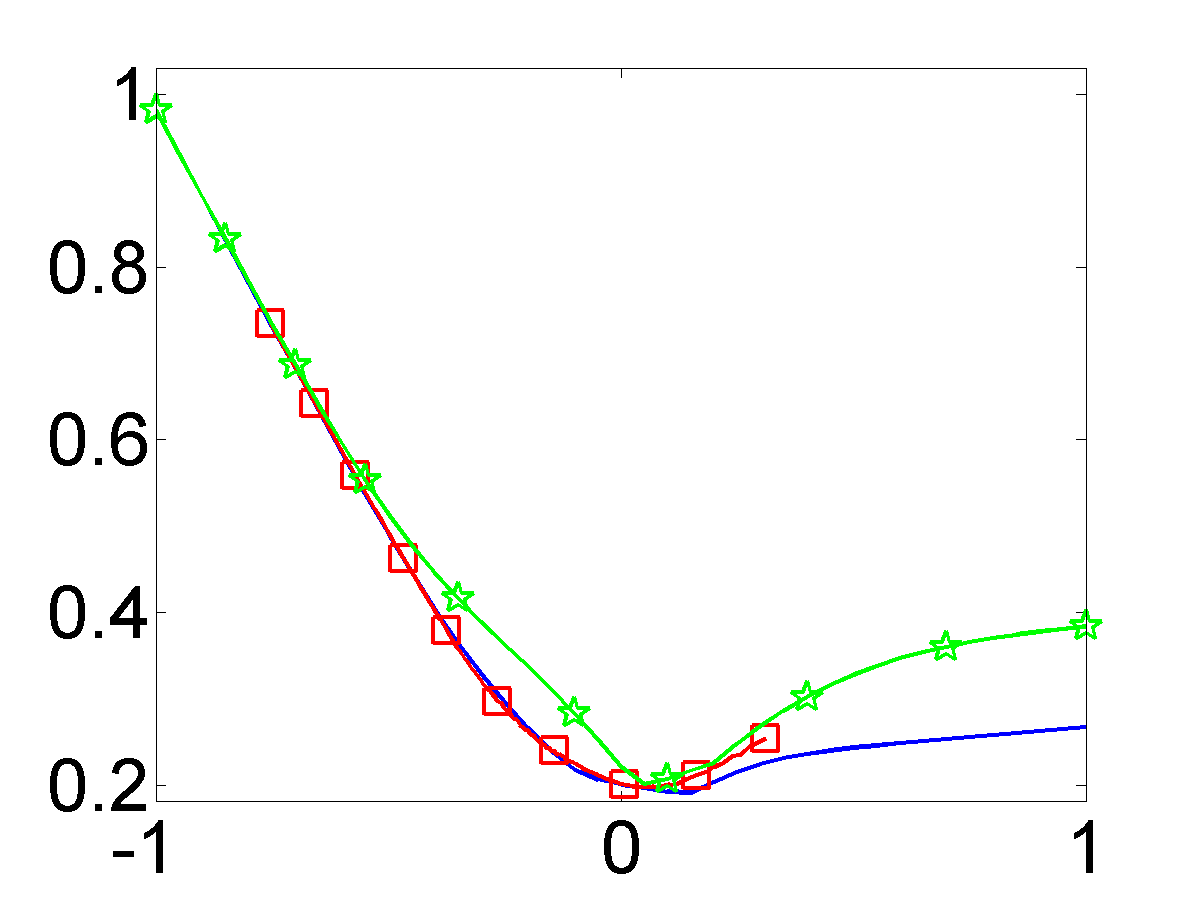}\hfill
       \includegraphics[width=0.475\textwidth]{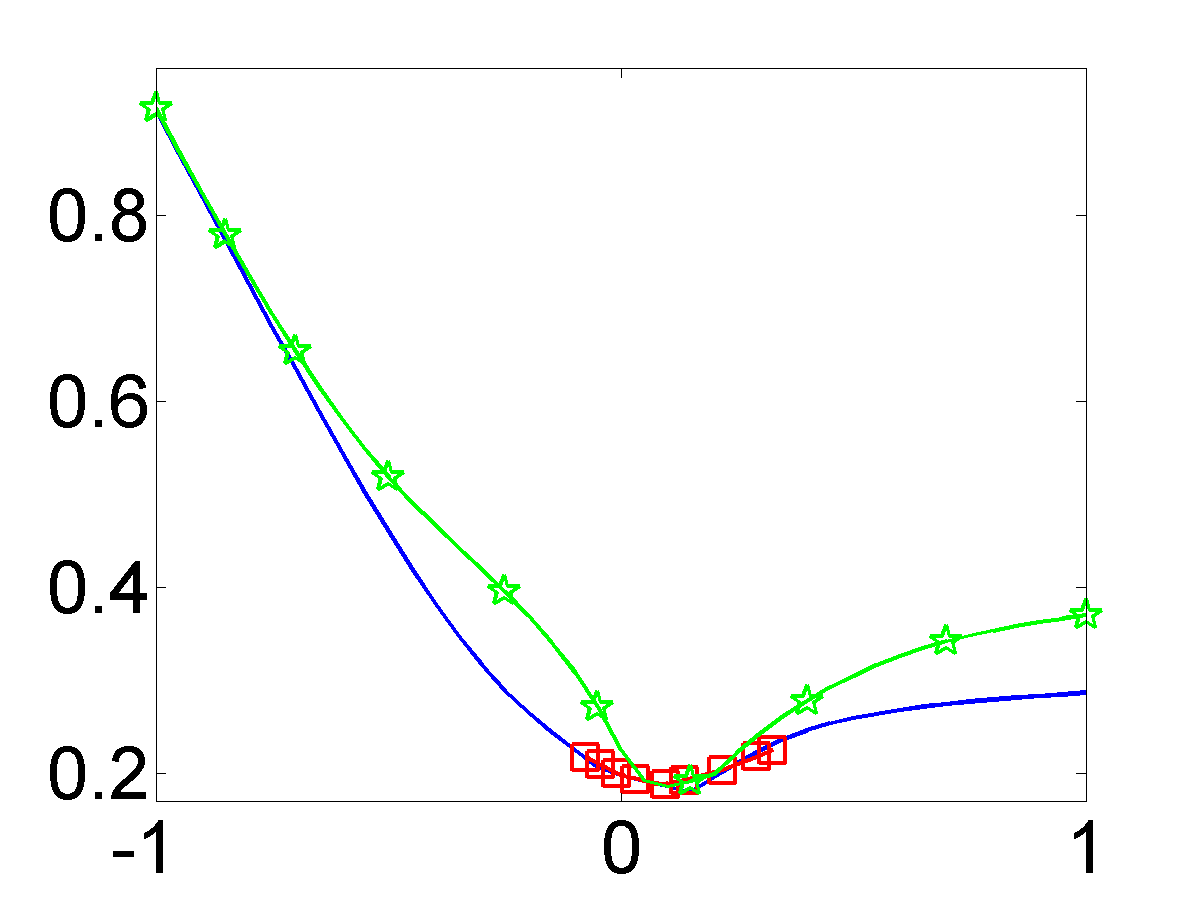}\hfill
       \includegraphics[width=0.475\textwidth]{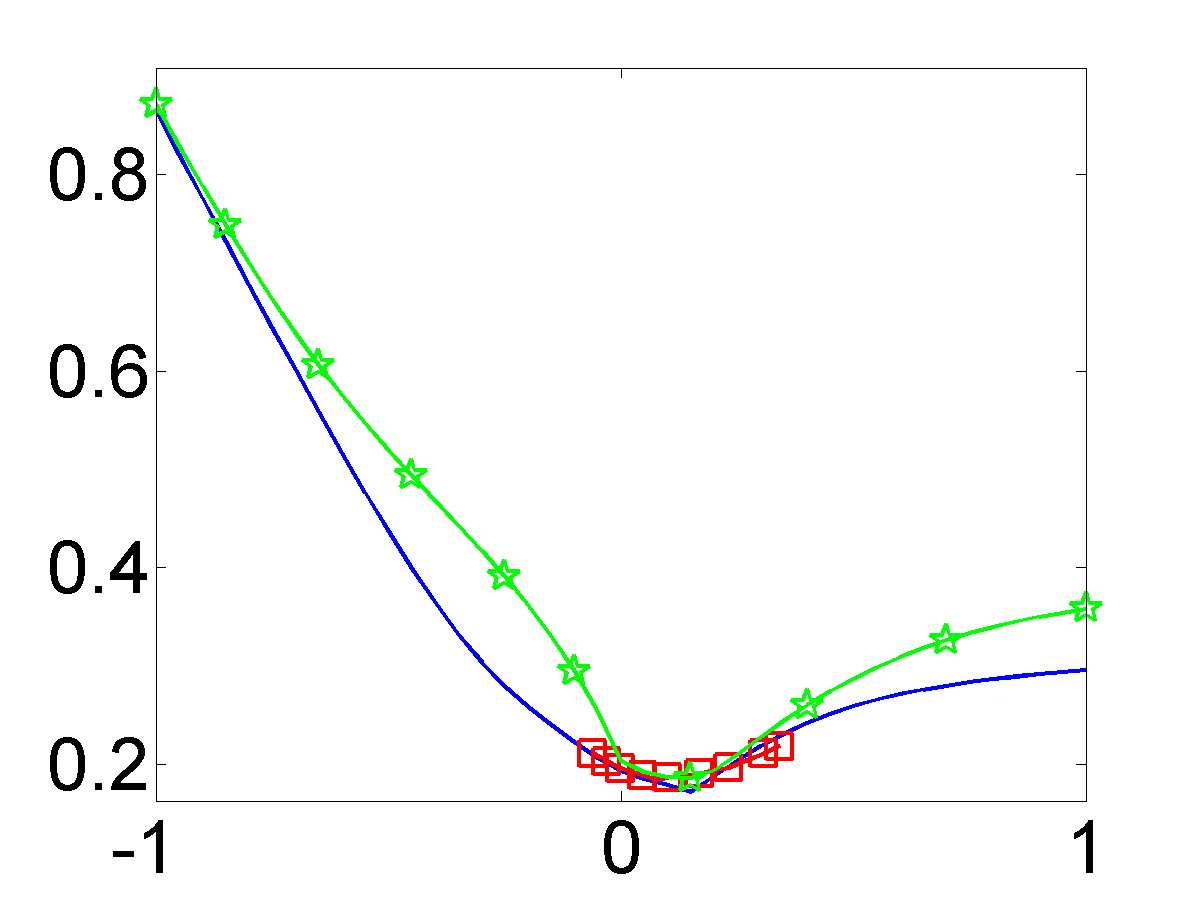}\hfill
       \includegraphics[width=0.475\textwidth]{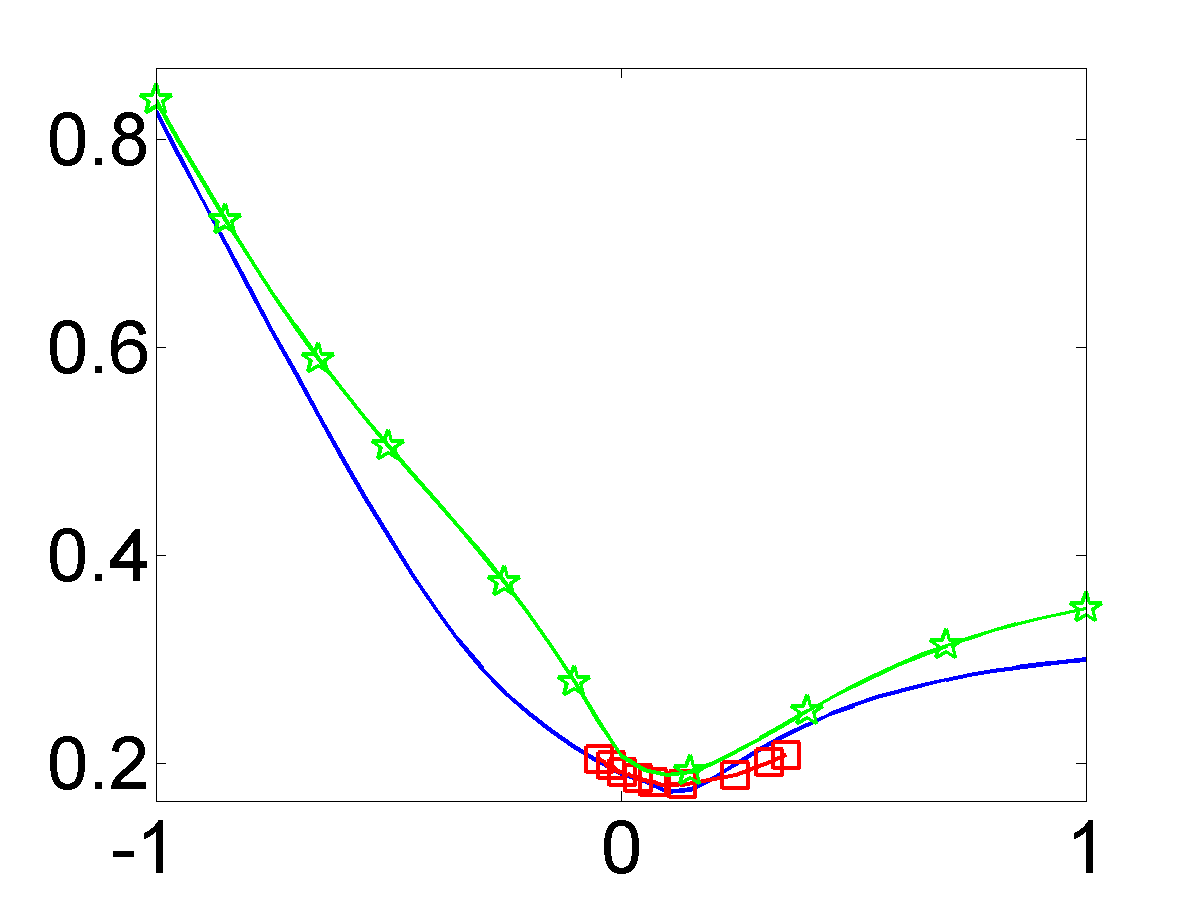}
\caption{WTI prices: completed data (green line with pentagram), sparse data (blue continuous line),
and market (red squares) implied volatilities.}
  \label{fig:com2_wti}
\end{minipage}%
\end{figure}

Figures~\ref{fig:com2_hh} and \ref{fig:com2_wti} present a comparison between the implied volatilities of both methods and the market ones, in order to assess how accurate  the reconstructions are. 
One of the main advantages of the local volatility model is the capability of fitting the market implied smile, 
which has an  important relationship with market risk. The implied volatilities were evaluated using 
the {\sc Matlab} function {\tt blsimpv}, and we used the interest rate as the dividend yield.

In all of these experiments, we have used the mesh widths $\Delta\tau = 1/365$ and $\Delta y = 0.05$,
the annualized risk-free interest rate was taken as $0.25\%$, and $b=0$ in~\eqref{x2}, since futures have no drift. 


Table~\ref{table_dc1} displays the 
parameters obtained in the tests of local volatility calibration with Henry Hub and WTI call prices, 
with scarce and completed data.
In this table, by residual, we mean the $\ell_2$-distance between the evaluated quantity and the data, 
normalized by the $\ell_2$-norm of the data.
\begin{table}[!ht]
\caption{Parameters obtained in the local volatility calibration with Henry Hub and WTI call prices using sparse data and completed data.}
\label{table_dc1}
\begin{tabular}{l|c|c|c|c|}
\cline{2-5}
 &\multicolumn{2}{c|}{WTI}& \multicolumn{2}{c|}{Henry Hub}\\
\cline{2-5}
 & Comp. Data & Sparse Data  & Comp. Data& Sparse Data \\
\hline
\multicolumn{1}{|l|}{$\alpha_0$}          & 1.0e4                 & 1.0e3              & 1.0e3              & 1.0e3 \\
\multicolumn{1}{|l|}{$\alpha_1=\alpha_2=\alpha_3$}& 4.5 & 1.0 & 1.3 & 1.0\\
\multicolumn{1}{|l|}{Price Residual}&  2.16e-2 & 3.21e-3 & 3.47e-2 & 2.14e-2 \\
\multicolumn{1}{|l|}{Implied Vol. Residual}&  1.26e-1 & 2.66e-2 & 9.61e-2 & 5.98e-2 \\
\hline
\end{tabular}
\end{table}

\paragraph{Discussion of the real commodity data results}

Observing the market implied volatilities and the implied volatilities obtained with both methods 
in Figures~\ref{fig:com2_hh} and~\ref{fig:com2_wti}, the results with scarce data present a much better smile adherence 
than when using completed data. So, completing the data can be seen as an unnecessary introduction of noise
or inconsistency. It can be better noticed when observing the implied volatilities at deep in-the-money ($y<-0.1$) 
and deep out-of-the-money ($y>0.1$) strikes. 
In these regions, for almost all cases, the results with scarce data practically matched the implied volatility, 
whereas with data completion, the resulting implied volatilities presented higher values. 
For financial market practitioners, higher implied volatilities can be translated to higher risk. 
So, using data completion could lead investors to be more conservative than necessary.


\section{Conclusions}
\label{sec:conclusions}

The questions of how to treat observed data in order to produce agreeable solutions,
and relatedly, how much to trust the quality of a given data set collected by another agent,
are prevalent in many nontrivial applied inverse problems. 
Standard theory appears to have little to contribute to their satisfactory resolution,
and more carefully assembled experience is required.
In this article we have highlighted some of the issues involved through an important
application in mathematical finance, and we have proposed methods that improve on
techniques available in the open literature.

The problem of constructing a local volatility surface has similarities with some
applications in areas such as geophysics and medical imaging, in that it boils down
to the calibration of a diffusive PDE problem, reconstructing a distributed parameter function,
i.e., a surface, rather than a few unrelated parameter values.
The difficulty of dealing with scarce data, which highlights the need for a careful
practical selection of a prior, is common, as is the uncertainty in data location,
although the latter appears here as dependence only on a scalar additional variable.

The finance problem considered here is distinguished from its more physical comrades in two important
aspects. One is the availability of real data: we have experimented here with three different
real applications, while most papers appearing in the applied mathematics literature
never get to deal with real data at all. The other aspect is the difficulty of assessing the resulting
recovered volatilities: there is no physical solution here, data sets are changing daily,
and experience rules.
In the present setting we had to determine the coefficients of the regularization operators,
for instance $\alpha_i$ in \eqref{x5}, \eqref{x7} and \eqref{S0var}, by trial and error.

Of particular interest is the regularization term involving $a_0$ and $\alpha_1$.
This term is a penalty on the sought local volatility surface
for straying away from a given constant, $a_0$, which in turn is estimated based
on the type of asset under consideration.
The question whether or not this should be done is a deep one
and touches upon the very foundations of the
model under consideration. It also has a practical unfolding, since
the possibility that the volatility, as a function of the asset price,
grows at a sufficiently fast rate may be connected (at least in some
similar models) to the presence of market bubbles~\cite{jakcpr}.

The EnKF method considered in Section~\ref{sec:enkf} is an adaptation of one of several
methods proposed in the literature~\cite{iglast,reco}, improved by adding smoothing penalty terms.
We have also experimented with a similar adaptation of the method in~\cite{caerso}.
In both cases the results are not consistently better than those obtained by the Tikhonov-type method
of Section~\ref{sec:scarce}, which in hindsight is not surprising. However, since such methods
are currently very much in vogue (especially in our target application area) we feel that our
results in this sense are important: one must have data of sufficient quality for such methods
to dominate. 




\section*{Acknowledgements}
VA acknowledges and thanks CNPq through grant 201644/2014-2.
UMA and XY acknowledge with thanks a {\em Ciencias Sem Fronteiras} (visiting scientist / postdoc) grant from CAPES, Brazil.
JPZ was supported by the Brazilian National Research Council (CNPq) under grant 307873.

\appendix
\section{EnKF details}
\label{sec:enkf_details}


The regularized weighted least-squares problem (\ref{enkf4}) has the solution 
\begin{eqnarray}
\tilde{a}_h&=&(D^{-1}+
H^T\Gamma^{-1}H+L_{\tau}^TD_{\tau}^{-1}
L_{\tau}+
L_y^TD_y^{-1}L_y)^{-1} \cdot\nonumber\\
& &(D^{-1}\hat{a}_0+H^T\Gamma^{-1}
d+L_{\tau}^TD_{\tau}^{-1}L_{\tau}\hat{a}_0
+L_y^TD_y^{-1}L_y\hat{a}_0).
\label{enkf3}
\end{eqnarray}
However, 
$D$ is unknown. 
To address this we apply EnKF and replace $D$ by the covariance matrix computed 
from generated samples. 
We also replace $L_{\tau}\hat{a}_0$ and $L_y\hat{a}_0$ in~\eqref{enkf4} by 
$r_{\tau}$ and $r_y$, sampled from $\mathcal{N}(L_{\tau}\hat{a}_0,D_{\tau})$ 
and $\mathcal{N}(L_y\hat{a}_0, D_y)$, respectively. 
In this way, $r_{\tau}$ and $r_y$ can be viewed as two observations like $d$. 
Next we explain the calculations in one iteration of the prediction and analysis steps
that appear in the EnKF algorithm stated at the end of Section~\ref{sec:enkf}. 

Having generated $\{a_h^{(0,j)}\}_{j=1}^J$ as the initial ensemble,
we calculate $Pu_h(a_h^{(0,j)})$ and define
$$\hat{a}_h^{(1,j)}=\begin{pmatrix}
a_h^{(0,j)} \\
Pu_h(a_h^{(0,j)})
\end{pmatrix}, \quad j=1,2,\ldots, J . $$
This is the first prediction step (i.e., we set $n = 0$ in the stated algorithm). 
The sample mean and covariance matrix are then given by
\begin{eqnarray*}
D_1=\frac{1}{J}\sum_{j=1}^J\hat{a}_h^{(1,j)}(\hat{a}_h^{(1,j)})^T - 
\bar{a}_h^{(1)}(\bar{a}_h^{(1)})^T, \quad \bar{a}_h^{(1)}=\frac{1}{J}\sum_{j=1}^J \hat{a}_h^{(1,j)}.\label{enkf5}
\end{eqnarray*}

We now calculate the Kalman gain three times to obtain 
$$
(D_1^{-1}+
H^T\Gamma^{-1}H+L_{\tau}^TD_{\tau}^{-1}
L_{\tau}+
L_y^TD_y^{-1}L_y)^{-1}.$$
Denoting 
$$U=(D_1^{-1}+
H^T\Gamma^{-1}H+L_{\tau}^TD_{\tau}^{-1}
L_{\tau})^{-1} ,$$
we have 
\begin{eqnarray*}
&&(D_1^{-1}+
H^T\Gamma^{-1}H+L_{\tau}^TD_{\tau}^{-1}
L_{\tau}+
L_y^TD_y^{-1}L_y)^{-1} \\
&=&(U^{-1}+L_y^TD_y^{-1}L_y)^{-1} \\
&=&U-UL_y^T(L_yUL_y^T+D_y)^{-1}L_yU \\
&=&(I-W_1L_y)U ,
\end{eqnarray*}
where $W_1=UL_y^T(L_yUL_y^T+D_y)^{-1}$.
We continue with the same procedure to calculate $U$. Let 
$V^{-1}=D_1^{-1}+H^T\Gamma^{-1}H$ and $W_2=VL_{\tau}^T(L_{\tau}VL_{\tau}^T+D_{\tau})^{-1}$.
Then we have 
\begin{eqnarray*}
U&=&(D_1^{-1}+
H^T\Gamma^{-1}H+L_{\tau}^TD_{\tau}^{-1}
L_{\tau})^{-1} \\
&=&(V^{-1}+L_{\tau}^TD_{\tau}^{-1}
L_{\tau})^{-1} \\
&=&V-VL_{\tau}^T(L_{\tau}VL_{\tau}^T+
D_{\tau})^{-1}L_{\tau}V \\
&=&(I-W_2L_{\tau})V .
\end{eqnarray*}
Hence, defining also $W_3=D_1H^T(HD_1H^T+\Gamma)^{-1}$ we have
\begin{eqnarray*}
V&=&(D_1^{-1}+H^T\Gamma^{-1}H)^{-1} \\
&=&D_1-D_1H^T(HD_1^T+\Gamma)^{-1}HD_1 \\
&=&(I-W_3H)D_1 .
\end{eqnarray*}

Now \eqref{enkf3} reads
\begin{eqnarray*}
\tilde{a}_h^{(1,j)}&=&U \left(D^{-1}\hat{a}_h^{(1,j)}+H^T\Gamma^{-1}
d+L_{\tau}^TD_{\tau}^{-1}L_{\tau}\hat{a}_h^{(1,j)}
+L_y^TD_y^{-1}L_y\hat{a}_h^{(1,j)}\right) \\
&=&U \left(U^{-1}\hat{a}_h^{(1,j)}+H^T\Gamma^{-1}(d-H\hat{a}_h^{(1,j)})+
L_{\tau}^TD_{\tau}^{-1}(r_{\tau}-L_{\tau}
\hat{a}_h^{(1,j)})+L_y^TD_y^{-1}(r_y-L_y\hat{a}^{(1,j)}_h)\right) \\
&=&\hat{a}_h^{(1,j)}+U \left(H^T\Gamma^{-1}(d-H\hat{a}_h^{(1,j)})+L_{\tau}^TD_{\tau}^{-1}(r_{\tau}-L_{\tau}
\hat{a}_h^{(1,j)})+L_y^TD_y^{-1}(r_y-L_y\hat{a}^{(1,j)}_h)\right) .
\end{eqnarray*} 
Reversing this procedure and adjusting notation, we obtain the analysis step 2(b)
of the algorithm in Section~\ref{sec:enkf} for $n = 0$.

In each iteration of the EnKF algorithm, 
the prediction step is used to map the samples into the data space. 
The analysis step then calculates the ``distance'' between the mapped ensemble and the noisy data.
Following the analysis step, the approximated local variance and estimated option prices are calculated, 
and this is used to compare with the original data. 
We compute the residual and compare to the one from the last iteration: 
Step 2 of the algorithm is repeated until the residual is  less than a certain threshold.

%


\newcommand{\etalchar}[1]{$^{#1}$}

\end{document}